\documentclass[preprint,3p,sort&compress]{elsarticle}
\journal{{\tt Journal of Computational Physics}}

\usepackage[dvips]{epsfig}
\usepackage{graphicx} 		
\usepackage{pdfpages}
\usepackage{latexsym}
\usepackage{verbatim}
\usepackage{amsmath,amsthm}
\usepackage[utf8]{inputenc}
\usepackage{amssymb}
\usepackage{bbm}
\usepackage{fonts}
\usepackage{enumerate}
\usepackage{placeins}
\usepackage{standalone}
\usepackage{paralist}
\usepackage{subcaption}
\usepackage{booktabs}
\usepackage{array}
\newcolumntype{C}[1]{>{\centering\let\newline\\\arraybackslash\hspace{0pt}}m{#1}}

\usepackage{algorithm,algorithmicx,algpseudocode}

\usepackage[]{hyperref} 

\usepackage{pgfplots}
\pgfplotsset{compat=newest}       
\usepackage[]{externaltikz}

 \usepackage{relinput}
\newtheorem{theorem}{Theorem}[section]
\newtheorem{definition}[theorem]{Definition}
\newtheorem{remark}[theorem]{Remark}

\newcounter{tikzsubfigcounter}[figure]
\renewcommand{\thetikzsubfigcounter}{\the\numexpr\value{figure}+1\relax\alph{tikzsubfigcounter}}

\newcounter{tikzsubfigcounterinvisible}[figure]
\renewcommand{\thetikzsubfigcounterinvisible}{\the\numexpr\value{figure}+1\relax\alph{tikzsubfigcounterinvisible}}

\newcommand{\refone}[1]{\textcolor{black}{#1}}
\newcommand{\reftwo}[1]{\textcolor{black}{#1}}
\newcommand{\refthree}[1]{\textcolor{black}{#1}}
\newcommand{\reffour}[1]{\textcolor{black}{#1}}

\numberwithin{equation}{section}

\title{\scheme{} approximation for hyperbolic conservation laws}

\author[ls]{Louisa Schlachter}
\address[ls]{Fachbereich Mathematik, TU Kaiserslautern, Erwin-Schr\"odinger-Str., 67663 Kaiserslautern, Germany, {\tt schlacht@mathematik.uni-kl.de}}

\author[fs]{Florian Schneider}
\address[fs]{Fachbereich Mathematik, TU Kaiserslautern, Erwin-Schr\"odinger-Str., 67663 Kaiserslautern, Germany, {\tt schneider@mathematik.uni-kl.de}}

\author[ok]{Oliver Kolb}
\address[ok]{Institut für Mathematik, Universität Mannheim, A5, 68131 Mannheim, Germany, {\tt kolb@uni-mannheim.de}}

\date{}

\definecolor{greenyellow}   {cmyk}{0.15, 0   , 0.69, 0   }
\definecolor{yellow}        {cmyk}{0   , 0   , 1   , 0   }
\definecolor{goldenrod}     {cmyk}{0   , 0.10, 0.84, 0   }
\definecolor{dandelion}     {cmyk}{0   , 0.29, 0.84, 0   }
\definecolor{apricot}       {cmyk}{0   , 0.32, 0.52, 0   }
\definecolor{peach}         {cmyk}{0   , 0.50, 0.70, 0   }
\definecolor{melon}         {cmyk}{0   , 0.46, 0.50, 0   }
\definecolor{yelloworange}  {cmyk}{0   , 0.42, 1   , 0   }
\definecolor{orange}        {cmyk}{0   , 0.61, 0.87, 0   }
\definecolor{burntorange}   {cmyk}{0   , 0.51, 1   , 0   }
\definecolor{bittersweet}   {cmyk}{0   , 0.75, 1   , 0.24}
\definecolor{redorange}     {cmyk}{0   , 0.77, 0.87, 0   }
\definecolor{mahogany}      {cmyk}{0   , 0.85, 0.87, 0.35}
\definecolor{maroon}        {cmyk}{0   , 0.87, 0.68, 0.32}
\definecolor{brickred}      {cmyk}{0   , 0.89, 0.94, 0.28}
\definecolor{red}           {cmyk}{0   , 1   , 1   , 0   }
\definecolor{orangered}     {cmyk}{0   , 1   , 0.50, 0   }
\definecolor{rubinered}     {cmyk}{0   , 1   , 0.13, 0   }
\definecolor{wildstrawberry}{cmyk}{0   , 0.96, 0.39, 0   }
\definecolor{salmon}        {cmyk}{0   , 0.53, 0.38, 0   }
\definecolor{carnationpink} {cmyk}{0   , 0.63, 0   , 0   }
\definecolor{magenta}       {cmyk}{0   , 1   , 0   , 0   }
\definecolor{violetred}     {cmyk}{0   , 0.81, 0   , 0   }
\definecolor{rhodamine}     {cmyk}{0   , 0.82, 0   , 0   }
\definecolor{mulberry}      {cmyk}{0.34, 0.90, 0   , 0.02}
\definecolor{redviolet}     {cmyk}{0.07, 0.90, 0   , 0.34}
\definecolor{fuchsia}       {cmyk}{0.47, 0.91, 0   , 0.08}
\definecolor{lavender}      {cmyk}{0   , 0.48, 0   , 0   }
\definecolor{thistle}       {cmyk}{0.12, 0.59, 0   , 0   }
\definecolor{orchid}        {cmyk}{0.32, 0.64, 0   , 0   }
\definecolor{darkorchid}    {cmyk}{0.40, 0.80, 0.20, 0   }
\definecolor{purple}        {cmyk}{0.45, 0.86, 0   , 0   }
\definecolor{plum}          {cmyk}{0.50, 1   , 0   , 0   }
\definecolor{violet}        {cmyk}{0.79, 0.88, 0   , 0   }
\definecolor{royalpurple}   {cmyk}{0.75, 0.90, 0   , 0   }
\definecolor{blueviolet}    {cmyk}{0.86, 0.91, 0   , 0.04}
\definecolor{periwinkle}    {cmyk}{0.57, 0.55, 0   , 0   }
\definecolor{cadetblue}     {cmyk}{0.62, 0.57, 0.23, 0   }
\definecolor{cornflowerblue}{cmyk}{0.65, 0.13, 0   , 0   }
\definecolor{midnightblue}  {cmyk}{0.98, 0.13, 0   , 0.43}
\definecolor{navyblue}      {cmyk}{0.94, 0.54, 0   , 0   }
\definecolor{royalblue}     {cmyk}{1   , 0.50, 0   , 0   }
\definecolor{blue}          {cmyk}{1   , 1   , 0   , 0   }
\definecolor{cerulean}      {cmyk}{0.94, 0.11, 0   , 0   }
\definecolor{cyan}          {cmyk}{1   , 0   , 0   , 0   }
\definecolor{processblue}   {cmyk}{0.96, 0   , 0   , 0   }
\definecolor{skyblue}       {cmyk}{0.62, 0   , 0.12, 0   }
\definecolor{turquoise}     {cmyk}{0.85, 0   , 0.20, 0   }
\definecolor{tealblue}      {cmyk}{0.86, 0   , 0.34, 0.02}
\definecolor{aquamarine}    {cmyk}{0.82, 0   , 0.30, 0   }
\definecolor{bluegreen}     {cmyk}{0.85, 0   , 0.33, 0   }
\definecolor{emerald}       {cmyk}{1   , 0   , 0.50, 0   }
\definecolor{junglegreen}   {cmyk}{0.99, 0   , 0.52, 0   }
\definecolor{seagreen}      {cmyk}{0.69, 0   , 0.50, 0   }
\definecolor{green}         {cmyk}{1   , 0   , 1   , 0   }
\definecolor{forestgreen}   {cmyk}{0.91, 0   , 0.88, 0.12}
\definecolor{pinegreen}     {cmyk}{0.92, 0   , 0.59, 0.25}
\definecolor{limegreen}     {cmyk}{0.50, 0   , 1   , 0   }
\definecolor{yellowgreen}   {cmyk}{0.44, 0   , 0.74, 0   }
\definecolor{springgreen}   {cmyk}{0.26, 0   , 0.76, 0   }
\definecolor{olivegreen}    {cmyk}{0.64, 0   , 0.95, 0.40}
\definecolor{rawsienna}     {cmyk}{0   , 0.72, 1   , 0.45}
\definecolor{sepia}         {cmyk}{0   , 0.83, 1   , 0.70}
\definecolor{brown}         {cmyk}{0   , 0.81, 1   , 0.60}
\definecolor{tan}           {cmyk}{0.14, 0.42, 0.56, 0   }
\definecolor{gray}          {cmyk}{0   , 0   , 0   , 0.50}
\definecolor{black}         {cmyk}{0   , 0   , 0   , 1   }
\definecolor{white}         {cmyk}{0   , 0   , 0   , 0   }
\definecolor{tuklblue}{RGB}{0,95,140}
\definecolor{tuklred}{RGB}{185,40,25} 

\usepgfplotslibrary{groupplots}

\pgfplotsset{
	colormap={jet}{
rgb(0.000000 pt)=(0.000000,0.000000,0.504000);
rgb(1.000000 pt)=(0.000000,0.000000,0.508000);
rgb(2.000000 pt)=(0.000000,0.000000,0.512000);
rgb(3.000000 pt)=(0.000000,0.000000,0.516000);
rgb(4.000000 pt)=(0.000000,0.000000,0.520000);
rgb(5.000000 pt)=(0.000000,0.000000,0.524000);
rgb(6.000000 pt)=(0.000000,0.000000,0.528000);
rgb(7.000000 pt)=(0.000000,0.000000,0.532000);
rgb(8.000000 pt)=(0.000000,0.000000,0.536000);
rgb(9.000000 pt)=(0.000000,0.000000,0.540000);
rgb(10.000000 pt)=(0.000000,0.000000,0.544000);
rgb(11.000000 pt)=(0.000000,0.000000,0.548000);
rgb(12.000000 pt)=(0.000000,0.000000,0.552000);
rgb(13.000000 pt)=(0.000000,0.000000,0.556000);
rgb(14.000000 pt)=(0.000000,0.000000,0.560000);
rgb(15.000000 pt)=(0.000000,0.000000,0.564000);
rgb(16.000000 pt)=(0.000000,0.000000,0.568000);
rgb(17.000000 pt)=(0.000000,0.000000,0.572000);
rgb(18.000000 pt)=(0.000000,0.000000,0.576000);
rgb(19.000000 pt)=(0.000000,0.000000,0.580000);
rgb(20.000000 pt)=(0.000000,0.000000,0.584000);
rgb(21.000000 pt)=(0.000000,0.000000,0.588000);
rgb(22.000000 pt)=(0.000000,0.000000,0.592000);
rgb(23.000000 pt)=(0.000000,0.000000,0.596000);
rgb(24.000000 pt)=(0.000000,0.000000,0.600000);
rgb(25.000000 pt)=(0.000000,0.000000,0.604000);
rgb(26.000000 pt)=(0.000000,0.000000,0.608000);
rgb(27.000000 pt)=(0.000000,0.000000,0.612000);
rgb(28.000000 pt)=(0.000000,0.000000,0.616000);
rgb(29.000000 pt)=(0.000000,0.000000,0.620000);
rgb(30.000000 pt)=(0.000000,0.000000,0.624000);
rgb(31.000000 pt)=(0.000000,0.000000,0.628000);
rgb(32.000000 pt)=(0.000000,0.000000,0.632000);
rgb(33.000000 pt)=(0.000000,0.000000,0.636000);
rgb(34.000000 pt)=(0.000000,0.000000,0.640000);
rgb(35.000000 pt)=(0.000000,0.000000,0.644000);
rgb(36.000000 pt)=(0.000000,0.000000,0.648000);
rgb(37.000000 pt)=(0.000000,0.000000,0.652000);
rgb(38.000000 pt)=(0.000000,0.000000,0.656000);
rgb(39.000000 pt)=(0.000000,0.000000,0.660000);
rgb(40.000000 pt)=(0.000000,0.000000,0.664000);
rgb(41.000000 pt)=(0.000000,0.000000,0.668000);
rgb(42.000000 pt)=(0.000000,0.000000,0.672000);
rgb(43.000000 pt)=(0.000000,0.000000,0.676000);
rgb(44.000000 pt)=(0.000000,0.000000,0.680000);
rgb(45.000000 pt)=(0.000000,0.000000,0.684000);
rgb(46.000000 pt)=(0.000000,0.000000,0.688000);
rgb(47.000000 pt)=(0.000000,0.000000,0.692000);
rgb(48.000000 pt)=(0.000000,0.000000,0.696000);
rgb(49.000000 pt)=(0.000000,0.000000,0.700000);
rgb(50.000000 pt)=(0.000000,0.000000,0.704000);
rgb(51.000000 pt)=(0.000000,0.000000,0.708000);
rgb(52.000000 pt)=(0.000000,0.000000,0.712000);
rgb(53.000000 pt)=(0.000000,0.000000,0.716000);
rgb(54.000000 pt)=(0.000000,0.000000,0.720000);
rgb(55.000000 pt)=(0.000000,0.000000,0.724000);
rgb(56.000000 pt)=(0.000000,0.000000,0.728000);
rgb(57.000000 pt)=(0.000000,0.000000,0.732000);
rgb(58.000000 pt)=(0.000000,0.000000,0.736000);
rgb(59.000000 pt)=(0.000000,0.000000,0.740000);
rgb(60.000000 pt)=(0.000000,0.000000,0.744000);
rgb(61.000000 pt)=(0.000000,0.000000,0.748000);
rgb(62.000000 pt)=(0.000000,0.000000,0.752000);
rgb(63.000000 pt)=(0.000000,0.000000,0.756000);
rgb(64.000000 pt)=(0.000000,0.000000,0.760000);
rgb(65.000000 pt)=(0.000000,0.000000,0.764000);
rgb(66.000000 pt)=(0.000000,0.000000,0.768000);
rgb(67.000000 pt)=(0.000000,0.000000,0.772000);
rgb(68.000000 pt)=(0.000000,0.000000,0.776000);
rgb(69.000000 pt)=(0.000000,0.000000,0.780000);
rgb(70.000000 pt)=(0.000000,0.000000,0.784000);
rgb(71.000000 pt)=(0.000000,0.000000,0.788000);
rgb(72.000000 pt)=(0.000000,0.000000,0.792000);
rgb(73.000000 pt)=(0.000000,0.000000,0.796000);
rgb(74.000000 pt)=(0.000000,0.000000,0.800000);
rgb(75.000000 pt)=(0.000000,0.000000,0.804000);
rgb(76.000000 pt)=(0.000000,0.000000,0.808000);
rgb(77.000000 pt)=(0.000000,0.000000,0.812000);
rgb(78.000000 pt)=(0.000000,0.000000,0.816000);
rgb(79.000000 pt)=(0.000000,0.000000,0.820000);
rgb(80.000000 pt)=(0.000000,0.000000,0.824000);
rgb(81.000000 pt)=(0.000000,0.000000,0.828000);
rgb(82.000000 pt)=(0.000000,0.000000,0.832000);
rgb(83.000000 pt)=(0.000000,0.000000,0.836000);
rgb(84.000000 pt)=(0.000000,0.000000,0.840000);
rgb(85.000000 pt)=(0.000000,0.000000,0.844000);
rgb(86.000000 pt)=(0.000000,0.000000,0.848000);
rgb(87.000000 pt)=(0.000000,0.000000,0.852000);
rgb(88.000000 pt)=(0.000000,0.000000,0.856000);
rgb(89.000000 pt)=(0.000000,0.000000,0.860000);
rgb(90.000000 pt)=(0.000000,0.000000,0.864000);
rgb(91.000000 pt)=(0.000000,0.000000,0.868000);
rgb(92.000000 pt)=(0.000000,0.000000,0.872000);
rgb(93.000000 pt)=(0.000000,0.000000,0.876000);
rgb(94.000000 pt)=(0.000000,0.000000,0.880000);
rgb(95.000000 pt)=(0.000000,0.000000,0.884000);
rgb(96.000000 pt)=(0.000000,0.000000,0.888000);
rgb(97.000000 pt)=(0.000000,0.000000,0.892000);
rgb(98.000000 pt)=(0.000000,0.000000,0.896000);
rgb(99.000000 pt)=(0.000000,0.000000,0.900000);
rgb(100.000000 pt)=(0.000000,0.000000,0.904000);
rgb(101.000000 pt)=(0.000000,0.000000,0.908000);
rgb(102.000000 pt)=(0.000000,0.000000,0.912000);
rgb(103.000000 pt)=(0.000000,0.000000,0.916000);
rgb(104.000000 pt)=(0.000000,0.000000,0.920000);
rgb(105.000000 pt)=(0.000000,0.000000,0.924000);
rgb(106.000000 pt)=(0.000000,0.000000,0.928000);
rgb(107.000000 pt)=(0.000000,0.000000,0.932000);
rgb(108.000000 pt)=(0.000000,0.000000,0.936000);
rgb(109.000000 pt)=(0.000000,0.000000,0.940000);
rgb(110.000000 pt)=(0.000000,0.000000,0.944000);
rgb(111.000000 pt)=(0.000000,0.000000,0.948000);
rgb(112.000000 pt)=(0.000000,0.000000,0.952000);
rgb(113.000000 pt)=(0.000000,0.000000,0.956000);
rgb(114.000000 pt)=(0.000000,0.000000,0.960000);
rgb(115.000000 pt)=(0.000000,0.000000,0.964000);
rgb(116.000000 pt)=(0.000000,0.000000,0.968000);
rgb(117.000000 pt)=(0.000000,0.000000,0.972000);
rgb(118.000000 pt)=(0.000000,0.000000,0.976000);
rgb(119.000000 pt)=(0.000000,0.000000,0.980000);
rgb(120.000000 pt)=(0.000000,0.000000,0.984000);
rgb(121.000000 pt)=(0.000000,0.000000,0.988000);
rgb(122.000000 pt)=(0.000000,0.000000,0.992000);
rgb(123.000000 pt)=(0.000000,0.000000,0.996000);
rgb(124.000000 pt)=(0.000000,0.000000,1.000000);
rgb(125.000000 pt)=(0.000000,0.004000,1.000000);
rgb(126.000000 pt)=(0.000000,0.008000,1.000000);
rgb(127.000000 pt)=(0.000000,0.012000,1.000000);
rgb(128.000000 pt)=(0.000000,0.016000,1.000000);
rgb(129.000000 pt)=(0.000000,0.020000,1.000000);
rgb(130.000000 pt)=(0.000000,0.024000,1.000000);
rgb(131.000000 pt)=(0.000000,0.028000,1.000000);
rgb(132.000000 pt)=(0.000000,0.032000,1.000000);
rgb(133.000000 pt)=(0.000000,0.036000,1.000000);
rgb(134.000000 pt)=(0.000000,0.040000,1.000000);
rgb(135.000000 pt)=(0.000000,0.044000,1.000000);
rgb(136.000000 pt)=(0.000000,0.048000,1.000000);
rgb(137.000000 pt)=(0.000000,0.052000,1.000000);
rgb(138.000000 pt)=(0.000000,0.056000,1.000000);
rgb(139.000000 pt)=(0.000000,0.060000,1.000000);
rgb(140.000000 pt)=(0.000000,0.064000,1.000000);
rgb(141.000000 pt)=(0.000000,0.068000,1.000000);
rgb(142.000000 pt)=(0.000000,0.072000,1.000000);
rgb(143.000000 pt)=(0.000000,0.076000,1.000000);
rgb(144.000000 pt)=(0.000000,0.080000,1.000000);
rgb(145.000000 pt)=(0.000000,0.084000,1.000000);
rgb(146.000000 pt)=(0.000000,0.088000,1.000000);
rgb(147.000000 pt)=(0.000000,0.092000,1.000000);
rgb(148.000000 pt)=(0.000000,0.096000,1.000000);
rgb(149.000000 pt)=(0.000000,0.100000,1.000000);
rgb(150.000000 pt)=(0.000000,0.104000,1.000000);
rgb(151.000000 pt)=(0.000000,0.108000,1.000000);
rgb(152.000000 pt)=(0.000000,0.112000,1.000000);
rgb(153.000000 pt)=(0.000000,0.116000,1.000000);
rgb(154.000000 pt)=(0.000000,0.120000,1.000000);
rgb(155.000000 pt)=(0.000000,0.124000,1.000000);
rgb(156.000000 pt)=(0.000000,0.128000,1.000000);
rgb(157.000000 pt)=(0.000000,0.132000,1.000000);
rgb(158.000000 pt)=(0.000000,0.136000,1.000000);
rgb(159.000000 pt)=(0.000000,0.140000,1.000000);
rgb(160.000000 pt)=(0.000000,0.144000,1.000000);
rgb(161.000000 pt)=(0.000000,0.148000,1.000000);
rgb(162.000000 pt)=(0.000000,0.152000,1.000000);
rgb(163.000000 pt)=(0.000000,0.156000,1.000000);
rgb(164.000000 pt)=(0.000000,0.160000,1.000000);
rgb(165.000000 pt)=(0.000000,0.164000,1.000000);
rgb(166.000000 pt)=(0.000000,0.168000,1.000000);
rgb(167.000000 pt)=(0.000000,0.172000,1.000000);
rgb(168.000000 pt)=(0.000000,0.176000,1.000000);
rgb(169.000000 pt)=(0.000000,0.180000,1.000000);
rgb(170.000000 pt)=(0.000000,0.184000,1.000000);
rgb(171.000000 pt)=(0.000000,0.188000,1.000000);
rgb(172.000000 pt)=(0.000000,0.192000,1.000000);
rgb(173.000000 pt)=(0.000000,0.196000,1.000000);
rgb(174.000000 pt)=(0.000000,0.200000,1.000000);
rgb(175.000000 pt)=(0.000000,0.204000,1.000000);
rgb(176.000000 pt)=(0.000000,0.208000,1.000000);
rgb(177.000000 pt)=(0.000000,0.212000,1.000000);
rgb(178.000000 pt)=(0.000000,0.216000,1.000000);
rgb(179.000000 pt)=(0.000000,0.220000,1.000000);
rgb(180.000000 pt)=(0.000000,0.224000,1.000000);
rgb(181.000000 pt)=(0.000000,0.228000,1.000000);
rgb(182.000000 pt)=(0.000000,0.232000,1.000000);
rgb(183.000000 pt)=(0.000000,0.236000,1.000000);
rgb(184.000000 pt)=(0.000000,0.240000,1.000000);
rgb(185.000000 pt)=(0.000000,0.244000,1.000000);
rgb(186.000000 pt)=(0.000000,0.248000,1.000000);
rgb(187.000000 pt)=(0.000000,0.252000,1.000000);
rgb(188.000000 pt)=(0.000000,0.256000,1.000000);
rgb(189.000000 pt)=(0.000000,0.260000,1.000000);
rgb(190.000000 pt)=(0.000000,0.264000,1.000000);
rgb(191.000000 pt)=(0.000000,0.268000,1.000000);
rgb(192.000000 pt)=(0.000000,0.272000,1.000000);
rgb(193.000000 pt)=(0.000000,0.276000,1.000000);
rgb(194.000000 pt)=(0.000000,0.280000,1.000000);
rgb(195.000000 pt)=(0.000000,0.284000,1.000000);
rgb(196.000000 pt)=(0.000000,0.288000,1.000000);
rgb(197.000000 pt)=(0.000000,0.292000,1.000000);
rgb(198.000000 pt)=(0.000000,0.296000,1.000000);
rgb(199.000000 pt)=(0.000000,0.300000,1.000000);
rgb(200.000000 pt)=(0.000000,0.304000,1.000000);
rgb(201.000000 pt)=(0.000000,0.308000,1.000000);
rgb(202.000000 pt)=(0.000000,0.312000,1.000000);
rgb(203.000000 pt)=(0.000000,0.316000,1.000000);
rgb(204.000000 pt)=(0.000000,0.320000,1.000000);
rgb(205.000000 pt)=(0.000000,0.324000,1.000000);
rgb(206.000000 pt)=(0.000000,0.328000,1.000000);
rgb(207.000000 pt)=(0.000000,0.332000,1.000000);
rgb(208.000000 pt)=(0.000000,0.336000,1.000000);
rgb(209.000000 pt)=(0.000000,0.340000,1.000000);
rgb(210.000000 pt)=(0.000000,0.344000,1.000000);
rgb(211.000000 pt)=(0.000000,0.348000,1.000000);
rgb(212.000000 pt)=(0.000000,0.352000,1.000000);
rgb(213.000000 pt)=(0.000000,0.356000,1.000000);
rgb(214.000000 pt)=(0.000000,0.360000,1.000000);
rgb(215.000000 pt)=(0.000000,0.364000,1.000000);
rgb(216.000000 pt)=(0.000000,0.368000,1.000000);
rgb(217.000000 pt)=(0.000000,0.372000,1.000000);
rgb(218.000000 pt)=(0.000000,0.376000,1.000000);
rgb(219.000000 pt)=(0.000000,0.380000,1.000000);
rgb(220.000000 pt)=(0.000000,0.384000,1.000000);
rgb(221.000000 pt)=(0.000000,0.388000,1.000000);
rgb(222.000000 pt)=(0.000000,0.392000,1.000000);
rgb(223.000000 pt)=(0.000000,0.396000,1.000000);
rgb(224.000000 pt)=(0.000000,0.400000,1.000000);
rgb(225.000000 pt)=(0.000000,0.404000,1.000000);
rgb(226.000000 pt)=(0.000000,0.408000,1.000000);
rgb(227.000000 pt)=(0.000000,0.412000,1.000000);
rgb(228.000000 pt)=(0.000000,0.416000,1.000000);
rgb(229.000000 pt)=(0.000000,0.420000,1.000000);
rgb(230.000000 pt)=(0.000000,0.424000,1.000000);
rgb(231.000000 pt)=(0.000000,0.428000,1.000000);
rgb(232.000000 pt)=(0.000000,0.432000,1.000000);
rgb(233.000000 pt)=(0.000000,0.436000,1.000000);
rgb(234.000000 pt)=(0.000000,0.440000,1.000000);
rgb(235.000000 pt)=(0.000000,0.444000,1.000000);
rgb(236.000000 pt)=(0.000000,0.448000,1.000000);
rgb(237.000000 pt)=(0.000000,0.452000,1.000000);
rgb(238.000000 pt)=(0.000000,0.456000,1.000000);
rgb(239.000000 pt)=(0.000000,0.460000,1.000000);
rgb(240.000000 pt)=(0.000000,0.464000,1.000000);
rgb(241.000000 pt)=(0.000000,0.468000,1.000000);
rgb(242.000000 pt)=(0.000000,0.472000,1.000000);
rgb(243.000000 pt)=(0.000000,0.476000,1.000000);
rgb(244.000000 pt)=(0.000000,0.480000,1.000000);
rgb(245.000000 pt)=(0.000000,0.484000,1.000000);
rgb(246.000000 pt)=(0.000000,0.488000,1.000000);
rgb(247.000000 pt)=(0.000000,0.492000,1.000000);
rgb(248.000000 pt)=(0.000000,0.496000,1.000000);
rgb(249.000000 pt)=(0.000000,0.500000,1.000000);
rgb(250.000000 pt)=(0.000000,0.504000,1.000000);
rgb(251.000000 pt)=(0.000000,0.508000,1.000000);
rgb(252.000000 pt)=(0.000000,0.512000,1.000000);
rgb(253.000000 pt)=(0.000000,0.516000,1.000000);
rgb(254.000000 pt)=(0.000000,0.520000,1.000000);
rgb(255.000000 pt)=(0.000000,0.524000,1.000000);
rgb(256.000000 pt)=(0.000000,0.528000,1.000000);
rgb(257.000000 pt)=(0.000000,0.532000,1.000000);
rgb(258.000000 pt)=(0.000000,0.536000,1.000000);
rgb(259.000000 pt)=(0.000000,0.540000,1.000000);
rgb(260.000000 pt)=(0.000000,0.544000,1.000000);
rgb(261.000000 pt)=(0.000000,0.548000,1.000000);
rgb(262.000000 pt)=(0.000000,0.552000,1.000000);
rgb(263.000000 pt)=(0.000000,0.556000,1.000000);
rgb(264.000000 pt)=(0.000000,0.560000,1.000000);
rgb(265.000000 pt)=(0.000000,0.564000,1.000000);
rgb(266.000000 pt)=(0.000000,0.568000,1.000000);
rgb(267.000000 pt)=(0.000000,0.572000,1.000000);
rgb(268.000000 pt)=(0.000000,0.576000,1.000000);
rgb(269.000000 pt)=(0.000000,0.580000,1.000000);
rgb(270.000000 pt)=(0.000000,0.584000,1.000000);
rgb(271.000000 pt)=(0.000000,0.588000,1.000000);
rgb(272.000000 pt)=(0.000000,0.592000,1.000000);
rgb(273.000000 pt)=(0.000000,0.596000,1.000000);
rgb(274.000000 pt)=(0.000000,0.600000,1.000000);
rgb(275.000000 pt)=(0.000000,0.604000,1.000000);
rgb(276.000000 pt)=(0.000000,0.608000,1.000000);
rgb(277.000000 pt)=(0.000000,0.612000,1.000000);
rgb(278.000000 pt)=(0.000000,0.616000,1.000000);
rgb(279.000000 pt)=(0.000000,0.620000,1.000000);
rgb(280.000000 pt)=(0.000000,0.624000,1.000000);
rgb(281.000000 pt)=(0.000000,0.628000,1.000000);
rgb(282.000000 pt)=(0.000000,0.632000,1.000000);
rgb(283.000000 pt)=(0.000000,0.636000,1.000000);
rgb(284.000000 pt)=(0.000000,0.640000,1.000000);
rgb(285.000000 pt)=(0.000000,0.644000,1.000000);
rgb(286.000000 pt)=(0.000000,0.648000,1.000000);
rgb(287.000000 pt)=(0.000000,0.652000,1.000000);
rgb(288.000000 pt)=(0.000000,0.656000,1.000000);
rgb(289.000000 pt)=(0.000000,0.660000,1.000000);
rgb(290.000000 pt)=(0.000000,0.664000,1.000000);
rgb(291.000000 pt)=(0.000000,0.668000,1.000000);
rgb(292.000000 pt)=(0.000000,0.672000,1.000000);
rgb(293.000000 pt)=(0.000000,0.676000,1.000000);
rgb(294.000000 pt)=(0.000000,0.680000,1.000000);
rgb(295.000000 pt)=(0.000000,0.684000,1.000000);
rgb(296.000000 pt)=(0.000000,0.688000,1.000000);
rgb(297.000000 pt)=(0.000000,0.692000,1.000000);
rgb(298.000000 pt)=(0.000000,0.696000,1.000000);
rgb(299.000000 pt)=(0.000000,0.700000,1.000000);
rgb(300.000000 pt)=(0.000000,0.704000,1.000000);
rgb(301.000000 pt)=(0.000000,0.708000,1.000000);
rgb(302.000000 pt)=(0.000000,0.712000,1.000000);
rgb(303.000000 pt)=(0.000000,0.716000,1.000000);
rgb(304.000000 pt)=(0.000000,0.720000,1.000000);
rgb(305.000000 pt)=(0.000000,0.724000,1.000000);
rgb(306.000000 pt)=(0.000000,0.728000,1.000000);
rgb(307.000000 pt)=(0.000000,0.732000,1.000000);
rgb(308.000000 pt)=(0.000000,0.736000,1.000000);
rgb(309.000000 pt)=(0.000000,0.740000,1.000000);
rgb(310.000000 pt)=(0.000000,0.744000,1.000000);
rgb(311.000000 pt)=(0.000000,0.748000,1.000000);
rgb(312.000000 pt)=(0.000000,0.752000,1.000000);
rgb(313.000000 pt)=(0.000000,0.756000,1.000000);
rgb(314.000000 pt)=(0.000000,0.760000,1.000000);
rgb(315.000000 pt)=(0.000000,0.764000,1.000000);
rgb(316.000000 pt)=(0.000000,0.768000,1.000000);
rgb(317.000000 pt)=(0.000000,0.772000,1.000000);
rgb(318.000000 pt)=(0.000000,0.776000,1.000000);
rgb(319.000000 pt)=(0.000000,0.780000,1.000000);
rgb(320.000000 pt)=(0.000000,0.784000,1.000000);
rgb(321.000000 pt)=(0.000000,0.788000,1.000000);
rgb(322.000000 pt)=(0.000000,0.792000,1.000000);
rgb(323.000000 pt)=(0.000000,0.796000,1.000000);
rgb(324.000000 pt)=(0.000000,0.800000,1.000000);
rgb(325.000000 pt)=(0.000000,0.804000,1.000000);
rgb(326.000000 pt)=(0.000000,0.808000,1.000000);
rgb(327.000000 pt)=(0.000000,0.812000,1.000000);
rgb(328.000000 pt)=(0.000000,0.816000,1.000000);
rgb(329.000000 pt)=(0.000000,0.820000,1.000000);
rgb(330.000000 pt)=(0.000000,0.824000,1.000000);
rgb(331.000000 pt)=(0.000000,0.828000,1.000000);
rgb(332.000000 pt)=(0.000000,0.832000,1.000000);
rgb(333.000000 pt)=(0.000000,0.836000,1.000000);
rgb(334.000000 pt)=(0.000000,0.840000,1.000000);
rgb(335.000000 pt)=(0.000000,0.844000,1.000000);
rgb(336.000000 pt)=(0.000000,0.848000,1.000000);
rgb(337.000000 pt)=(0.000000,0.852000,1.000000);
rgb(338.000000 pt)=(0.000000,0.856000,1.000000);
rgb(339.000000 pt)=(0.000000,0.860000,1.000000);
rgb(340.000000 pt)=(0.000000,0.864000,1.000000);
rgb(341.000000 pt)=(0.000000,0.868000,1.000000);
rgb(342.000000 pt)=(0.000000,0.872000,1.000000);
rgb(343.000000 pt)=(0.000000,0.876000,1.000000);
rgb(344.000000 pt)=(0.000000,0.880000,1.000000);
rgb(345.000000 pt)=(0.000000,0.884000,1.000000);
rgb(346.000000 pt)=(0.000000,0.888000,1.000000);
rgb(347.000000 pt)=(0.000000,0.892000,1.000000);
rgb(348.000000 pt)=(0.000000,0.896000,1.000000);
rgb(349.000000 pt)=(0.000000,0.900000,1.000000);
rgb(350.000000 pt)=(0.000000,0.904000,1.000000);
rgb(351.000000 pt)=(0.000000,0.908000,1.000000);
rgb(352.000000 pt)=(0.000000,0.912000,1.000000);
rgb(353.000000 pt)=(0.000000,0.916000,1.000000);
rgb(354.000000 pt)=(0.000000,0.920000,1.000000);
rgb(355.000000 pt)=(0.000000,0.924000,1.000000);
rgb(356.000000 pt)=(0.000000,0.928000,1.000000);
rgb(357.000000 pt)=(0.000000,0.932000,1.000000);
rgb(358.000000 pt)=(0.000000,0.936000,1.000000);
rgb(359.000000 pt)=(0.000000,0.940000,1.000000);
rgb(360.000000 pt)=(0.000000,0.944000,1.000000);
rgb(361.000000 pt)=(0.000000,0.948000,1.000000);
rgb(362.000000 pt)=(0.000000,0.952000,1.000000);
rgb(363.000000 pt)=(0.000000,0.956000,1.000000);
rgb(364.000000 pt)=(0.000000,0.960000,1.000000);
rgb(365.000000 pt)=(0.000000,0.964000,1.000000);
rgb(366.000000 pt)=(0.000000,0.968000,1.000000);
rgb(367.000000 pt)=(0.000000,0.972000,1.000000);
rgb(368.000000 pt)=(0.000000,0.976000,1.000000);
rgb(369.000000 pt)=(0.000000,0.980000,1.000000);
rgb(370.000000 pt)=(0.000000,0.984000,1.000000);
rgb(371.000000 pt)=(0.000000,0.988000,1.000000);
rgb(372.000000 pt)=(0.000000,0.992000,1.000000);
rgb(373.000000 pt)=(0.000000,0.996000,1.000000);
rgb(374.000000 pt)=(0.000000,1.000000,1.000000);
rgb(375.000000 pt)=(0.004000,1.000000,0.996000);
rgb(376.000000 pt)=(0.008000,1.000000,0.992000);
rgb(377.000000 pt)=(0.012000,1.000000,0.988000);
rgb(378.000000 pt)=(0.016000,1.000000,0.984000);
rgb(379.000000 pt)=(0.020000,1.000000,0.980000);
rgb(380.000000 pt)=(0.024000,1.000000,0.976000);
rgb(381.000000 pt)=(0.028000,1.000000,0.972000);
rgb(382.000000 pt)=(0.032000,1.000000,0.968000);
rgb(383.000000 pt)=(0.036000,1.000000,0.964000);
rgb(384.000000 pt)=(0.040000,1.000000,0.960000);
rgb(385.000000 pt)=(0.044000,1.000000,0.956000);
rgb(386.000000 pt)=(0.048000,1.000000,0.952000);
rgb(387.000000 pt)=(0.052000,1.000000,0.948000);
rgb(388.000000 pt)=(0.056000,1.000000,0.944000);
rgb(389.000000 pt)=(0.060000,1.000000,0.940000);
rgb(390.000000 pt)=(0.064000,1.000000,0.936000);
rgb(391.000000 pt)=(0.068000,1.000000,0.932000);
rgb(392.000000 pt)=(0.072000,1.000000,0.928000);
rgb(393.000000 pt)=(0.076000,1.000000,0.924000);
rgb(394.000000 pt)=(0.080000,1.000000,0.920000);
rgb(395.000000 pt)=(0.084000,1.000000,0.916000);
rgb(396.000000 pt)=(0.088000,1.000000,0.912000);
rgb(397.000000 pt)=(0.092000,1.000000,0.908000);
rgb(398.000000 pt)=(0.096000,1.000000,0.904000);
rgb(399.000000 pt)=(0.100000,1.000000,0.900000);
rgb(400.000000 pt)=(0.104000,1.000000,0.896000);
rgb(401.000000 pt)=(0.108000,1.000000,0.892000);
rgb(402.000000 pt)=(0.112000,1.000000,0.888000);
rgb(403.000000 pt)=(0.116000,1.000000,0.884000);
rgb(404.000000 pt)=(0.120000,1.000000,0.880000);
rgb(405.000000 pt)=(0.124000,1.000000,0.876000);
rgb(406.000000 pt)=(0.128000,1.000000,0.872000);
rgb(407.000000 pt)=(0.132000,1.000000,0.868000);
rgb(408.000000 pt)=(0.136000,1.000000,0.864000);
rgb(409.000000 pt)=(0.140000,1.000000,0.860000);
rgb(410.000000 pt)=(0.144000,1.000000,0.856000);
rgb(411.000000 pt)=(0.148000,1.000000,0.852000);
rgb(412.000000 pt)=(0.152000,1.000000,0.848000);
rgb(413.000000 pt)=(0.156000,1.000000,0.844000);
rgb(414.000000 pt)=(0.160000,1.000000,0.840000);
rgb(415.000000 pt)=(0.164000,1.000000,0.836000);
rgb(416.000000 pt)=(0.168000,1.000000,0.832000);
rgb(417.000000 pt)=(0.172000,1.000000,0.828000);
rgb(418.000000 pt)=(0.176000,1.000000,0.824000);
rgb(419.000000 pt)=(0.180000,1.000000,0.820000);
rgb(420.000000 pt)=(0.184000,1.000000,0.816000);
rgb(421.000000 pt)=(0.188000,1.000000,0.812000);
rgb(422.000000 pt)=(0.192000,1.000000,0.808000);
rgb(423.000000 pt)=(0.196000,1.000000,0.804000);
rgb(424.000000 pt)=(0.200000,1.000000,0.800000);
rgb(425.000000 pt)=(0.204000,1.000000,0.796000);
rgb(426.000000 pt)=(0.208000,1.000000,0.792000);
rgb(427.000000 pt)=(0.212000,1.000000,0.788000);
rgb(428.000000 pt)=(0.216000,1.000000,0.784000);
rgb(429.000000 pt)=(0.220000,1.000000,0.780000);
rgb(430.000000 pt)=(0.224000,1.000000,0.776000);
rgb(431.000000 pt)=(0.228000,1.000000,0.772000);
rgb(432.000000 pt)=(0.232000,1.000000,0.768000);
rgb(433.000000 pt)=(0.236000,1.000000,0.764000);
rgb(434.000000 pt)=(0.240000,1.000000,0.760000);
rgb(435.000000 pt)=(0.244000,1.000000,0.756000);
rgb(436.000000 pt)=(0.248000,1.000000,0.752000);
rgb(437.000000 pt)=(0.252000,1.000000,0.748000);
rgb(438.000000 pt)=(0.256000,1.000000,0.744000);
rgb(439.000000 pt)=(0.260000,1.000000,0.740000);
rgb(440.000000 pt)=(0.264000,1.000000,0.736000);
rgb(441.000000 pt)=(0.268000,1.000000,0.732000);
rgb(442.000000 pt)=(0.272000,1.000000,0.728000);
rgb(443.000000 pt)=(0.276000,1.000000,0.724000);
rgb(444.000000 pt)=(0.280000,1.000000,0.720000);
rgb(445.000000 pt)=(0.284000,1.000000,0.716000);
rgb(446.000000 pt)=(0.288000,1.000000,0.712000);
rgb(447.000000 pt)=(0.292000,1.000000,0.708000);
rgb(448.000000 pt)=(0.296000,1.000000,0.704000);
rgb(449.000000 pt)=(0.300000,1.000000,0.700000);
rgb(450.000000 pt)=(0.304000,1.000000,0.696000);
rgb(451.000000 pt)=(0.308000,1.000000,0.692000);
rgb(452.000000 pt)=(0.312000,1.000000,0.688000);
rgb(453.000000 pt)=(0.316000,1.000000,0.684000);
rgb(454.000000 pt)=(0.320000,1.000000,0.680000);
rgb(455.000000 pt)=(0.324000,1.000000,0.676000);
rgb(456.000000 pt)=(0.328000,1.000000,0.672000);
rgb(457.000000 pt)=(0.332000,1.000000,0.668000);
rgb(458.000000 pt)=(0.336000,1.000000,0.664000);
rgb(459.000000 pt)=(0.340000,1.000000,0.660000);
rgb(460.000000 pt)=(0.344000,1.000000,0.656000);
rgb(461.000000 pt)=(0.348000,1.000000,0.652000);
rgb(462.000000 pt)=(0.352000,1.000000,0.648000);
rgb(463.000000 pt)=(0.356000,1.000000,0.644000);
rgb(464.000000 pt)=(0.360000,1.000000,0.640000);
rgb(465.000000 pt)=(0.364000,1.000000,0.636000);
rgb(466.000000 pt)=(0.368000,1.000000,0.632000);
rgb(467.000000 pt)=(0.372000,1.000000,0.628000);
rgb(468.000000 pt)=(0.376000,1.000000,0.624000);
rgb(469.000000 pt)=(0.380000,1.000000,0.620000);
rgb(470.000000 pt)=(0.384000,1.000000,0.616000);
rgb(471.000000 pt)=(0.388000,1.000000,0.612000);
rgb(472.000000 pt)=(0.392000,1.000000,0.608000);
rgb(473.000000 pt)=(0.396000,1.000000,0.604000);
rgb(474.000000 pt)=(0.400000,1.000000,0.600000);
rgb(475.000000 pt)=(0.404000,1.000000,0.596000);
rgb(476.000000 pt)=(0.408000,1.000000,0.592000);
rgb(477.000000 pt)=(0.412000,1.000000,0.588000);
rgb(478.000000 pt)=(0.416000,1.000000,0.584000);
rgb(479.000000 pt)=(0.420000,1.000000,0.580000);
rgb(480.000000 pt)=(0.424000,1.000000,0.576000);
rgb(481.000000 pt)=(0.428000,1.000000,0.572000);
rgb(482.000000 pt)=(0.432000,1.000000,0.568000);
rgb(483.000000 pt)=(0.436000,1.000000,0.564000);
rgb(484.000000 pt)=(0.440000,1.000000,0.560000);
rgb(485.000000 pt)=(0.444000,1.000000,0.556000);
rgb(486.000000 pt)=(0.448000,1.000000,0.552000);
rgb(487.000000 pt)=(0.452000,1.000000,0.548000);
rgb(488.000000 pt)=(0.456000,1.000000,0.544000);
rgb(489.000000 pt)=(0.460000,1.000000,0.540000);
rgb(490.000000 pt)=(0.464000,1.000000,0.536000);
rgb(491.000000 pt)=(0.468000,1.000000,0.532000);
rgb(492.000000 pt)=(0.472000,1.000000,0.528000);
rgb(493.000000 pt)=(0.476000,1.000000,0.524000);
rgb(494.000000 pt)=(0.480000,1.000000,0.520000);
rgb(495.000000 pt)=(0.484000,1.000000,0.516000);
rgb(496.000000 pt)=(0.488000,1.000000,0.512000);
rgb(497.000000 pt)=(0.492000,1.000000,0.508000);
rgb(498.000000 pt)=(0.496000,1.000000,0.504000);
rgb(499.000000 pt)=(0.500000,1.000000,0.500000);
rgb(500.000000 pt)=(0.504000,1.000000,0.496000);
rgb(501.000000 pt)=(0.508000,1.000000,0.492000);
rgb(502.000000 pt)=(0.512000,1.000000,0.488000);
rgb(503.000000 pt)=(0.516000,1.000000,0.484000);
rgb(504.000000 pt)=(0.520000,1.000000,0.480000);
rgb(505.000000 pt)=(0.524000,1.000000,0.476000);
rgb(506.000000 pt)=(0.528000,1.000000,0.472000);
rgb(507.000000 pt)=(0.532000,1.000000,0.468000);
rgb(508.000000 pt)=(0.536000,1.000000,0.464000);
rgb(509.000000 pt)=(0.540000,1.000000,0.460000);
rgb(510.000000 pt)=(0.544000,1.000000,0.456000);
rgb(511.000000 pt)=(0.548000,1.000000,0.452000);
rgb(512.000000 pt)=(0.552000,1.000000,0.448000);
rgb(513.000000 pt)=(0.556000,1.000000,0.444000);
rgb(514.000000 pt)=(0.560000,1.000000,0.440000);
rgb(515.000000 pt)=(0.564000,1.000000,0.436000);
rgb(516.000000 pt)=(0.568000,1.000000,0.432000);
rgb(517.000000 pt)=(0.572000,1.000000,0.428000);
rgb(518.000000 pt)=(0.576000,1.000000,0.424000);
rgb(519.000000 pt)=(0.580000,1.000000,0.420000);
rgb(520.000000 pt)=(0.584000,1.000000,0.416000);
rgb(521.000000 pt)=(0.588000,1.000000,0.412000);
rgb(522.000000 pt)=(0.592000,1.000000,0.408000);
rgb(523.000000 pt)=(0.596000,1.000000,0.404000);
rgb(524.000000 pt)=(0.600000,1.000000,0.400000);
rgb(525.000000 pt)=(0.604000,1.000000,0.396000);
rgb(526.000000 pt)=(0.608000,1.000000,0.392000);
rgb(527.000000 pt)=(0.612000,1.000000,0.388000);
rgb(528.000000 pt)=(0.616000,1.000000,0.384000);
rgb(529.000000 pt)=(0.620000,1.000000,0.380000);
rgb(530.000000 pt)=(0.624000,1.000000,0.376000);
rgb(531.000000 pt)=(0.628000,1.000000,0.372000);
rgb(532.000000 pt)=(0.632000,1.000000,0.368000);
rgb(533.000000 pt)=(0.636000,1.000000,0.364000);
rgb(534.000000 pt)=(0.640000,1.000000,0.360000);
rgb(535.000000 pt)=(0.644000,1.000000,0.356000);
rgb(536.000000 pt)=(0.648000,1.000000,0.352000);
rgb(537.000000 pt)=(0.652000,1.000000,0.348000);
rgb(538.000000 pt)=(0.656000,1.000000,0.344000);
rgb(539.000000 pt)=(0.660000,1.000000,0.340000);
rgb(540.000000 pt)=(0.664000,1.000000,0.336000);
rgb(541.000000 pt)=(0.668000,1.000000,0.332000);
rgb(542.000000 pt)=(0.672000,1.000000,0.328000);
rgb(543.000000 pt)=(0.676000,1.000000,0.324000);
rgb(544.000000 pt)=(0.680000,1.000000,0.320000);
rgb(545.000000 pt)=(0.684000,1.000000,0.316000);
rgb(546.000000 pt)=(0.688000,1.000000,0.312000);
rgb(547.000000 pt)=(0.692000,1.000000,0.308000);
rgb(548.000000 pt)=(0.696000,1.000000,0.304000);
rgb(549.000000 pt)=(0.700000,1.000000,0.300000);
rgb(550.000000 pt)=(0.704000,1.000000,0.296000);
rgb(551.000000 pt)=(0.708000,1.000000,0.292000);
rgb(552.000000 pt)=(0.712000,1.000000,0.288000);
rgb(553.000000 pt)=(0.716000,1.000000,0.284000);
rgb(554.000000 pt)=(0.720000,1.000000,0.280000);
rgb(555.000000 pt)=(0.724000,1.000000,0.276000);
rgb(556.000000 pt)=(0.728000,1.000000,0.272000);
rgb(557.000000 pt)=(0.732000,1.000000,0.268000);
rgb(558.000000 pt)=(0.736000,1.000000,0.264000);
rgb(559.000000 pt)=(0.740000,1.000000,0.260000);
rgb(560.000000 pt)=(0.744000,1.000000,0.256000);
rgb(561.000000 pt)=(0.748000,1.000000,0.252000);
rgb(562.000000 pt)=(0.752000,1.000000,0.248000);
rgb(563.000000 pt)=(0.756000,1.000000,0.244000);
rgb(564.000000 pt)=(0.760000,1.000000,0.240000);
rgb(565.000000 pt)=(0.764000,1.000000,0.236000);
rgb(566.000000 pt)=(0.768000,1.000000,0.232000);
rgb(567.000000 pt)=(0.772000,1.000000,0.228000);
rgb(568.000000 pt)=(0.776000,1.000000,0.224000);
rgb(569.000000 pt)=(0.780000,1.000000,0.220000);
rgb(570.000000 pt)=(0.784000,1.000000,0.216000);
rgb(571.000000 pt)=(0.788000,1.000000,0.212000);
rgb(572.000000 pt)=(0.792000,1.000000,0.208000);
rgb(573.000000 pt)=(0.796000,1.000000,0.204000);
rgb(574.000000 pt)=(0.800000,1.000000,0.200000);
rgb(575.000000 pt)=(0.804000,1.000000,0.196000);
rgb(576.000000 pt)=(0.808000,1.000000,0.192000);
rgb(577.000000 pt)=(0.812000,1.000000,0.188000);
rgb(578.000000 pt)=(0.816000,1.000000,0.184000);
rgb(579.000000 pt)=(0.820000,1.000000,0.180000);
rgb(580.000000 pt)=(0.824000,1.000000,0.176000);
rgb(581.000000 pt)=(0.828000,1.000000,0.172000);
rgb(582.000000 pt)=(0.832000,1.000000,0.168000);
rgb(583.000000 pt)=(0.836000,1.000000,0.164000);
rgb(584.000000 pt)=(0.840000,1.000000,0.160000);
rgb(585.000000 pt)=(0.844000,1.000000,0.156000);
rgb(586.000000 pt)=(0.848000,1.000000,0.152000);
rgb(587.000000 pt)=(0.852000,1.000000,0.148000);
rgb(588.000000 pt)=(0.856000,1.000000,0.144000);
rgb(589.000000 pt)=(0.860000,1.000000,0.140000);
rgb(590.000000 pt)=(0.864000,1.000000,0.136000);
rgb(591.000000 pt)=(0.868000,1.000000,0.132000);
rgb(592.000000 pt)=(0.872000,1.000000,0.128000);
rgb(593.000000 pt)=(0.876000,1.000000,0.124000);
rgb(594.000000 pt)=(0.880000,1.000000,0.120000);
rgb(595.000000 pt)=(0.884000,1.000000,0.116000);
rgb(596.000000 pt)=(0.888000,1.000000,0.112000);
rgb(597.000000 pt)=(0.892000,1.000000,0.108000);
rgb(598.000000 pt)=(0.896000,1.000000,0.104000);
rgb(599.000000 pt)=(0.900000,1.000000,0.100000);
rgb(600.000000 pt)=(0.904000,1.000000,0.096000);
rgb(601.000000 pt)=(0.908000,1.000000,0.092000);
rgb(602.000000 pt)=(0.912000,1.000000,0.088000);
rgb(603.000000 pt)=(0.916000,1.000000,0.084000);
rgb(604.000000 pt)=(0.920000,1.000000,0.080000);
rgb(605.000000 pt)=(0.924000,1.000000,0.076000);
rgb(606.000000 pt)=(0.928000,1.000000,0.072000);
rgb(607.000000 pt)=(0.932000,1.000000,0.068000);
rgb(608.000000 pt)=(0.936000,1.000000,0.064000);
rgb(609.000000 pt)=(0.940000,1.000000,0.060000);
rgb(610.000000 pt)=(0.944000,1.000000,0.056000);
rgb(611.000000 pt)=(0.948000,1.000000,0.052000);
rgb(612.000000 pt)=(0.952000,1.000000,0.048000);
rgb(613.000000 pt)=(0.956000,1.000000,0.044000);
rgb(614.000000 pt)=(0.960000,1.000000,0.040000);
rgb(615.000000 pt)=(0.964000,1.000000,0.036000);
rgb(616.000000 pt)=(0.968000,1.000000,0.032000);
rgb(617.000000 pt)=(0.972000,1.000000,0.028000);
rgb(618.000000 pt)=(0.976000,1.000000,0.024000);
rgb(619.000000 pt)=(0.980000,1.000000,0.020000);
rgb(620.000000 pt)=(0.984000,1.000000,0.016000);
rgb(621.000000 pt)=(0.988000,1.000000,0.012000);
rgb(622.000000 pt)=(0.992000,1.000000,0.008000);
rgb(623.000000 pt)=(0.996000,1.000000,0.004000);
rgb(624.000000 pt)=(1.000000,1.000000,0.000000);
rgb(625.000000 pt)=(1.000000,0.996000,0.000000);
rgb(626.000000 pt)=(1.000000,0.992000,0.000000);
rgb(627.000000 pt)=(1.000000,0.988000,0.000000);
rgb(628.000000 pt)=(1.000000,0.984000,0.000000);
rgb(629.000000 pt)=(1.000000,0.980000,0.000000);
rgb(630.000000 pt)=(1.000000,0.976000,0.000000);
rgb(631.000000 pt)=(1.000000,0.972000,0.000000);
rgb(632.000000 pt)=(1.000000,0.968000,0.000000);
rgb(633.000000 pt)=(1.000000,0.964000,0.000000);
rgb(634.000000 pt)=(1.000000,0.960000,0.000000);
rgb(635.000000 pt)=(1.000000,0.956000,0.000000);
rgb(636.000000 pt)=(1.000000,0.952000,0.000000);
rgb(637.000000 pt)=(1.000000,0.948000,0.000000);
rgb(638.000000 pt)=(1.000000,0.944000,0.000000);
rgb(639.000000 pt)=(1.000000,0.940000,0.000000);
rgb(640.000000 pt)=(1.000000,0.936000,0.000000);
rgb(641.000000 pt)=(1.000000,0.932000,0.000000);
rgb(642.000000 pt)=(1.000000,0.928000,0.000000);
rgb(643.000000 pt)=(1.000000,0.924000,0.000000);
rgb(644.000000 pt)=(1.000000,0.920000,0.000000);
rgb(645.000000 pt)=(1.000000,0.916000,0.000000);
rgb(646.000000 pt)=(1.000000,0.912000,0.000000);
rgb(647.000000 pt)=(1.000000,0.908000,0.000000);
rgb(648.000000 pt)=(1.000000,0.904000,0.000000);
rgb(649.000000 pt)=(1.000000,0.900000,0.000000);
rgb(650.000000 pt)=(1.000000,0.896000,0.000000);
rgb(651.000000 pt)=(1.000000,0.892000,0.000000);
rgb(652.000000 pt)=(1.000000,0.888000,0.000000);
rgb(653.000000 pt)=(1.000000,0.884000,0.000000);
rgb(654.000000 pt)=(1.000000,0.880000,0.000000);
rgb(655.000000 pt)=(1.000000,0.876000,0.000000);
rgb(656.000000 pt)=(1.000000,0.872000,0.000000);
rgb(657.000000 pt)=(1.000000,0.868000,0.000000);
rgb(658.000000 pt)=(1.000000,0.864000,0.000000);
rgb(659.000000 pt)=(1.000000,0.860000,0.000000);
rgb(660.000000 pt)=(1.000000,0.856000,0.000000);
rgb(661.000000 pt)=(1.000000,0.852000,0.000000);
rgb(662.000000 pt)=(1.000000,0.848000,0.000000);
rgb(663.000000 pt)=(1.000000,0.844000,0.000000);
rgb(664.000000 pt)=(1.000000,0.840000,0.000000);
rgb(665.000000 pt)=(1.000000,0.836000,0.000000);
rgb(666.000000 pt)=(1.000000,0.832000,0.000000);
rgb(667.000000 pt)=(1.000000,0.828000,0.000000);
rgb(668.000000 pt)=(1.000000,0.824000,0.000000);
rgb(669.000000 pt)=(1.000000,0.820000,0.000000);
rgb(670.000000 pt)=(1.000000,0.816000,0.000000);
rgb(671.000000 pt)=(1.000000,0.812000,0.000000);
rgb(672.000000 pt)=(1.000000,0.808000,0.000000);
rgb(673.000000 pt)=(1.000000,0.804000,0.000000);
rgb(674.000000 pt)=(1.000000,0.800000,0.000000);
rgb(675.000000 pt)=(1.000000,0.796000,0.000000);
rgb(676.000000 pt)=(1.000000,0.792000,0.000000);
rgb(677.000000 pt)=(1.000000,0.788000,0.000000);
rgb(678.000000 pt)=(1.000000,0.784000,0.000000);
rgb(679.000000 pt)=(1.000000,0.780000,0.000000);
rgb(680.000000 pt)=(1.000000,0.776000,0.000000);
rgb(681.000000 pt)=(1.000000,0.772000,0.000000);
rgb(682.000000 pt)=(1.000000,0.768000,0.000000);
rgb(683.000000 pt)=(1.000000,0.764000,0.000000);
rgb(684.000000 pt)=(1.000000,0.760000,0.000000);
rgb(685.000000 pt)=(1.000000,0.756000,0.000000);
rgb(686.000000 pt)=(1.000000,0.752000,0.000000);
rgb(687.000000 pt)=(1.000000,0.748000,0.000000);
rgb(688.000000 pt)=(1.000000,0.744000,0.000000);
rgb(689.000000 pt)=(1.000000,0.740000,0.000000);
rgb(690.000000 pt)=(1.000000,0.736000,0.000000);
rgb(691.000000 pt)=(1.000000,0.732000,0.000000);
rgb(692.000000 pt)=(1.000000,0.728000,0.000000);
rgb(693.000000 pt)=(1.000000,0.724000,0.000000);
rgb(694.000000 pt)=(1.000000,0.720000,0.000000);
rgb(695.000000 pt)=(1.000000,0.716000,0.000000);
rgb(696.000000 pt)=(1.000000,0.712000,0.000000);
rgb(697.000000 pt)=(1.000000,0.708000,0.000000);
rgb(698.000000 pt)=(1.000000,0.704000,0.000000);
rgb(699.000000 pt)=(1.000000,0.700000,0.000000);
rgb(700.000000 pt)=(1.000000,0.696000,0.000000);
rgb(701.000000 pt)=(1.000000,0.692000,0.000000);
rgb(702.000000 pt)=(1.000000,0.688000,0.000000);
rgb(703.000000 pt)=(1.000000,0.684000,0.000000);
rgb(704.000000 pt)=(1.000000,0.680000,0.000000);
rgb(705.000000 pt)=(1.000000,0.676000,0.000000);
rgb(706.000000 pt)=(1.000000,0.672000,0.000000);
rgb(707.000000 pt)=(1.000000,0.668000,0.000000);
rgb(708.000000 pt)=(1.000000,0.664000,0.000000);
rgb(709.000000 pt)=(1.000000,0.660000,0.000000);
rgb(710.000000 pt)=(1.000000,0.656000,0.000000);
rgb(711.000000 pt)=(1.000000,0.652000,0.000000);
rgb(712.000000 pt)=(1.000000,0.648000,0.000000);
rgb(713.000000 pt)=(1.000000,0.644000,0.000000);
rgb(714.000000 pt)=(1.000000,0.640000,0.000000);
rgb(715.000000 pt)=(1.000000,0.636000,0.000000);
rgb(716.000000 pt)=(1.000000,0.632000,0.000000);
rgb(717.000000 pt)=(1.000000,0.628000,0.000000);
rgb(718.000000 pt)=(1.000000,0.624000,0.000000);
rgb(719.000000 pt)=(1.000000,0.620000,0.000000);
rgb(720.000000 pt)=(1.000000,0.616000,0.000000);
rgb(721.000000 pt)=(1.000000,0.612000,0.000000);
rgb(722.000000 pt)=(1.000000,0.608000,0.000000);
rgb(723.000000 pt)=(1.000000,0.604000,0.000000);
rgb(724.000000 pt)=(1.000000,0.600000,0.000000);
rgb(725.000000 pt)=(1.000000,0.596000,0.000000);
rgb(726.000000 pt)=(1.000000,0.592000,0.000000);
rgb(727.000000 pt)=(1.000000,0.588000,0.000000);
rgb(728.000000 pt)=(1.000000,0.584000,0.000000);
rgb(729.000000 pt)=(1.000000,0.580000,0.000000);
rgb(730.000000 pt)=(1.000000,0.576000,0.000000);
rgb(731.000000 pt)=(1.000000,0.572000,0.000000);
rgb(732.000000 pt)=(1.000000,0.568000,0.000000);
rgb(733.000000 pt)=(1.000000,0.564000,0.000000);
rgb(734.000000 pt)=(1.000000,0.560000,0.000000);
rgb(735.000000 pt)=(1.000000,0.556000,0.000000);
rgb(736.000000 pt)=(1.000000,0.552000,0.000000);
rgb(737.000000 pt)=(1.000000,0.548000,0.000000);
rgb(738.000000 pt)=(1.000000,0.544000,0.000000);
rgb(739.000000 pt)=(1.000000,0.540000,0.000000);
rgb(740.000000 pt)=(1.000000,0.536000,0.000000);
rgb(741.000000 pt)=(1.000000,0.532000,0.000000);
rgb(742.000000 pt)=(1.000000,0.528000,0.000000);
rgb(743.000000 pt)=(1.000000,0.524000,0.000000);
rgb(744.000000 pt)=(1.000000,0.520000,0.000000);
rgb(745.000000 pt)=(1.000000,0.516000,0.000000);
rgb(746.000000 pt)=(1.000000,0.512000,0.000000);
rgb(747.000000 pt)=(1.000000,0.508000,0.000000);
rgb(748.000000 pt)=(1.000000,0.504000,0.000000);
rgb(749.000000 pt)=(1.000000,0.500000,0.000000);
rgb(750.000000 pt)=(1.000000,0.496000,0.000000);
rgb(751.000000 pt)=(1.000000,0.492000,0.000000);
rgb(752.000000 pt)=(1.000000,0.488000,0.000000);
rgb(753.000000 pt)=(1.000000,0.484000,0.000000);
rgb(754.000000 pt)=(1.000000,0.480000,0.000000);
rgb(755.000000 pt)=(1.000000,0.476000,0.000000);
rgb(756.000000 pt)=(1.000000,0.472000,0.000000);
rgb(757.000000 pt)=(1.000000,0.468000,0.000000);
rgb(758.000000 pt)=(1.000000,0.464000,0.000000);
rgb(759.000000 pt)=(1.000000,0.460000,0.000000);
rgb(760.000000 pt)=(1.000000,0.456000,0.000000);
rgb(761.000000 pt)=(1.000000,0.452000,0.000000);
rgb(762.000000 pt)=(1.000000,0.448000,0.000000);
rgb(763.000000 pt)=(1.000000,0.444000,0.000000);
rgb(764.000000 pt)=(1.000000,0.440000,0.000000);
rgb(765.000000 pt)=(1.000000,0.436000,0.000000);
rgb(766.000000 pt)=(1.000000,0.432000,0.000000);
rgb(767.000000 pt)=(1.000000,0.428000,0.000000);
rgb(768.000000 pt)=(1.000000,0.424000,0.000000);
rgb(769.000000 pt)=(1.000000,0.420000,0.000000);
rgb(770.000000 pt)=(1.000000,0.416000,0.000000);
rgb(771.000000 pt)=(1.000000,0.412000,0.000000);
rgb(772.000000 pt)=(1.000000,0.408000,0.000000);
rgb(773.000000 pt)=(1.000000,0.404000,0.000000);
rgb(774.000000 pt)=(1.000000,0.400000,0.000000);
rgb(775.000000 pt)=(1.000000,0.396000,0.000000);
rgb(776.000000 pt)=(1.000000,0.392000,0.000000);
rgb(777.000000 pt)=(1.000000,0.388000,0.000000);
rgb(778.000000 pt)=(1.000000,0.384000,0.000000);
rgb(779.000000 pt)=(1.000000,0.380000,0.000000);
rgb(780.000000 pt)=(1.000000,0.376000,0.000000);
rgb(781.000000 pt)=(1.000000,0.372000,0.000000);
rgb(782.000000 pt)=(1.000000,0.368000,0.000000);
rgb(783.000000 pt)=(1.000000,0.364000,0.000000);
rgb(784.000000 pt)=(1.000000,0.360000,0.000000);
rgb(785.000000 pt)=(1.000000,0.356000,0.000000);
rgb(786.000000 pt)=(1.000000,0.352000,0.000000);
rgb(787.000000 pt)=(1.000000,0.348000,0.000000);
rgb(788.000000 pt)=(1.000000,0.344000,0.000000);
rgb(789.000000 pt)=(1.000000,0.340000,0.000000);
rgb(790.000000 pt)=(1.000000,0.336000,0.000000);
rgb(791.000000 pt)=(1.000000,0.332000,0.000000);
rgb(792.000000 pt)=(1.000000,0.328000,0.000000);
rgb(793.000000 pt)=(1.000000,0.324000,0.000000);
rgb(794.000000 pt)=(1.000000,0.320000,0.000000);
rgb(795.000000 pt)=(1.000000,0.316000,0.000000);
rgb(796.000000 pt)=(1.000000,0.312000,0.000000);
rgb(797.000000 pt)=(1.000000,0.308000,0.000000);
rgb(798.000000 pt)=(1.000000,0.304000,0.000000);
rgb(799.000000 pt)=(1.000000,0.300000,0.000000);
rgb(800.000000 pt)=(1.000000,0.296000,0.000000);
rgb(801.000000 pt)=(1.000000,0.292000,0.000000);
rgb(802.000000 pt)=(1.000000,0.288000,0.000000);
rgb(803.000000 pt)=(1.000000,0.284000,0.000000);
rgb(804.000000 pt)=(1.000000,0.280000,0.000000);
rgb(805.000000 pt)=(1.000000,0.276000,0.000000);
rgb(806.000000 pt)=(1.000000,0.272000,0.000000);
rgb(807.000000 pt)=(1.000000,0.268000,0.000000);
rgb(808.000000 pt)=(1.000000,0.264000,0.000000);
rgb(809.000000 pt)=(1.000000,0.260000,0.000000);
rgb(810.000000 pt)=(1.000000,0.256000,0.000000);
rgb(811.000000 pt)=(1.000000,0.252000,0.000000);
rgb(812.000000 pt)=(1.000000,0.248000,0.000000);
rgb(813.000000 pt)=(1.000000,0.244000,0.000000);
rgb(814.000000 pt)=(1.000000,0.240000,0.000000);
rgb(815.000000 pt)=(1.000000,0.236000,0.000000);
rgb(816.000000 pt)=(1.000000,0.232000,0.000000);
rgb(817.000000 pt)=(1.000000,0.228000,0.000000);
rgb(818.000000 pt)=(1.000000,0.224000,0.000000);
rgb(819.000000 pt)=(1.000000,0.220000,0.000000);
rgb(820.000000 pt)=(1.000000,0.216000,0.000000);
rgb(821.000000 pt)=(1.000000,0.212000,0.000000);
rgb(822.000000 pt)=(1.000000,0.208000,0.000000);
rgb(823.000000 pt)=(1.000000,0.204000,0.000000);
rgb(824.000000 pt)=(1.000000,0.200000,0.000000);
rgb(825.000000 pt)=(1.000000,0.196000,0.000000);
rgb(826.000000 pt)=(1.000000,0.192000,0.000000);
rgb(827.000000 pt)=(1.000000,0.188000,0.000000);
rgb(828.000000 pt)=(1.000000,0.184000,0.000000);
rgb(829.000000 pt)=(1.000000,0.180000,0.000000);
rgb(830.000000 pt)=(1.000000,0.176000,0.000000);
rgb(831.000000 pt)=(1.000000,0.172000,0.000000);
rgb(832.000000 pt)=(1.000000,0.168000,0.000000);
rgb(833.000000 pt)=(1.000000,0.164000,0.000000);
rgb(834.000000 pt)=(1.000000,0.160000,0.000000);
rgb(835.000000 pt)=(1.000000,0.156000,0.000000);
rgb(836.000000 pt)=(1.000000,0.152000,0.000000);
rgb(837.000000 pt)=(1.000000,0.148000,0.000000);
rgb(838.000000 pt)=(1.000000,0.144000,0.000000);
rgb(839.000000 pt)=(1.000000,0.140000,0.000000);
rgb(840.000000 pt)=(1.000000,0.136000,0.000000);
rgb(841.000000 pt)=(1.000000,0.132000,0.000000);
rgb(842.000000 pt)=(1.000000,0.128000,0.000000);
rgb(843.000000 pt)=(1.000000,0.124000,0.000000);
rgb(844.000000 pt)=(1.000000,0.120000,0.000000);
rgb(845.000000 pt)=(1.000000,0.116000,0.000000);
rgb(846.000000 pt)=(1.000000,0.112000,0.000000);
rgb(847.000000 pt)=(1.000000,0.108000,0.000000);
rgb(848.000000 pt)=(1.000000,0.104000,0.000000);
rgb(849.000000 pt)=(1.000000,0.100000,0.000000);
rgb(850.000000 pt)=(1.000000,0.096000,0.000000);
rgb(851.000000 pt)=(1.000000,0.092000,0.000000);
rgb(852.000000 pt)=(1.000000,0.088000,0.000000);
rgb(853.000000 pt)=(1.000000,0.084000,0.000000);
rgb(854.000000 pt)=(1.000000,0.080000,0.000000);
rgb(855.000000 pt)=(1.000000,0.076000,0.000000);
rgb(856.000000 pt)=(1.000000,0.072000,0.000000);
rgb(857.000000 pt)=(1.000000,0.068000,0.000000);
rgb(858.000000 pt)=(1.000000,0.064000,0.000000);
rgb(859.000000 pt)=(1.000000,0.060000,0.000000);
rgb(860.000000 pt)=(1.000000,0.056000,0.000000);
rgb(861.000000 pt)=(1.000000,0.052000,0.000000);
rgb(862.000000 pt)=(1.000000,0.048000,0.000000);
rgb(863.000000 pt)=(1.000000,0.044000,0.000000);
rgb(864.000000 pt)=(1.000000,0.040000,0.000000);
rgb(865.000000 pt)=(1.000000,0.036000,0.000000);
rgb(866.000000 pt)=(1.000000,0.032000,0.000000);
rgb(867.000000 pt)=(1.000000,0.028000,0.000000);
rgb(868.000000 pt)=(1.000000,0.024000,0.000000);
rgb(869.000000 pt)=(1.000000,0.020000,0.000000);
rgb(870.000000 pt)=(1.000000,0.016000,0.000000);
rgb(871.000000 pt)=(1.000000,0.012000,0.000000);
rgb(872.000000 pt)=(1.000000,0.008000,0.000000);
rgb(873.000000 pt)=(1.000000,0.004000,0.000000);
rgb(874.000000 pt)=(1.000000,0.000000,0.000000);
rgb(875.000000 pt)=(0.996000,0.000000,0.000000);
rgb(876.000000 pt)=(0.992000,0.000000,0.000000);
rgb(877.000000 pt)=(0.988000,0.000000,0.000000);
rgb(878.000000 pt)=(0.984000,0.000000,0.000000);
rgb(879.000000 pt)=(0.980000,0.000000,0.000000);
rgb(880.000000 pt)=(0.976000,0.000000,0.000000);
rgb(881.000000 pt)=(0.972000,0.000000,0.000000);
rgb(882.000000 pt)=(0.968000,0.000000,0.000000);
rgb(883.000000 pt)=(0.964000,0.000000,0.000000);
rgb(884.000000 pt)=(0.960000,0.000000,0.000000);
rgb(885.000000 pt)=(0.956000,0.000000,0.000000);
rgb(886.000000 pt)=(0.952000,0.000000,0.000000);
rgb(887.000000 pt)=(0.948000,0.000000,0.000000);
rgb(888.000000 pt)=(0.944000,0.000000,0.000000);
rgb(889.000000 pt)=(0.940000,0.000000,0.000000);
rgb(890.000000 pt)=(0.936000,0.000000,0.000000);
rgb(891.000000 pt)=(0.932000,0.000000,0.000000);
rgb(892.000000 pt)=(0.928000,0.000000,0.000000);
rgb(893.000000 pt)=(0.924000,0.000000,0.000000);
rgb(894.000000 pt)=(0.920000,0.000000,0.000000);
rgb(895.000000 pt)=(0.916000,0.000000,0.000000);
rgb(896.000000 pt)=(0.912000,0.000000,0.000000);
rgb(897.000000 pt)=(0.908000,0.000000,0.000000);
rgb(898.000000 pt)=(0.904000,0.000000,0.000000);
rgb(899.000000 pt)=(0.900000,0.000000,0.000000);
rgb(900.000000 pt)=(0.896000,0.000000,0.000000);
rgb(901.000000 pt)=(0.892000,0.000000,0.000000);
rgb(902.000000 pt)=(0.888000,0.000000,0.000000);
rgb(903.000000 pt)=(0.884000,0.000000,0.000000);
rgb(904.000000 pt)=(0.880000,0.000000,0.000000);
rgb(905.000000 pt)=(0.876000,0.000000,0.000000);
rgb(906.000000 pt)=(0.872000,0.000000,0.000000);
rgb(907.000000 pt)=(0.868000,0.000000,0.000000);
rgb(908.000000 pt)=(0.864000,0.000000,0.000000);
rgb(909.000000 pt)=(0.860000,0.000000,0.000000);
rgb(910.000000 pt)=(0.856000,0.000000,0.000000);
rgb(911.000000 pt)=(0.852000,0.000000,0.000000);
rgb(912.000000 pt)=(0.848000,0.000000,0.000000);
rgb(913.000000 pt)=(0.844000,0.000000,0.000000);
rgb(914.000000 pt)=(0.840000,0.000000,0.000000);
rgb(915.000000 pt)=(0.836000,0.000000,0.000000);
rgb(916.000000 pt)=(0.832000,0.000000,0.000000);
rgb(917.000000 pt)=(0.828000,0.000000,0.000000);
rgb(918.000000 pt)=(0.824000,0.000000,0.000000);
rgb(919.000000 pt)=(0.820000,0.000000,0.000000);
rgb(920.000000 pt)=(0.816000,0.000000,0.000000);
rgb(921.000000 pt)=(0.812000,0.000000,0.000000);
rgb(922.000000 pt)=(0.808000,0.000000,0.000000);
rgb(923.000000 pt)=(0.804000,0.000000,0.000000);
rgb(924.000000 pt)=(0.800000,0.000000,0.000000);
rgb(925.000000 pt)=(0.796000,0.000000,0.000000);
rgb(926.000000 pt)=(0.792000,0.000000,0.000000);
rgb(927.000000 pt)=(0.788000,0.000000,0.000000);
rgb(928.000000 pt)=(0.784000,0.000000,0.000000);
rgb(929.000000 pt)=(0.780000,0.000000,0.000000);
rgb(930.000000 pt)=(0.776000,0.000000,0.000000);
rgb(931.000000 pt)=(0.772000,0.000000,0.000000);
rgb(932.000000 pt)=(0.768000,0.000000,0.000000);
rgb(933.000000 pt)=(0.764000,0.000000,0.000000);
rgb(934.000000 pt)=(0.760000,0.000000,0.000000);
rgb(935.000000 pt)=(0.756000,0.000000,0.000000);
rgb(936.000000 pt)=(0.752000,0.000000,0.000000);
rgb(937.000000 pt)=(0.748000,0.000000,0.000000);
rgb(938.000000 pt)=(0.744000,0.000000,0.000000);
rgb(939.000000 pt)=(0.740000,0.000000,0.000000);
rgb(940.000000 pt)=(0.736000,0.000000,0.000000);
rgb(941.000000 pt)=(0.732000,0.000000,0.000000);
rgb(942.000000 pt)=(0.728000,0.000000,0.000000);
rgb(943.000000 pt)=(0.724000,0.000000,0.000000);
rgb(944.000000 pt)=(0.720000,0.000000,0.000000);
rgb(945.000000 pt)=(0.716000,0.000000,0.000000);
rgb(946.000000 pt)=(0.712000,0.000000,0.000000);
rgb(947.000000 pt)=(0.708000,0.000000,0.000000);
rgb(948.000000 pt)=(0.704000,0.000000,0.000000);
rgb(949.000000 pt)=(0.700000,0.000000,0.000000);
rgb(950.000000 pt)=(0.696000,0.000000,0.000000);
rgb(951.000000 pt)=(0.692000,0.000000,0.000000);
rgb(952.000000 pt)=(0.688000,0.000000,0.000000);
rgb(953.000000 pt)=(0.684000,0.000000,0.000000);
rgb(954.000000 pt)=(0.680000,0.000000,0.000000);
rgb(955.000000 pt)=(0.676000,0.000000,0.000000);
rgb(956.000000 pt)=(0.672000,0.000000,0.000000);
rgb(957.000000 pt)=(0.668000,0.000000,0.000000);
rgb(958.000000 pt)=(0.664000,0.000000,0.000000);
rgb(959.000000 pt)=(0.660000,0.000000,0.000000);
rgb(960.000000 pt)=(0.656000,0.000000,0.000000);
rgb(961.000000 pt)=(0.652000,0.000000,0.000000);
rgb(962.000000 pt)=(0.648000,0.000000,0.000000);
rgb(963.000000 pt)=(0.644000,0.000000,0.000000);
rgb(964.000000 pt)=(0.640000,0.000000,0.000000);
rgb(965.000000 pt)=(0.636000,0.000000,0.000000);
rgb(966.000000 pt)=(0.632000,0.000000,0.000000);
rgb(967.000000 pt)=(0.628000,0.000000,0.000000);
rgb(968.000000 pt)=(0.624000,0.000000,0.000000);
rgb(969.000000 pt)=(0.620000,0.000000,0.000000);
rgb(970.000000 pt)=(0.616000,0.000000,0.000000);
rgb(971.000000 pt)=(0.612000,0.000000,0.000000);
rgb(972.000000 pt)=(0.608000,0.000000,0.000000);
rgb(973.000000 pt)=(0.604000,0.000000,0.000000);
rgb(974.000000 pt)=(0.600000,0.000000,0.000000);
rgb(975.000000 pt)=(0.596000,0.000000,0.000000);
rgb(976.000000 pt)=(0.592000,0.000000,0.000000);
rgb(977.000000 pt)=(0.588000,0.000000,0.000000);
rgb(978.000000 pt)=(0.584000,0.000000,0.000000);
rgb(979.000000 pt)=(0.580000,0.000000,0.000000);
rgb(980.000000 pt)=(0.576000,0.000000,0.000000);
rgb(981.000000 pt)=(0.572000,0.000000,0.000000);
rgb(982.000000 pt)=(0.568000,0.000000,0.000000);
rgb(983.000000 pt)=(0.564000,0.000000,0.000000);
rgb(984.000000 pt)=(0.560000,0.000000,0.000000);
rgb(985.000000 pt)=(0.556000,0.000000,0.000000);
rgb(986.000000 pt)=(0.552000,0.000000,0.000000);
rgb(987.000000 pt)=(0.548000,0.000000,0.000000);
rgb(988.000000 pt)=(0.544000,0.000000,0.000000);
rgb(989.000000 pt)=(0.540000,0.000000,0.000000);
rgb(990.000000 pt)=(0.536000,0.000000,0.000000);
rgb(991.000000 pt)=(0.532000,0.000000,0.000000);
rgb(992.000000 pt)=(0.528000,0.000000,0.000000);
rgb(993.000000 pt)=(0.524000,0.000000,0.000000);
rgb(994.000000 pt)=(0.520000,0.000000,0.000000);
rgb(995.000000 pt)=(0.516000,0.000000,0.000000);
rgb(996.000000 pt)=(0.512000,0.000000,0.000000);
rgb(997.000000 pt)=(0.508000,0.000000,0.000000);
rgb(998.000000 pt)=(0.504000,0.000000,0.000000);
rgb(999.000000 pt)=(0.500000,0.000000,0.000000);
}}

\pgfplotsset{
	colormap={myhot}{
		rgb(0pt)=(1.000000,1.000000,1.000000);
		rgb(1pt)=(1.000000,1.000000,0.996000);
		rgb(2pt)=(1.000000,1.000000,0.992000);
		rgb(3pt)=(1.000000,1.000000,0.988000);
		rgb(4pt)=(1.000000,1.000000,0.984000);
		rgb(5pt)=(1.000000,1.000000,0.980000);
		rgb(6pt)=(1.000000,1.000000,0.976000);
		rgb(7pt)=(1.000000,1.000000,0.972000);
		rgb(8pt)=(1.000000,1.000000,0.968000);
		rgb(9pt)=(1.000000,1.000000,0.964000);
		rgb(10pt)=(1.000000,1.000000,0.960000);
		rgb(11pt)=(1.000000,1.000000,0.956000);
		rgb(12pt)=(1.000000,1.000000,0.952000);
		rgb(13pt)=(1.000000,1.000000,0.948000);
		rgb(14pt)=(1.000000,1.000000,0.944000);
		rgb(15pt)=(1.000000,1.000000,0.940000);
		rgb(16pt)=(1.000000,1.000000,0.936000);
		rgb(17pt)=(1.000000,1.000000,0.932000);
		rgb(18pt)=(1.000000,1.000000,0.928000);
		rgb(19pt)=(1.000000,1.000000,0.924000);
		rgb(20pt)=(1.000000,1.000000,0.920000);
		rgb(21pt)=(1.000000,1.000000,0.916000);
		rgb(22pt)=(1.000000,1.000000,0.912000);
		rgb(23pt)=(1.000000,1.000000,0.908000);
		rgb(24pt)=(1.000000,1.000000,0.904000);
		rgb(25pt)=(1.000000,1.000000,0.900000);
		rgb(26pt)=(1.000000,1.000000,0.896000);
		rgb(27pt)=(1.000000,1.000000,0.892000);
		rgb(28pt)=(1.000000,1.000000,0.888000);
		rgb(29pt)=(1.000000,1.000000,0.884000);
		rgb(30pt)=(1.000000,1.000000,0.880000);
		rgb(31pt)=(1.000000,1.000000,0.876000);
		rgb(32pt)=(1.000000,1.000000,0.872000);
		rgb(33pt)=(1.000000,1.000000,0.868000);
		rgb(34pt)=(1.000000,1.000000,0.864000);
		rgb(35pt)=(1.000000,1.000000,0.860000);
		rgb(36pt)=(1.000000,1.000000,0.856000);
		rgb(37pt)=(1.000000,1.000000,0.852000);
		rgb(38pt)=(1.000000,1.000000,0.848000);
		rgb(39pt)=(1.000000,1.000000,0.844000);
		rgb(40pt)=(1.000000,1.000000,0.840000);
		rgb(41pt)=(1.000000,1.000000,0.836000);
		rgb(42pt)=(1.000000,1.000000,0.832000);
		rgb(43pt)=(1.000000,1.000000,0.828000);
		rgb(44pt)=(1.000000,1.000000,0.824000);
		rgb(45pt)=(1.000000,1.000000,0.820000);
		rgb(46pt)=(1.000000,1.000000,0.816000);
		rgb(47pt)=(1.000000,1.000000,0.812000);
		rgb(48pt)=(1.000000,1.000000,0.808000);
		rgb(49pt)=(1.000000,1.000000,0.804000);
		rgb(50pt)=(1.000000,1.000000,0.800000);
		rgb(51pt)=(1.000000,1.000000,0.796000);
		rgb(52pt)=(1.000000,1.000000,0.792000);
		rgb(53pt)=(1.000000,1.000000,0.788000);
		rgb(54pt)=(1.000000,1.000000,0.784000);
		rgb(55pt)=(1.000000,1.000000,0.780000);
		rgb(56pt)=(1.000000,1.000000,0.776000);
		rgb(57pt)=(1.000000,1.000000,0.772000);
		rgb(58pt)=(1.000000,1.000000,0.768000);
		rgb(59pt)=(1.000000,1.000000,0.764000);
		rgb(60pt)=(1.000000,1.000000,0.760000);
		rgb(61pt)=(1.000000,1.000000,0.756000);
		rgb(62pt)=(1.000000,1.000000,0.752000);
		rgb(63pt)=(1.000000,1.000000,0.748000);
		rgb(64pt)=(1.000000,1.000000,0.744000);
		rgb(65pt)=(1.000000,1.000000,0.740000);
		rgb(66pt)=(1.000000,1.000000,0.736000);
		rgb(67pt)=(1.000000,1.000000,0.732000);
		rgb(68pt)=(1.000000,1.000000,0.728000);
		rgb(69pt)=(1.000000,1.000000,0.724000);
		rgb(70pt)=(1.000000,1.000000,0.720000);
		rgb(71pt)=(1.000000,1.000000,0.716000);
		rgb(72pt)=(1.000000,1.000000,0.712000);
		rgb(73pt)=(1.000000,1.000000,0.708000);
		rgb(74pt)=(1.000000,1.000000,0.704000);
		rgb(75pt)=(1.000000,1.000000,0.700000);
		rgb(76pt)=(1.000000,1.000000,0.696000);
		rgb(77pt)=(1.000000,1.000000,0.692000);
		rgb(78pt)=(1.000000,1.000000,0.688000);
		rgb(79pt)=(1.000000,1.000000,0.684000);
		rgb(80pt)=(1.000000,1.000000,0.680000);
		rgb(81pt)=(1.000000,1.000000,0.676000);
		rgb(82pt)=(1.000000,1.000000,0.672000);
		rgb(83pt)=(1.000000,1.000000,0.668000);
		rgb(84pt)=(1.000000,1.000000,0.664000);
		rgb(85pt)=(1.000000,1.000000,0.660000);
		rgb(86pt)=(1.000000,1.000000,0.656000);
		rgb(87pt)=(1.000000,1.000000,0.652000);
		rgb(88pt)=(1.000000,1.000000,0.648000);
		rgb(89pt)=(1.000000,1.000000,0.644000);
		rgb(90pt)=(1.000000,1.000000,0.640000);
		rgb(91pt)=(1.000000,1.000000,0.636000);
		rgb(92pt)=(1.000000,1.000000,0.632000);
		rgb(93pt)=(1.000000,1.000000,0.628000);
		rgb(94pt)=(1.000000,1.000000,0.624000);
		rgb(95pt)=(1.000000,1.000000,0.620000);
		rgb(96pt)=(1.000000,1.000000,0.616000);
		rgb(97pt)=(1.000000,1.000000,0.612000);
		rgb(98pt)=(1.000000,1.000000,0.608000);
		rgb(99pt)=(1.000000,1.000000,0.604000);
		rgb(100pt)=(1.000000,1.000000,0.600000);
		rgb(101pt)=(1.000000,1.000000,0.596000);
		rgb(102pt)=(1.000000,1.000000,0.592000);
		rgb(103pt)=(1.000000,1.000000,0.588000);
		rgb(104pt)=(1.000000,1.000000,0.584000);
		rgb(105pt)=(1.000000,1.000000,0.580000);
		rgb(106pt)=(1.000000,1.000000,0.576000);
		rgb(107pt)=(1.000000,1.000000,0.572000);
		rgb(108pt)=(1.000000,1.000000,0.568000);
		rgb(109pt)=(1.000000,1.000000,0.564000);
		rgb(110pt)=(1.000000,1.000000,0.560000);
		rgb(111pt)=(1.000000,1.000000,0.556000);
		rgb(112pt)=(1.000000,1.000000,0.552000);
		rgb(113pt)=(1.000000,1.000000,0.548000);
		rgb(114pt)=(1.000000,1.000000,0.544000);
		rgb(115pt)=(1.000000,1.000000,0.540000);
		rgb(116pt)=(1.000000,1.000000,0.536000);
		rgb(117pt)=(1.000000,1.000000,0.532000);
		rgb(118pt)=(1.000000,1.000000,0.528000);
		rgb(119pt)=(1.000000,1.000000,0.524000);
		rgb(120pt)=(1.000000,1.000000,0.520000);
		rgb(121pt)=(1.000000,1.000000,0.516000);
		rgb(122pt)=(1.000000,1.000000,0.512000);
		rgb(123pt)=(1.000000,1.000000,0.508000);
		rgb(124pt)=(1.000000,1.000000,0.504000);
		rgb(125pt)=(1.000000,1.000000,0.500000);
		rgb(126pt)=(1.000000,1.000000,0.496000);
		rgb(127pt)=(1.000000,1.000000,0.492000);
		rgb(128pt)=(1.000000,1.000000,0.488000);
		rgb(129pt)=(1.000000,1.000000,0.484000);
		rgb(130pt)=(1.000000,1.000000,0.480000);
		rgb(131pt)=(1.000000,1.000000,0.476000);
		rgb(132pt)=(1.000000,1.000000,0.472000);
		rgb(133pt)=(1.000000,1.000000,0.468000);
		rgb(134pt)=(1.000000,1.000000,0.464000);
		rgb(135pt)=(1.000000,1.000000,0.460000);
		rgb(136pt)=(1.000000,1.000000,0.456000);
		rgb(137pt)=(1.000000,1.000000,0.452000);
		rgb(138pt)=(1.000000,1.000000,0.448000);
		rgb(139pt)=(1.000000,1.000000,0.444000);
		rgb(140pt)=(1.000000,1.000000,0.440000);
		rgb(141pt)=(1.000000,1.000000,0.436000);
		rgb(142pt)=(1.000000,1.000000,0.432000);
		rgb(143pt)=(1.000000,1.000000,0.428000);
		rgb(144pt)=(1.000000,1.000000,0.424000);
		rgb(145pt)=(1.000000,1.000000,0.420000);
		rgb(146pt)=(1.000000,1.000000,0.416000);
		rgb(147pt)=(1.000000,1.000000,0.412000);
		rgb(148pt)=(1.000000,1.000000,0.408000);
		rgb(149pt)=(1.000000,1.000000,0.404000);
		rgb(150pt)=(1.000000,1.000000,0.400000);
		rgb(151pt)=(1.000000,1.000000,0.396000);
		rgb(152pt)=(1.000000,1.000000,0.392000);
		rgb(153pt)=(1.000000,1.000000,0.388000);
		rgb(154pt)=(1.000000,1.000000,0.384000);
		rgb(155pt)=(1.000000,1.000000,0.380000);
		rgb(156pt)=(1.000000,1.000000,0.376000);
		rgb(157pt)=(1.000000,1.000000,0.372000);
		rgb(158pt)=(1.000000,1.000000,0.368000);
		rgb(159pt)=(1.000000,1.000000,0.364000);
		rgb(160pt)=(1.000000,1.000000,0.360000);
		rgb(161pt)=(1.000000,1.000000,0.356000);
		rgb(162pt)=(1.000000,1.000000,0.352000);
		rgb(163pt)=(1.000000,1.000000,0.348000);
		rgb(164pt)=(1.000000,1.000000,0.344000);
		rgb(165pt)=(1.000000,1.000000,0.340000);
		rgb(166pt)=(1.000000,1.000000,0.336000);
		rgb(167pt)=(1.000000,1.000000,0.332000);
		rgb(168pt)=(1.000000,1.000000,0.328000);
		rgb(169pt)=(1.000000,1.000000,0.324000);
		rgb(170pt)=(1.000000,1.000000,0.320000);
		rgb(171pt)=(1.000000,1.000000,0.316000);
		rgb(172pt)=(1.000000,1.000000,0.312000);
		rgb(173pt)=(1.000000,1.000000,0.308000);
		rgb(174pt)=(1.000000,1.000000,0.304000);
		rgb(175pt)=(1.000000,1.000000,0.300000);
		rgb(176pt)=(1.000000,1.000000,0.296000);
		rgb(177pt)=(1.000000,1.000000,0.292000);
		rgb(178pt)=(1.000000,1.000000,0.288000);
		rgb(179pt)=(1.000000,1.000000,0.284000);
		rgb(180pt)=(1.000000,1.000000,0.280000);
		rgb(181pt)=(1.000000,1.000000,0.276000);
		rgb(182pt)=(1.000000,1.000000,0.272000);
		rgb(183pt)=(1.000000,1.000000,0.268000);
		rgb(184pt)=(1.000000,1.000000,0.264000);
		rgb(185pt)=(1.000000,1.000000,0.260000);
		rgb(186pt)=(1.000000,1.000000,0.256000);
		rgb(187pt)=(1.000000,1.000000,0.252000);
		rgb(188pt)=(1.000000,1.000000,0.248000);
		rgb(189pt)=(1.000000,1.000000,0.244000);
		rgb(190pt)=(1.000000,1.000000,0.240000);
		rgb(191pt)=(1.000000,1.000000,0.236000);
		rgb(192pt)=(1.000000,1.000000,0.232000);
		rgb(193pt)=(1.000000,1.000000,0.228000);
		rgb(194pt)=(1.000000,1.000000,0.224000);
		rgb(195pt)=(1.000000,1.000000,0.220000);
		rgb(196pt)=(1.000000,1.000000,0.216000);
		rgb(197pt)=(1.000000,1.000000,0.212000);
		rgb(198pt)=(1.000000,1.000000,0.208000);
		rgb(199pt)=(1.000000,1.000000,0.204000);
		rgb(200pt)=(1.000000,1.000000,0.200000);
		rgb(201pt)=(1.000000,1.000000,0.196000);
		rgb(202pt)=(1.000000,1.000000,0.192000);
		rgb(203pt)=(1.000000,1.000000,0.188000);
		rgb(204pt)=(1.000000,1.000000,0.184000);
		rgb(205pt)=(1.000000,1.000000,0.180000);
		rgb(206pt)=(1.000000,1.000000,0.176000);
		rgb(207pt)=(1.000000,1.000000,0.172000);
		rgb(208pt)=(1.000000,1.000000,0.168000);
		rgb(209pt)=(1.000000,1.000000,0.164000);
		rgb(210pt)=(1.000000,1.000000,0.160000);
		rgb(211pt)=(1.000000,1.000000,0.156000);
		rgb(212pt)=(1.000000,1.000000,0.152000);
		rgb(213pt)=(1.000000,1.000000,0.148000);
		rgb(214pt)=(1.000000,1.000000,0.144000);
		rgb(215pt)=(1.000000,1.000000,0.140000);
		rgb(216pt)=(1.000000,1.000000,0.136000);
		rgb(217pt)=(1.000000,1.000000,0.132000);
		rgb(218pt)=(1.000000,1.000000,0.128000);
		rgb(219pt)=(1.000000,1.000000,0.124000);
		rgb(220pt)=(1.000000,1.000000,0.120000);
		rgb(221pt)=(1.000000,1.000000,0.116000);
		rgb(222pt)=(1.000000,1.000000,0.112000);
		rgb(223pt)=(1.000000,1.000000,0.108000);
		rgb(224pt)=(1.000000,1.000000,0.104000);
		rgb(225pt)=(1.000000,1.000000,0.100000);
		rgb(226pt)=(1.000000,1.000000,0.096000);
		rgb(227pt)=(1.000000,1.000000,0.092000);
		rgb(228pt)=(1.000000,1.000000,0.088000);
		rgb(229pt)=(1.000000,1.000000,0.084000);
		rgb(230pt)=(1.000000,1.000000,0.080000);
		rgb(231pt)=(1.000000,1.000000,0.076000);
		rgb(232pt)=(1.000000,1.000000,0.072000);
		rgb(233pt)=(1.000000,1.000000,0.068000);
		rgb(234pt)=(1.000000,1.000000,0.064000);
		rgb(235pt)=(1.000000,1.000000,0.060000);
		rgb(236pt)=(1.000000,1.000000,0.056000);
		rgb(237pt)=(1.000000,1.000000,0.052000);
		rgb(238pt)=(1.000000,1.000000,0.048000);
		rgb(239pt)=(1.000000,1.000000,0.044000);
		rgb(240pt)=(1.000000,1.000000,0.040000);
		rgb(241pt)=(1.000000,1.000000,0.036000);
		rgb(242pt)=(1.000000,1.000000,0.032000);
		rgb(243pt)=(1.000000,1.000000,0.028000);
		rgb(244pt)=(1.000000,1.000000,0.024000);
		rgb(245pt)=(1.000000,1.000000,0.020000);
		rgb(246pt)=(1.000000,1.000000,0.016000);
		rgb(247pt)=(1.000000,1.000000,0.012000);
		rgb(248pt)=(1.000000,1.000000,0.008000);
		rgb(249pt)=(1.000000,1.000000,0.004000);
		rgb(250pt)=(1.000000,1.000000,0.000000);
		rgb(251pt)=(1.000000,0.997333,0.000000);
		rgb(252pt)=(1.000000,0.994667,0.000000);
		rgb(253pt)=(1.000000,0.992000,0.000000);
		rgb(254pt)=(1.000000,0.989333,0.000000);
		rgb(255pt)=(1.000000,0.986667,0.000000);
		rgb(256pt)=(1.000000,0.984000,0.000000);
		rgb(257pt)=(1.000000,0.981333,0.000000);
		rgb(258pt)=(1.000000,0.978667,0.000000);
		rgb(259pt)=(1.000000,0.976000,0.000000);
		rgb(260pt)=(1.000000,0.973333,0.000000);
		rgb(261pt)=(1.000000,0.970667,0.000000);
		rgb(262pt)=(1.000000,0.968000,0.000000);
		rgb(263pt)=(1.000000,0.965333,0.000000);
		rgb(264pt)=(1.000000,0.962667,0.000000);
		rgb(265pt)=(1.000000,0.960000,0.000000);
		rgb(266pt)=(1.000000,0.957333,0.000000);
		rgb(267pt)=(1.000000,0.954667,0.000000);
		rgb(268pt)=(1.000000,0.952000,0.000000);
		rgb(269pt)=(1.000000,0.949333,0.000000);
		rgb(270pt)=(1.000000,0.946667,0.000000);
		rgb(271pt)=(1.000000,0.944000,0.000000);
		rgb(272pt)=(1.000000,0.941333,0.000000);
		rgb(273pt)=(1.000000,0.938667,0.000000);
		rgb(274pt)=(1.000000,0.936000,0.000000);
		rgb(275pt)=(1.000000,0.933333,0.000000);
		rgb(276pt)=(1.000000,0.930667,0.000000);
		rgb(277pt)=(1.000000,0.928000,0.000000);
		rgb(278pt)=(1.000000,0.925333,0.000000);
		rgb(279pt)=(1.000000,0.922667,0.000000);
		rgb(280pt)=(1.000000,0.920000,0.000000);
		rgb(281pt)=(1.000000,0.917333,0.000000);
		rgb(282pt)=(1.000000,0.914667,0.000000);
		rgb(283pt)=(1.000000,0.912000,0.000000);
		rgb(284pt)=(1.000000,0.909333,0.000000);
		rgb(285pt)=(1.000000,0.906667,0.000000);
		rgb(286pt)=(1.000000,0.904000,0.000000);
		rgb(287pt)=(1.000000,0.901333,0.000000);
		rgb(288pt)=(1.000000,0.898667,0.000000);
		rgb(289pt)=(1.000000,0.896000,0.000000);
		rgb(290pt)=(1.000000,0.893333,0.000000);
		rgb(291pt)=(1.000000,0.890667,0.000000);
		rgb(292pt)=(1.000000,0.888000,0.000000);
		rgb(293pt)=(1.000000,0.885333,0.000000);
		rgb(294pt)=(1.000000,0.882667,0.000000);
		rgb(295pt)=(1.000000,0.880000,0.000000);
		rgb(296pt)=(1.000000,0.877333,0.000000);
		rgb(297pt)=(1.000000,0.874667,0.000000);
		rgb(298pt)=(1.000000,0.872000,0.000000);
		rgb(299pt)=(1.000000,0.869333,0.000000);
		rgb(300pt)=(1.000000,0.866667,0.000000);
		rgb(301pt)=(1.000000,0.864000,0.000000);
		rgb(302pt)=(1.000000,0.861333,0.000000);
		rgb(303pt)=(1.000000,0.858667,0.000000);
		rgb(304pt)=(1.000000,0.856000,0.000000);
		rgb(305pt)=(1.000000,0.853333,0.000000);
		rgb(306pt)=(1.000000,0.850667,0.000000);
		rgb(307pt)=(1.000000,0.848000,0.000000);
		rgb(308pt)=(1.000000,0.845333,0.000000);
		rgb(309pt)=(1.000000,0.842667,0.000000);
		rgb(310pt)=(1.000000,0.840000,0.000000);
		rgb(311pt)=(1.000000,0.837333,0.000000);
		rgb(312pt)=(1.000000,0.834667,0.000000);
		rgb(313pt)=(1.000000,0.832000,0.000000);
		rgb(314pt)=(1.000000,0.829333,0.000000);
		rgb(315pt)=(1.000000,0.826667,0.000000);
		rgb(316pt)=(1.000000,0.824000,0.000000);
		rgb(317pt)=(1.000000,0.821333,0.000000);
		rgb(318pt)=(1.000000,0.818667,0.000000);
		rgb(319pt)=(1.000000,0.816000,0.000000);
		rgb(320pt)=(1.000000,0.813333,0.000000);
		rgb(321pt)=(1.000000,0.810667,0.000000);
		rgb(322pt)=(1.000000,0.808000,0.000000);
		rgb(323pt)=(1.000000,0.805333,0.000000);
		rgb(324pt)=(1.000000,0.802667,0.000000);
		rgb(325pt)=(1.000000,0.800000,0.000000);
		rgb(326pt)=(1.000000,0.797333,0.000000);
		rgb(327pt)=(1.000000,0.794667,0.000000);
		rgb(328pt)=(1.000000,0.792000,0.000000);
		rgb(329pt)=(1.000000,0.789333,0.000000);
		rgb(330pt)=(1.000000,0.786667,0.000000);
		rgb(331pt)=(1.000000,0.784000,0.000000);
		rgb(332pt)=(1.000000,0.781333,0.000000);
		rgb(333pt)=(1.000000,0.778667,0.000000);
		rgb(334pt)=(1.000000,0.776000,0.000000);
		rgb(335pt)=(1.000000,0.773333,0.000000);
		rgb(336pt)=(1.000000,0.770667,0.000000);
		rgb(337pt)=(1.000000,0.768000,0.000000);
		rgb(338pt)=(1.000000,0.765333,0.000000);
		rgb(339pt)=(1.000000,0.762667,0.000000);
		rgb(340pt)=(1.000000,0.760000,0.000000);
		rgb(341pt)=(1.000000,0.757333,0.000000);
		rgb(342pt)=(1.000000,0.754667,0.000000);
		rgb(343pt)=(1.000000,0.752000,0.000000);
		rgb(344pt)=(1.000000,0.749333,0.000000);
		rgb(345pt)=(1.000000,0.746667,0.000000);
		rgb(346pt)=(1.000000,0.744000,0.000000);
		rgb(347pt)=(1.000000,0.741333,0.000000);
		rgb(348pt)=(1.000000,0.738667,0.000000);
		rgb(349pt)=(1.000000,0.736000,0.000000);
		rgb(350pt)=(1.000000,0.733333,0.000000);
		rgb(351pt)=(1.000000,0.730667,0.000000);
		rgb(352pt)=(1.000000,0.728000,0.000000);
		rgb(353pt)=(1.000000,0.725333,0.000000);
		rgb(354pt)=(1.000000,0.722667,0.000000);
		rgb(355pt)=(1.000000,0.720000,0.000000);
		rgb(356pt)=(1.000000,0.717333,0.000000);
		rgb(357pt)=(1.000000,0.714667,0.000000);
		rgb(358pt)=(1.000000,0.712000,0.000000);
		rgb(359pt)=(1.000000,0.709333,0.000000);
		rgb(360pt)=(1.000000,0.706667,0.000000);
		rgb(361pt)=(1.000000,0.704000,0.000000);
		rgb(362pt)=(1.000000,0.701333,0.000000);
		rgb(363pt)=(1.000000,0.698667,0.000000);
		rgb(364pt)=(1.000000,0.696000,0.000000);
		rgb(365pt)=(1.000000,0.693333,0.000000);
		rgb(366pt)=(1.000000,0.690667,0.000000);
		rgb(367pt)=(1.000000,0.688000,0.000000);
		rgb(368pt)=(1.000000,0.685333,0.000000);
		rgb(369pt)=(1.000000,0.682667,0.000000);
		rgb(370pt)=(1.000000,0.680000,0.000000);
		rgb(371pt)=(1.000000,0.677333,0.000000);
		rgb(372pt)=(1.000000,0.674667,0.000000);
		rgb(373pt)=(1.000000,0.672000,0.000000);
		rgb(374pt)=(1.000000,0.669333,0.000000);
		rgb(375pt)=(1.000000,0.666667,0.000000);
		rgb(376pt)=(1.000000,0.664000,0.000000);
		rgb(377pt)=(1.000000,0.661333,0.000000);
		rgb(378pt)=(1.000000,0.658667,0.000000);
		rgb(379pt)=(1.000000,0.656000,0.000000);
		rgb(380pt)=(1.000000,0.653333,0.000000);
		rgb(381pt)=(1.000000,0.650667,0.000000);
		rgb(382pt)=(1.000000,0.648000,0.000000);
		rgb(383pt)=(1.000000,0.645333,0.000000);
		rgb(384pt)=(1.000000,0.642667,0.000000);
		rgb(385pt)=(1.000000,0.640000,0.000000);
		rgb(386pt)=(1.000000,0.637333,0.000000);
		rgb(387pt)=(1.000000,0.634667,0.000000);
		rgb(388pt)=(1.000000,0.632000,0.000000);
		rgb(389pt)=(1.000000,0.629333,0.000000);
		rgb(390pt)=(1.000000,0.626667,0.000000);
		rgb(391pt)=(1.000000,0.624000,0.000000);
		rgb(392pt)=(1.000000,0.621333,0.000000);
		rgb(393pt)=(1.000000,0.618667,0.000000);
		rgb(394pt)=(1.000000,0.616000,0.000000);
		rgb(395pt)=(1.000000,0.613333,0.000000);
		rgb(396pt)=(1.000000,0.610667,0.000000);
		rgb(397pt)=(1.000000,0.608000,0.000000);
		rgb(398pt)=(1.000000,0.605333,0.000000);
		rgb(399pt)=(1.000000,0.602667,0.000000);
		rgb(400pt)=(1.000000,0.600000,0.000000);
		rgb(401pt)=(1.000000,0.597333,0.000000);
		rgb(402pt)=(1.000000,0.594667,0.000000);
		rgb(403pt)=(1.000000,0.592000,0.000000);
		rgb(404pt)=(1.000000,0.589333,0.000000);
		rgb(405pt)=(1.000000,0.586667,0.000000);
		rgb(406pt)=(1.000000,0.584000,0.000000);
		rgb(407pt)=(1.000000,0.581333,0.000000);
		rgb(408pt)=(1.000000,0.578667,0.000000);
		rgb(409pt)=(1.000000,0.576000,0.000000);
		rgb(410pt)=(1.000000,0.573333,0.000000);
		rgb(411pt)=(1.000000,0.570667,0.000000);
		rgb(412pt)=(1.000000,0.568000,0.000000);
		rgb(413pt)=(1.000000,0.565333,0.000000);
		rgb(414pt)=(1.000000,0.562667,0.000000);
		rgb(415pt)=(1.000000,0.560000,0.000000);
		rgb(416pt)=(1.000000,0.557333,0.000000);
		rgb(417pt)=(1.000000,0.554667,0.000000);
		rgb(418pt)=(1.000000,0.552000,0.000000);
		rgb(419pt)=(1.000000,0.549333,0.000000);
		rgb(420pt)=(1.000000,0.546667,0.000000);
		rgb(421pt)=(1.000000,0.544000,0.000000);
		rgb(422pt)=(1.000000,0.541333,0.000000);
		rgb(423pt)=(1.000000,0.538667,0.000000);
		rgb(424pt)=(1.000000,0.536000,0.000000);
		rgb(425pt)=(1.000000,0.533333,0.000000);
		rgb(426pt)=(1.000000,0.530667,0.000000);
		rgb(427pt)=(1.000000,0.528000,0.000000);
		rgb(428pt)=(1.000000,0.525333,0.000000);
		rgb(429pt)=(1.000000,0.522667,0.000000);
		rgb(430pt)=(1.000000,0.520000,0.000000);
		rgb(431pt)=(1.000000,0.517333,0.000000);
		rgb(432pt)=(1.000000,0.514667,0.000000);
		rgb(433pt)=(1.000000,0.512000,0.000000);
		rgb(434pt)=(1.000000,0.509333,0.000000);
		rgb(435pt)=(1.000000,0.506667,0.000000);
		rgb(436pt)=(1.000000,0.504000,0.000000);
		rgb(437pt)=(1.000000,0.501333,0.000000);
		rgb(438pt)=(1.000000,0.498667,0.000000);
		rgb(439pt)=(1.000000,0.496000,0.000000);
		rgb(440pt)=(1.000000,0.493333,0.000000);
		rgb(441pt)=(1.000000,0.490667,0.000000);
		rgb(442pt)=(1.000000,0.488000,0.000000);
		rgb(443pt)=(1.000000,0.485333,0.000000);
		rgb(444pt)=(1.000000,0.482667,0.000000);
		rgb(445pt)=(1.000000,0.480000,0.000000);
		rgb(446pt)=(1.000000,0.477333,0.000000);
		rgb(447pt)=(1.000000,0.474667,0.000000);
		rgb(448pt)=(1.000000,0.472000,0.000000);
		rgb(449pt)=(1.000000,0.469333,0.000000);
		rgb(450pt)=(1.000000,0.466667,0.000000);
		rgb(451pt)=(1.000000,0.464000,0.000000);
		rgb(452pt)=(1.000000,0.461333,0.000000);
		rgb(453pt)=(1.000000,0.458667,0.000000);
		rgb(454pt)=(1.000000,0.456000,0.000000);
		rgb(455pt)=(1.000000,0.453333,0.000000);
		rgb(456pt)=(1.000000,0.450667,0.000000);
		rgb(457pt)=(1.000000,0.448000,0.000000);
		rgb(458pt)=(1.000000,0.445333,0.000000);
		rgb(459pt)=(1.000000,0.442667,0.000000);
		rgb(460pt)=(1.000000,0.440000,0.000000);
		rgb(461pt)=(1.000000,0.437333,0.000000);
		rgb(462pt)=(1.000000,0.434667,0.000000);
		rgb(463pt)=(1.000000,0.432000,0.000000);
		rgb(464pt)=(1.000000,0.429333,0.000000);
		rgb(465pt)=(1.000000,0.426667,0.000000);
		rgb(466pt)=(1.000000,0.424000,0.000000);
		rgb(467pt)=(1.000000,0.421333,0.000000);
		rgb(468pt)=(1.000000,0.418667,0.000000);
		rgb(469pt)=(1.000000,0.416000,0.000000);
		rgb(470pt)=(1.000000,0.413333,0.000000);
		rgb(471pt)=(1.000000,0.410667,0.000000);
		rgb(472pt)=(1.000000,0.408000,0.000000);
		rgb(473pt)=(1.000000,0.405333,0.000000);
		rgb(474pt)=(1.000000,0.402667,0.000000);
		rgb(475pt)=(1.000000,0.400000,0.000000);
		rgb(476pt)=(1.000000,0.397333,0.000000);
		rgb(477pt)=(1.000000,0.394667,0.000000);
		rgb(478pt)=(1.000000,0.392000,0.000000);
		rgb(479pt)=(1.000000,0.389333,0.000000);
		rgb(480pt)=(1.000000,0.386667,0.000000);
		rgb(481pt)=(1.000000,0.384000,0.000000);
		rgb(482pt)=(1.000000,0.381333,0.000000);
		rgb(483pt)=(1.000000,0.378667,0.000000);
		rgb(484pt)=(1.000000,0.376000,0.000000);
		rgb(485pt)=(1.000000,0.373333,0.000000);
		rgb(486pt)=(1.000000,0.370667,0.000000);
		rgb(487pt)=(1.000000,0.368000,0.000000);
		rgb(488pt)=(1.000000,0.365333,0.000000);
		rgb(489pt)=(1.000000,0.362667,0.000000);
		rgb(490pt)=(1.000000,0.360000,0.000000);
		rgb(491pt)=(1.000000,0.357333,0.000000);
		rgb(492pt)=(1.000000,0.354667,0.000000);
		rgb(493pt)=(1.000000,0.352000,0.000000);
		rgb(494pt)=(1.000000,0.349333,0.000000);
		rgb(495pt)=(1.000000,0.346667,0.000000);
		rgb(496pt)=(1.000000,0.344000,0.000000);
		rgb(497pt)=(1.000000,0.341333,0.000000);
		rgb(498pt)=(1.000000,0.338667,0.000000);
		rgb(499pt)=(1.000000,0.336000,0.000000);
		rgb(500pt)=(1.000000,0.333333,0.000000);
		rgb(501pt)=(1.000000,0.330667,0.000000);
		rgb(502pt)=(1.000000,0.328000,0.000000);
		rgb(503pt)=(1.000000,0.325333,0.000000);
		rgb(504pt)=(1.000000,0.322667,0.000000);
		rgb(505pt)=(1.000000,0.320000,0.000000);
		rgb(506pt)=(1.000000,0.317333,0.000000);
		rgb(507pt)=(1.000000,0.314667,0.000000);
		rgb(508pt)=(1.000000,0.312000,0.000000);
		rgb(509pt)=(1.000000,0.309333,0.000000);
		rgb(510pt)=(1.000000,0.306667,0.000000);
		rgb(511pt)=(1.000000,0.304000,0.000000);
		rgb(512pt)=(1.000000,0.301333,0.000000);
		rgb(513pt)=(1.000000,0.298667,0.000000);
		rgb(514pt)=(1.000000,0.296000,0.000000);
		rgb(515pt)=(1.000000,0.293333,0.000000);
		rgb(516pt)=(1.000000,0.290667,0.000000);
		rgb(517pt)=(1.000000,0.288000,0.000000);
		rgb(518pt)=(1.000000,0.285333,0.000000);
		rgb(519pt)=(1.000000,0.282667,0.000000);
		rgb(520pt)=(1.000000,0.280000,0.000000);
		rgb(521pt)=(1.000000,0.277333,0.000000);
		rgb(522pt)=(1.000000,0.274667,0.000000);
		rgb(523pt)=(1.000000,0.272000,0.000000);
		rgb(524pt)=(1.000000,0.269333,0.000000);
		rgb(525pt)=(1.000000,0.266667,0.000000);
		rgb(526pt)=(1.000000,0.264000,0.000000);
		rgb(527pt)=(1.000000,0.261333,0.000000);
		rgb(528pt)=(1.000000,0.258667,0.000000);
		rgb(529pt)=(1.000000,0.256000,0.000000);
		rgb(530pt)=(1.000000,0.253333,0.000000);
		rgb(531pt)=(1.000000,0.250667,0.000000);
		rgb(532pt)=(1.000000,0.248000,0.000000);
		rgb(533pt)=(1.000000,0.245333,0.000000);
		rgb(534pt)=(1.000000,0.242667,0.000000);
		rgb(535pt)=(1.000000,0.240000,0.000000);
		rgb(536pt)=(1.000000,0.237333,0.000000);
		rgb(537pt)=(1.000000,0.234667,0.000000);
		rgb(538pt)=(1.000000,0.232000,0.000000);
		rgb(539pt)=(1.000000,0.229333,0.000000);
		rgb(540pt)=(1.000000,0.226667,0.000000);
		rgb(541pt)=(1.000000,0.224000,0.000000);
		rgb(542pt)=(1.000000,0.221333,0.000000);
		rgb(543pt)=(1.000000,0.218667,0.000000);
		rgb(544pt)=(1.000000,0.216000,0.000000);
		rgb(545pt)=(1.000000,0.213333,0.000000);
		rgb(546pt)=(1.000000,0.210667,0.000000);
		rgb(547pt)=(1.000000,0.208000,0.000000);
		rgb(548pt)=(1.000000,0.205333,0.000000);
		rgb(549pt)=(1.000000,0.202667,0.000000);
		rgb(550pt)=(1.000000,0.200000,0.000000);
		rgb(551pt)=(1.000000,0.197333,0.000000);
		rgb(552pt)=(1.000000,0.194667,0.000000);
		rgb(553pt)=(1.000000,0.192000,0.000000);
		rgb(554pt)=(1.000000,0.189333,0.000000);
		rgb(555pt)=(1.000000,0.186667,0.000000);
		rgb(556pt)=(1.000000,0.184000,0.000000);
		rgb(557pt)=(1.000000,0.181333,0.000000);
		rgb(558pt)=(1.000000,0.178667,0.000000);
		rgb(559pt)=(1.000000,0.176000,0.000000);
		rgb(560pt)=(1.000000,0.173333,0.000000);
		rgb(561pt)=(1.000000,0.170667,0.000000);
		rgb(562pt)=(1.000000,0.168000,0.000000);
		rgb(563pt)=(1.000000,0.165333,0.000000);
		rgb(564pt)=(1.000000,0.162667,0.000000);
		rgb(565pt)=(1.000000,0.160000,0.000000);
		rgb(566pt)=(1.000000,0.157333,0.000000);
		rgb(567pt)=(1.000000,0.154667,0.000000);
		rgb(568pt)=(1.000000,0.152000,0.000000);
		rgb(569pt)=(1.000000,0.149333,0.000000);
		rgb(570pt)=(1.000000,0.146667,0.000000);
		rgb(571pt)=(1.000000,0.144000,0.000000);
		rgb(572pt)=(1.000000,0.141333,0.000000);
		rgb(573pt)=(1.000000,0.138667,0.000000);
		rgb(574pt)=(1.000000,0.136000,0.000000);
		rgb(575pt)=(1.000000,0.133333,0.000000);
		rgb(576pt)=(1.000000,0.130667,0.000000);
		rgb(577pt)=(1.000000,0.128000,0.000000);
		rgb(578pt)=(1.000000,0.125333,0.000000);
		rgb(579pt)=(1.000000,0.122667,0.000000);
		rgb(580pt)=(1.000000,0.120000,0.000000);
		rgb(581pt)=(1.000000,0.117333,0.000000);
		rgb(582pt)=(1.000000,0.114667,0.000000);
		rgb(583pt)=(1.000000,0.112000,0.000000);
		rgb(584pt)=(1.000000,0.109333,0.000000);
		rgb(585pt)=(1.000000,0.106667,0.000000);
		rgb(586pt)=(1.000000,0.104000,0.000000);
		rgb(587pt)=(1.000000,0.101333,0.000000);
		rgb(588pt)=(1.000000,0.098667,0.000000);
		rgb(589pt)=(1.000000,0.096000,0.000000);
		rgb(590pt)=(1.000000,0.093333,0.000000);
		rgb(591pt)=(1.000000,0.090667,0.000000);
		rgb(592pt)=(1.000000,0.088000,0.000000);
		rgb(593pt)=(1.000000,0.085333,0.000000);
		rgb(594pt)=(1.000000,0.082667,0.000000);
		rgb(595pt)=(1.000000,0.080000,0.000000);
		rgb(596pt)=(1.000000,0.077333,0.000000);
		rgb(597pt)=(1.000000,0.074667,0.000000);
		rgb(598pt)=(1.000000,0.072000,0.000000);
		rgb(599pt)=(1.000000,0.069333,0.000000);
		rgb(600pt)=(1.000000,0.066667,0.000000);
		rgb(601pt)=(1.000000,0.064000,0.000000);
		rgb(602pt)=(1.000000,0.061333,0.000000);
		rgb(603pt)=(1.000000,0.058667,0.000000);
		rgb(604pt)=(1.000000,0.056000,0.000000);
		rgb(605pt)=(1.000000,0.053333,0.000000);
		rgb(606pt)=(1.000000,0.050667,0.000000);
		rgb(607pt)=(1.000000,0.048000,0.000000);
		rgb(608pt)=(1.000000,0.045333,0.000000);
		rgb(609pt)=(1.000000,0.042667,0.000000);
		rgb(610pt)=(1.000000,0.040000,0.000000);
		rgb(611pt)=(1.000000,0.037333,0.000000);
		rgb(612pt)=(1.000000,0.034667,0.000000);
		rgb(613pt)=(1.000000,0.032000,0.000000);
		rgb(614pt)=(1.000000,0.029333,0.000000);
		rgb(615pt)=(1.000000,0.026667,0.000000);
		rgb(616pt)=(1.000000,0.024000,0.000000);
		rgb(617pt)=(1.000000,0.021333,0.000000);
		rgb(618pt)=(1.000000,0.018667,0.000000);
		rgb(619pt)=(1.000000,0.016000,0.000000);
		rgb(620pt)=(1.000000,0.013333,0.000000);
		rgb(621pt)=(1.000000,0.010667,0.000000);
		rgb(622pt)=(1.000000,0.008000,0.000000);
		rgb(623pt)=(1.000000,0.005333,0.000000);
		rgb(624pt)=(1.000000,0.002667,0.000000);
		rgb(625pt)=(1.000000,0.000000,0.000000);
		rgb(626pt)=(0.997333,0.000000,0.000000);
		rgb(627pt)=(0.994667,0.000000,0.000000);
		rgb(628pt)=(0.992000,0.000000,0.000000);
		rgb(629pt)=(0.989333,0.000000,0.000000);
		rgb(630pt)=(0.986667,0.000000,0.000000);
		rgb(631pt)=(0.984000,0.000000,0.000000);
		rgb(632pt)=(0.981333,0.000000,0.000000);
		rgb(633pt)=(0.978667,0.000000,0.000000);
		rgb(634pt)=(0.976000,0.000000,0.000000);
		rgb(635pt)=(0.973333,0.000000,0.000000);
		rgb(636pt)=(0.970667,0.000000,0.000000);
		rgb(637pt)=(0.968000,0.000000,0.000000);
		rgb(638pt)=(0.965333,0.000000,0.000000);
		rgb(639pt)=(0.962667,0.000000,0.000000);
		rgb(640pt)=(0.960000,0.000000,0.000000);
		rgb(641pt)=(0.957333,0.000000,0.000000);
		rgb(642pt)=(0.954667,0.000000,0.000000);
		rgb(643pt)=(0.952000,0.000000,0.000000);
		rgb(644pt)=(0.949333,0.000000,0.000000);
		rgb(645pt)=(0.946667,0.000000,0.000000);
		rgb(646pt)=(0.944000,0.000000,0.000000);
		rgb(647pt)=(0.941333,0.000000,0.000000);
		rgb(648pt)=(0.938667,0.000000,0.000000);
		rgb(649pt)=(0.936000,0.000000,0.000000);
		rgb(650pt)=(0.933333,0.000000,0.000000);
		rgb(651pt)=(0.930667,0.000000,0.000000);
		rgb(652pt)=(0.928000,0.000000,0.000000);
		rgb(653pt)=(0.925333,0.000000,0.000000);
		rgb(654pt)=(0.922667,0.000000,0.000000);
		rgb(655pt)=(0.920000,0.000000,0.000000);
		rgb(656pt)=(0.917333,0.000000,0.000000);
		rgb(657pt)=(0.914667,0.000000,0.000000);
		rgb(658pt)=(0.912000,0.000000,0.000000);
		rgb(659pt)=(0.909333,0.000000,0.000000);
		rgb(660pt)=(0.906667,0.000000,0.000000);
		rgb(661pt)=(0.904000,0.000000,0.000000);
		rgb(662pt)=(0.901333,0.000000,0.000000);
		rgb(663pt)=(0.898667,0.000000,0.000000);
		rgb(664pt)=(0.896000,0.000000,0.000000);
		rgb(665pt)=(0.893333,0.000000,0.000000);
		rgb(666pt)=(0.890667,0.000000,0.000000);
		rgb(667pt)=(0.888000,0.000000,0.000000);
		rgb(668pt)=(0.885333,0.000000,0.000000);
		rgb(669pt)=(0.882667,0.000000,0.000000);
		rgb(670pt)=(0.880000,0.000000,0.000000);
		rgb(671pt)=(0.877333,0.000000,0.000000);
		rgb(672pt)=(0.874667,0.000000,0.000000);
		rgb(673pt)=(0.872000,0.000000,0.000000);
		rgb(674pt)=(0.869333,0.000000,0.000000);
		rgb(675pt)=(0.866667,0.000000,0.000000);
		rgb(676pt)=(0.864000,0.000000,0.000000);
		rgb(677pt)=(0.861333,0.000000,0.000000);
		rgb(678pt)=(0.858667,0.000000,0.000000);
		rgb(679pt)=(0.856000,0.000000,0.000000);
		rgb(680pt)=(0.853333,0.000000,0.000000);
		rgb(681pt)=(0.850667,0.000000,0.000000);
		rgb(682pt)=(0.848000,0.000000,0.000000);
		rgb(683pt)=(0.845333,0.000000,0.000000);
		rgb(684pt)=(0.842667,0.000000,0.000000);
		rgb(685pt)=(0.840000,0.000000,0.000000);
		rgb(686pt)=(0.837333,0.000000,0.000000);
		rgb(687pt)=(0.834667,0.000000,0.000000);
		rgb(688pt)=(0.832000,0.000000,0.000000);
		rgb(689pt)=(0.829333,0.000000,0.000000);
		rgb(690pt)=(0.826667,0.000000,0.000000);
		rgb(691pt)=(0.824000,0.000000,0.000000);
		rgb(692pt)=(0.821333,0.000000,0.000000);
		rgb(693pt)=(0.818667,0.000000,0.000000);
		rgb(694pt)=(0.816000,0.000000,0.000000);
		rgb(695pt)=(0.813333,0.000000,0.000000);
		rgb(696pt)=(0.810667,0.000000,0.000000);
		rgb(697pt)=(0.808000,0.000000,0.000000);
		rgb(698pt)=(0.805333,0.000000,0.000000);
		rgb(699pt)=(0.802667,0.000000,0.000000);
		rgb(700pt)=(0.800000,0.000000,0.000000);
		rgb(701pt)=(0.797333,0.000000,0.000000);
		rgb(702pt)=(0.794667,0.000000,0.000000);
		rgb(703pt)=(0.792000,0.000000,0.000000);
		rgb(704pt)=(0.789333,0.000000,0.000000);
		rgb(705pt)=(0.786667,0.000000,0.000000);
		rgb(706pt)=(0.784000,0.000000,0.000000);
		rgb(707pt)=(0.781333,0.000000,0.000000);
		rgb(708pt)=(0.778667,0.000000,0.000000);
		rgb(709pt)=(0.776000,0.000000,0.000000);
		rgb(710pt)=(0.773333,0.000000,0.000000);
		rgb(711pt)=(0.770667,0.000000,0.000000);
		rgb(712pt)=(0.768000,0.000000,0.000000);
		rgb(713pt)=(0.765333,0.000000,0.000000);
		rgb(714pt)=(0.762667,0.000000,0.000000);
		rgb(715pt)=(0.760000,0.000000,0.000000);
		rgb(716pt)=(0.757333,0.000000,0.000000);
		rgb(717pt)=(0.754667,0.000000,0.000000);
		rgb(718pt)=(0.752000,0.000000,0.000000);
		rgb(719pt)=(0.749333,0.000000,0.000000);
		rgb(720pt)=(0.746667,0.000000,0.000000);
		rgb(721pt)=(0.744000,0.000000,0.000000);
		rgb(722pt)=(0.741333,0.000000,0.000000);
		rgb(723pt)=(0.738667,0.000000,0.000000);
		rgb(724pt)=(0.736000,0.000000,0.000000);
		rgb(725pt)=(0.733333,0.000000,0.000000);
		rgb(726pt)=(0.730667,0.000000,0.000000);
		rgb(727pt)=(0.728000,0.000000,0.000000);
		rgb(728pt)=(0.725333,0.000000,0.000000);
		rgb(729pt)=(0.722667,0.000000,0.000000);
		rgb(730pt)=(0.720000,0.000000,0.000000);
		rgb(731pt)=(0.717333,0.000000,0.000000);
		rgb(732pt)=(0.714667,0.000000,0.000000);
		rgb(733pt)=(0.712000,0.000000,0.000000);
		rgb(734pt)=(0.709333,0.000000,0.000000);
		rgb(735pt)=(0.706667,0.000000,0.000000);
		rgb(736pt)=(0.704000,0.000000,0.000000);
		rgb(737pt)=(0.701333,0.000000,0.000000);
		rgb(738pt)=(0.698667,0.000000,0.000000);
		rgb(739pt)=(0.696000,0.000000,0.000000);
		rgb(740pt)=(0.693333,0.000000,0.000000);
		rgb(741pt)=(0.690667,0.000000,0.000000);
		rgb(742pt)=(0.688000,0.000000,0.000000);
		rgb(743pt)=(0.685333,0.000000,0.000000);
		rgb(744pt)=(0.682667,0.000000,0.000000);
		rgb(745pt)=(0.680000,0.000000,0.000000);
		rgb(746pt)=(0.677333,0.000000,0.000000);
		rgb(747pt)=(0.674667,0.000000,0.000000);
		rgb(748pt)=(0.672000,0.000000,0.000000);
		rgb(749pt)=(0.669333,0.000000,0.000000);
		rgb(750pt)=(0.666667,0.000000,0.000000);
		rgb(751pt)=(0.664000,0.000000,0.000000);
		rgb(752pt)=(0.661333,0.000000,0.000000);
		rgb(753pt)=(0.658667,0.000000,0.000000);
		rgb(754pt)=(0.656000,0.000000,0.000000);
		rgb(755pt)=(0.653333,0.000000,0.000000);
		rgb(756pt)=(0.650667,0.000000,0.000000);
		rgb(757pt)=(0.648000,0.000000,0.000000);
		rgb(758pt)=(0.645333,0.000000,0.000000);
		rgb(759pt)=(0.642667,0.000000,0.000000);
		rgb(760pt)=(0.640000,0.000000,0.000000);
		rgb(761pt)=(0.637333,0.000000,0.000000);
		rgb(762pt)=(0.634667,0.000000,0.000000);
		rgb(763pt)=(0.632000,0.000000,0.000000);
		rgb(764pt)=(0.629333,0.000000,0.000000);
		rgb(765pt)=(0.626667,0.000000,0.000000);
		rgb(766pt)=(0.624000,0.000000,0.000000);
		rgb(767pt)=(0.621333,0.000000,0.000000);
		rgb(768pt)=(0.618667,0.000000,0.000000);
		rgb(769pt)=(0.616000,0.000000,0.000000);
		rgb(770pt)=(0.613333,0.000000,0.000000);
		rgb(771pt)=(0.610667,0.000000,0.000000);
		rgb(772pt)=(0.608000,0.000000,0.000000);
		rgb(773pt)=(0.605333,0.000000,0.000000);
		rgb(774pt)=(0.602667,0.000000,0.000000);
		rgb(775pt)=(0.600000,0.000000,0.000000);
		rgb(776pt)=(0.597333,0.000000,0.000000);
		rgb(777pt)=(0.594667,0.000000,0.000000);
		rgb(778pt)=(0.592000,0.000000,0.000000);
		rgb(779pt)=(0.589333,0.000000,0.000000);
		rgb(780pt)=(0.586667,0.000000,0.000000);
		rgb(781pt)=(0.584000,0.000000,0.000000);
		rgb(782pt)=(0.581333,0.000000,0.000000);
		rgb(783pt)=(0.578667,0.000000,0.000000);
		rgb(784pt)=(0.576000,0.000000,0.000000);
		rgb(785pt)=(0.573333,0.000000,0.000000);
		rgb(786pt)=(0.570667,0.000000,0.000000);
		rgb(787pt)=(0.568000,0.000000,0.000000);
		rgb(788pt)=(0.565333,0.000000,0.000000);
		rgb(789pt)=(0.562667,0.000000,0.000000);
		rgb(790pt)=(0.560000,0.000000,0.000000);
		rgb(791pt)=(0.557333,0.000000,0.000000);
		rgb(792pt)=(0.554667,0.000000,0.000000);
		rgb(793pt)=(0.552000,0.000000,0.000000);
		rgb(794pt)=(0.549333,0.000000,0.000000);
		rgb(795pt)=(0.546667,0.000000,0.000000);
		rgb(796pt)=(0.544000,0.000000,0.000000);
		rgb(797pt)=(0.541333,0.000000,0.000000);
		rgb(798pt)=(0.538667,0.000000,0.000000);
		rgb(799pt)=(0.536000,0.000000,0.000000);
		rgb(800pt)=(0.533333,0.000000,0.000000);
		rgb(801pt)=(0.530667,0.000000,0.000000);
		rgb(802pt)=(0.528000,0.000000,0.000000);
		rgb(803pt)=(0.525333,0.000000,0.000000);
		rgb(804pt)=(0.522667,0.000000,0.000000);
		rgb(805pt)=(0.520000,0.000000,0.000000);
		rgb(806pt)=(0.517333,0.000000,0.000000);
		rgb(807pt)=(0.514667,0.000000,0.000000);
		rgb(808pt)=(0.512000,0.000000,0.000000);
		rgb(809pt)=(0.509333,0.000000,0.000000);
		rgb(810pt)=(0.506667,0.000000,0.000000);
		rgb(811pt)=(0.504000,0.000000,0.000000);
		rgb(812pt)=(0.501333,0.000000,0.000000);
		rgb(813pt)=(0.498667,0.000000,0.000000);
		rgb(814pt)=(0.496000,0.000000,0.000000);
		rgb(815pt)=(0.493333,0.000000,0.000000);
		rgb(816pt)=(0.490667,0.000000,0.000000);
		rgb(817pt)=(0.488000,0.000000,0.000000);
		rgb(818pt)=(0.485333,0.000000,0.000000);
		rgb(819pt)=(0.482667,0.000000,0.000000);
		rgb(820pt)=(0.480000,0.000000,0.000000);
		rgb(821pt)=(0.477333,0.000000,0.000000);
		rgb(822pt)=(0.474667,0.000000,0.000000);
		rgb(823pt)=(0.472000,0.000000,0.000000);
		rgb(824pt)=(0.469333,0.000000,0.000000);
		rgb(825pt)=(0.466667,0.000000,0.000000);
		rgb(826pt)=(0.464000,0.000000,0.000000);
		rgb(827pt)=(0.461333,0.000000,0.000000);
		rgb(828pt)=(0.458667,0.000000,0.000000);
		rgb(829pt)=(0.456000,0.000000,0.000000);
		rgb(830pt)=(0.453333,0.000000,0.000000);
		rgb(831pt)=(0.450667,0.000000,0.000000);
		rgb(832pt)=(0.448000,0.000000,0.000000);
		rgb(833pt)=(0.445333,0.000000,0.000000);
		rgb(834pt)=(0.442667,0.000000,0.000000);
		rgb(835pt)=(0.440000,0.000000,0.000000);
		rgb(836pt)=(0.437333,0.000000,0.000000);
		rgb(837pt)=(0.434667,0.000000,0.000000);
		rgb(838pt)=(0.432000,0.000000,0.000000);
		rgb(839pt)=(0.429333,0.000000,0.000000);
		rgb(840pt)=(0.426667,0.000000,0.000000);
		rgb(841pt)=(0.424000,0.000000,0.000000);
		rgb(842pt)=(0.421333,0.000000,0.000000);
		rgb(843pt)=(0.418667,0.000000,0.000000);
		rgb(844pt)=(0.416000,0.000000,0.000000);
		rgb(845pt)=(0.413333,0.000000,0.000000);
		rgb(846pt)=(0.410667,0.000000,0.000000);
		rgb(847pt)=(0.408000,0.000000,0.000000);
		rgb(848pt)=(0.405333,0.000000,0.000000);
		rgb(849pt)=(0.402667,0.000000,0.000000);
		rgb(850pt)=(0.400000,0.000000,0.000000);
		rgb(851pt)=(0.397333,0.000000,0.000000);
		rgb(852pt)=(0.394667,0.000000,0.000000);
		rgb(853pt)=(0.392000,0.000000,0.000000);
		rgb(854pt)=(0.389333,0.000000,0.000000);
		rgb(855pt)=(0.386667,0.000000,0.000000);
		rgb(856pt)=(0.384000,0.000000,0.000000);
		rgb(857pt)=(0.381333,0.000000,0.000000);
		rgb(858pt)=(0.378667,0.000000,0.000000);
		rgb(859pt)=(0.376000,0.000000,0.000000);
		rgb(860pt)=(0.373333,0.000000,0.000000);
		rgb(861pt)=(0.370667,0.000000,0.000000);
		rgb(862pt)=(0.368000,0.000000,0.000000);
		rgb(863pt)=(0.365333,0.000000,0.000000);
		rgb(864pt)=(0.362667,0.000000,0.000000);
		rgb(865pt)=(0.360000,0.000000,0.000000);
		rgb(866pt)=(0.357333,0.000000,0.000000);
		rgb(867pt)=(0.354667,0.000000,0.000000);
		rgb(868pt)=(0.352000,0.000000,0.000000);
		rgb(869pt)=(0.349333,0.000000,0.000000);
		rgb(870pt)=(0.346667,0.000000,0.000000);
		rgb(871pt)=(0.344000,0.000000,0.000000);
		rgb(872pt)=(0.341333,0.000000,0.000000);
		rgb(873pt)=(0.338667,0.000000,0.000000);
		rgb(874pt)=(0.336000,0.000000,0.000000);
		rgb(875pt)=(0.333333,0.000000,0.000000);
		rgb(876pt)=(0.330667,0.000000,0.000000);
		rgb(877pt)=(0.328000,0.000000,0.000000);
		rgb(878pt)=(0.325333,0.000000,0.000000);
		rgb(879pt)=(0.322667,0.000000,0.000000);
		rgb(880pt)=(0.320000,0.000000,0.000000);
		rgb(881pt)=(0.317333,0.000000,0.000000);
		rgb(882pt)=(0.314667,0.000000,0.000000);
		rgb(883pt)=(0.312000,0.000000,0.000000);
		rgb(884pt)=(0.309333,0.000000,0.000000);
		rgb(885pt)=(0.306667,0.000000,0.000000);
		rgb(886pt)=(0.304000,0.000000,0.000000);
		rgb(887pt)=(0.301333,0.000000,0.000000);
		rgb(888pt)=(0.298667,0.000000,0.000000);
		rgb(889pt)=(0.296000,0.000000,0.000000);
		rgb(890pt)=(0.293333,0.000000,0.000000);
		rgb(891pt)=(0.290667,0.000000,0.000000);
		rgb(892pt)=(0.288000,0.000000,0.000000);
		rgb(893pt)=(0.285333,0.000000,0.000000);
		rgb(894pt)=(0.282667,0.000000,0.000000);
		rgb(895pt)=(0.280000,0.000000,0.000000);
		rgb(896pt)=(0.277333,0.000000,0.000000);
		rgb(897pt)=(0.274667,0.000000,0.000000);
		rgb(898pt)=(0.272000,0.000000,0.000000);
		rgb(899pt)=(0.269333,0.000000,0.000000);
		rgb(900pt)=(0.266667,0.000000,0.000000);
		rgb(901pt)=(0.264000,0.000000,0.000000);
		rgb(902pt)=(0.261333,0.000000,0.000000);
		rgb(903pt)=(0.258667,0.000000,0.000000);
		rgb(904pt)=(0.256000,0.000000,0.000000);
		rgb(905pt)=(0.253333,0.000000,0.000000);
		rgb(906pt)=(0.250667,0.000000,0.000000);
		rgb(907pt)=(0.248000,0.000000,0.000000);
		rgb(908pt)=(0.245333,0.000000,0.000000);
		rgb(909pt)=(0.242667,0.000000,0.000000);
		rgb(910pt)=(0.240000,0.000000,0.000000);
		rgb(911pt)=(0.237333,0.000000,0.000000);
		rgb(912pt)=(0.234667,0.000000,0.000000);
		rgb(913pt)=(0.232000,0.000000,0.000000);
		rgb(914pt)=(0.229333,0.000000,0.000000);
		rgb(915pt)=(0.226667,0.000000,0.000000);
		rgb(916pt)=(0.224000,0.000000,0.000000);
		rgb(917pt)=(0.221333,0.000000,0.000000);
		rgb(918pt)=(0.218667,0.000000,0.000000);
		rgb(919pt)=(0.216000,0.000000,0.000000);
		rgb(920pt)=(0.213333,0.000000,0.000000);
		rgb(921pt)=(0.210667,0.000000,0.000000);
		rgb(922pt)=(0.208000,0.000000,0.000000);
		rgb(923pt)=(0.205333,0.000000,0.000000);
		rgb(924pt)=(0.202667,0.000000,0.000000);
		rgb(925pt)=(0.200000,0.000000,0.000000);
		rgb(926pt)=(0.197333,0.000000,0.000000);
		rgb(927pt)=(0.194667,0.000000,0.000000);
		rgb(928pt)=(0.192000,0.000000,0.000000);
		rgb(929pt)=(0.189333,0.000000,0.000000);
		rgb(930pt)=(0.186667,0.000000,0.000000);
		rgb(931pt)=(0.184000,0.000000,0.000000);
		rgb(932pt)=(0.181333,0.000000,0.000000);
		rgb(933pt)=(0.178667,0.000000,0.000000);
		rgb(934pt)=(0.176000,0.000000,0.000000);
		rgb(935pt)=(0.173333,0.000000,0.000000);
		rgb(936pt)=(0.170667,0.000000,0.000000);
		rgb(937pt)=(0.168000,0.000000,0.000000);
		rgb(938pt)=(0.165333,0.000000,0.000000);
		rgb(939pt)=(0.162667,0.000000,0.000000);
		rgb(940pt)=(0.160000,0.000000,0.000000);
		rgb(941pt)=(0.157333,0.000000,0.000000);
		rgb(942pt)=(0.154667,0.000000,0.000000);
		rgb(943pt)=(0.152000,0.000000,0.000000);
		rgb(944pt)=(0.149333,0.000000,0.000000);
		rgb(945pt)=(0.146667,0.000000,0.000000);
		rgb(946pt)=(0.144000,0.000000,0.000000);
		rgb(947pt)=(0.141333,0.000000,0.000000);
		rgb(948pt)=(0.138667,0.000000,0.000000);
		rgb(949pt)=(0.136000,0.000000,0.000000);
		rgb(950pt)=(0.133333,0.000000,0.000000);
		rgb(951pt)=(0.130667,0.000000,0.000000);
		rgb(952pt)=(0.128000,0.000000,0.000000);
		rgb(953pt)=(0.125333,0.000000,0.000000);
		rgb(954pt)=(0.122667,0.000000,0.000000);
		rgb(955pt)=(0.120000,0.000000,0.000000);
		rgb(956pt)=(0.117333,0.000000,0.000000);
		rgb(957pt)=(0.114667,0.000000,0.000000);
		rgb(958pt)=(0.112000,0.000000,0.000000);
		rgb(959pt)=(0.109333,0.000000,0.000000);
		rgb(960pt)=(0.106667,0.000000,0.000000);
		rgb(961pt)=(0.104000,0.000000,0.000000);
		rgb(962pt)=(0.101333,0.000000,0.000000);
		rgb(963pt)=(0.098667,0.000000,0.000000);
		rgb(964pt)=(0.096000,0.000000,0.000000);
		rgb(965pt)=(0.093333,0.000000,0.000000);
		rgb(966pt)=(0.090667,0.000000,0.000000);
		rgb(967pt)=(0.088000,0.000000,0.000000);
		rgb(968pt)=(0.085333,0.000000,0.000000);
		rgb(969pt)=(0.082667,0.000000,0.000000);
		rgb(970pt)=(0.080000,0.000000,0.000000);
		rgb(971pt)=(0.077333,0.000000,0.000000);
		rgb(972pt)=(0.074667,0.000000,0.000000);
		rgb(973pt)=(0.072000,0.000000,0.000000);
		rgb(974pt)=(0.069333,0.000000,0.000000);
		rgb(975pt)=(0.066667,0.000000,0.000000);
		rgb(976pt)=(0.064000,0.000000,0.000000);
		rgb(977pt)=(0.061333,0.000000,0.000000);
		rgb(978pt)=(0.058667,0.000000,0.000000);
		rgb(979pt)=(0.056000,0.000000,0.000000);
		rgb(980pt)=(0.053333,0.000000,0.000000);
		rgb(981pt)=(0.050667,0.000000,0.000000);
		rgb(982pt)=(0.048000,0.000000,0.000000);
		rgb(983pt)=(0.045333,0.000000,0.000000);
		rgb(984pt)=(0.042667,0.000000,0.000000);
		rgb(985pt)=(0.040000,0.000000,0.000000);
		rgb(986pt)=(0.037333,0.000000,0.000000);
		rgb(987pt)=(0.034667,0.000000,0.000000);
		rgb(988pt)=(0.032000,0.000000,0.000000);
		rgb(989pt)=(0.029333,0.000000,0.000000);
		rgb(990pt)=(0.026667,0.000000,0.000000);
		rgb(991pt)=(0.024000,0.000000,0.000000);
		rgb(992pt)=(0.021333,0.000000,0.000000);
		rgb(993pt)=(0.018667,0.000000,0.000000);
		rgb(994pt)=(0.016000,0.000000,0.000000);
		rgb(995pt)=(0.013333,0.000000,0.000000);
		rgb(996pt)=(0.010667,0.000000,0.000000);
		rgb(997pt)=(0.008000,0.000000,0.000000);
		rgb(998pt)=(0.005333,0.000000,0.000000);
		rgb(999pt)=(0.002667,0.000000,0.000000);
}}

\newlength{\figureheight}
\newlength{\figurewidth}
\pgfplotsset{
	colormap={parula}{
		rgb(0pt)=(0.2081,0.1663,0.5292);
		rgb(1pt)=(0.208355,0.16778,0.532238);
		rgb(2pt)=(0.208611,0.169261,0.535275);
		rgb(3pt)=(0.208866,0.170741,0.538313);
		rgb(4pt)=(0.209121,0.172222,0.54135);
		rgb(5pt)=(0.209376,0.173702,0.544388);
		rgb(6pt)=(0.209632,0.175183,0.547425);
		rgb(7pt)=(0.209887,0.176663,0.550463);
		rgb(8pt)=(0.210134,0.178144,0.553505);
		rgb(9pt)=(0.210338,0.179624,0.556568);
		rgb(10pt)=(0.210542,0.181105,0.559631);
		rgb(11pt)=(0.210746,0.182585,0.562694);
		rgb(12pt)=(0.210944,0.184066,0.565763);
		rgb(13pt)=(0.211123,0.185546,0.568852);
		rgb(14pt)=(0.211302,0.187027,0.57194);
		rgb(15pt)=(0.21148,0.188507,0.575029);
		rgb(16pt)=(0.211642,0.189996,0.578117);
		rgb(17pt)=(0.21177,0.191502,0.581206);
		rgb(18pt)=(0.211897,0.193008,0.584295);
		rgb(19pt)=(0.212025,0.194514,0.587383);
		rgb(20pt)=(0.212132,0.19602,0.590472);
		rgb(21pt)=(0.212208,0.197526,0.59356);
		rgb(22pt)=(0.212285,0.199032,0.596649);
		rgb(23pt)=(0.212361,0.200538,0.599738);
		rgb(24pt)=(0.212413,0.202044,0.602839);
		rgb(25pt)=(0.212438,0.20355,0.605953);
		rgb(26pt)=(0.212464,0.205056,0.609067);
		rgb(27pt)=(0.212489,0.206562,0.612181);
		rgb(28pt)=(0.212471,0.208083,0.61531);
		rgb(29pt)=(0.21242,0.209614,0.61845);
		rgb(30pt)=(0.212368,0.211146,0.621589);
		rgb(31pt)=(0.212317,0.212677,0.624729);
		rgb(32pt)=(0.212216,0.214209,0.627868);
		rgb(33pt)=(0.212088,0.215741,0.631008);
		rgb(34pt)=(0.211961,0.217272,0.634148);
		rgb(35pt)=(0.211833,0.218804,0.637287);
		rgb(36pt)=(0.211668,0.220354,0.640446);
		rgb(37pt)=(0.211489,0.221911,0.643611);
		rgb(38pt)=(0.21131,0.223468,0.646776);
		rgb(39pt)=(0.211132,0.225025,0.649941);
		rgb(40pt)=(0.210848,0.226603,0.653107);
		rgb(41pt)=(0.210541,0.228186,0.656272);
		rgb(42pt)=(0.210235,0.229768,0.659437);
		rgb(43pt)=(0.209929,0.231351,0.662602);
		rgb(44pt)=(0.209553,0.232934,0.665767);
		rgb(45pt)=(0.20917,0.234516,0.668932);
		rgb(46pt)=(0.208787,0.236099,0.672098);
		rgb(47pt)=(0.208405,0.237681,0.675263);
		rgb(48pt)=(0.20787,0.239289,0.678453);
		rgb(49pt)=(0.207334,0.240897,0.681644);
		rgb(50pt)=(0.206798,0.242505,0.684835);
		rgb(51pt)=(0.206255,0.244114,0.688025);
		rgb(52pt)=(0.205617,0.245722,0.691216);
		rgb(53pt)=(0.204979,0.24733,0.694407);
		rgb(54pt)=(0.204341,0.248938,0.697597);
		rgb(55pt)=(0.203675,0.250554,0.700792);
		rgb(56pt)=(0.202858,0.252213,0.704008);
		rgb(57pt)=(0.202041,0.253872,0.707224);
		rgb(58pt)=(0.201225,0.255531,0.710441);
		rgb(59pt)=(0.200372,0.257184,0.713657);
		rgb(60pt)=(0.199402,0.258818,0.716873);
		rgb(61pt)=(0.198432,0.260452,0.720089);
		rgb(62pt)=(0.197462,0.262085,0.723305);
		rgb(63pt)=(0.196419,0.263735,0.726522);
		rgb(64pt)=(0.195219,0.26542,0.729738);
		rgb(65pt)=(0.19402,0.267105,0.732954);
		rgb(66pt)=(0.19282,0.268789,0.73617);
		rgb(67pt)=(0.191549,0.270474,0.739386);
		rgb(68pt)=(0.19017,0.272159,0.742603);
		rgb(69pt)=(0.188792,0.273843,0.745819);
		rgb(70pt)=(0.187414,0.275528,0.749035);
		rgb(71pt)=(0.1859,0.277237,0.752264);
		rgb(72pt)=(0.184241,0.278973,0.755505);
		rgb(73pt)=(0.182581,0.280709,0.758747);
		rgb(74pt)=(0.180922,0.282444,0.761989);
		rgb(75pt)=(0.179133,0.284209,0.765245);
		rgb(76pt)=(0.177244,0.285996,0.768512);
		rgb(77pt)=(0.175356,0.287783,0.77178);
		rgb(78pt)=(0.173467,0.289569,0.775047);
		rgb(79pt)=(0.171363,0.291406,0.778314);
		rgb(80pt)=(0.169142,0.293269,0.781581);
		rgb(81pt)=(0.166922,0.295132,0.784849);
		rgb(82pt)=(0.164701,0.296996,0.788116);
		rgb(83pt)=(0.162238,0.298934,0.791365);
		rgb(84pt)=(0.159686,0.300899,0.794606);
		rgb(85pt)=(0.157133,0.302865,0.797848);
		rgb(86pt)=(0.15458,0.30483,0.80109);
		rgb(87pt)=(0.151738,0.306858,0.804352);
		rgb(88pt)=(0.148828,0.3089,0.80762);
		rgb(89pt)=(0.145918,0.310942,0.810887);
		rgb(90pt)=(0.143008,0.312984,0.814154);
		rgb(91pt)=(0.139687,0.31514,0.81733);
		rgb(92pt)=(0.136318,0.31731,0.820495);
		rgb(93pt)=(0.132949,0.319479,0.82366);
		rgb(94pt)=(0.129579,0.321649,0.826826);
		rgb(95pt)=(0.125811,0.323918,0.829841);
		rgb(96pt)=(0.122033,0.32619,0.832853);
		rgb(97pt)=(0.118256,0.328462,0.835865);
		rgb(98pt)=(0.114458,0.330737,0.838862);
		rgb(99pt)=(0.110349,0.333059,0.841619);
		rgb(100pt)=(0.106239,0.335382,0.844376);
		rgb(101pt)=(0.102129,0.337705,0.847132);
		rgb(102pt)=(0.0979874,0.340021,0.849835);
		rgb(103pt)=(0.093648,0.342292,0.852209);
		rgb(104pt)=(0.0893087,0.344564,0.854583);
		rgb(105pt)=(0.0849694,0.346836,0.856957);
		rgb(106pt)=(0.08063,0.349091,0.859234);
		rgb(107pt)=(0.0762907,0.351286,0.861174);
		rgb(108pt)=(0.0719514,0.353481,0.863114);
		rgb(109pt)=(0.067612,0.355676,0.865053);
		rgb(110pt)=(0.0633195,0.357817,0.866853);
		rgb(111pt)=(0.0591333,0.359833,0.868333);
		rgb(112pt)=(0.0549471,0.36185,0.869814);
		rgb(113pt)=(0.050761,0.363866,0.871294);
		rgb(114pt)=(0.0466838,0.365823,0.872626);
		rgb(115pt)=(0.0427784,0.367687,0.873724);
		rgb(116pt)=(0.038873,0.36955,0.874821);
		rgb(117pt)=(0.0349676,0.371414,0.875919);
		rgb(118pt)=(0.0315066,0.373217,0.876872);
		rgb(119pt)=(0.0285456,0.374953,0.877664);
		rgb(120pt)=(0.0255847,0.376688,0.878455);
		rgb(121pt)=(0.0226237,0.378424,0.879246);
		rgb(122pt)=(0.0202132,0.380061,0.879868);
		rgb(123pt)=(0.0182477,0.381618,0.880353);
		rgb(124pt)=(0.0162823,0.383175,0.880838);
		rgb(125pt)=(0.0143168,0.384732,0.881323);
		rgb(126pt)=(0.0127892,0.386241,0.881695);
		rgb(127pt)=(0.0115129,0.387721,0.882001);
		rgb(128pt)=(0.0102366,0.389202,0.882307);
		rgb(129pt)=(0.00896036,0.390682,0.882614);
		rgb(130pt)=(0.00812372,0.392089,0.88281);
		rgb(131pt)=(0.00746006,0.393468,0.882963);
		rgb(132pt)=(0.0067964,0.394846,0.883116);
		rgb(133pt)=(0.00613273,0.396224,0.883269);
		rgb(134pt)=(0.00581622,0.397562,0.88332);
		rgb(135pt)=(0.00558649,0.398889,0.883346);
		rgb(136pt)=(0.00535676,0.400217,0.883371);
		rgb(137pt)=(0.00512703,0.401544,0.883397);
		rgb(138pt)=(0.00516757,0.402804,0.883332);
		rgb(139pt)=(0.00524414,0.404054,0.883256);
		rgb(140pt)=(0.00532072,0.405305,0.883179);
		rgb(141pt)=(0.0053973,0.406556,0.883103);
		rgb(142pt)=(0.00572012,0.407757,0.882952);
		rgb(143pt)=(0.00605195,0.408957,0.882799);
		rgb(144pt)=(0.00638378,0.410157,0.882646);
		rgb(145pt)=(0.00672643,0.411355,0.882489);
		rgb(146pt)=(0.00728799,0.412529,0.882259);
		rgb(147pt)=(0.00784955,0.413704,0.88203);
		rgb(148pt)=(0.00841111,0.414878,0.8818);
		rgb(149pt)=(0.00898919,0.416045,0.881564);
		rgb(150pt)=(0.00967838,0.417168,0.881283);
		rgb(151pt)=(0.0103676,0.418292,0.881002);
		rgb(152pt)=(0.0110568,0.419415,0.880721);
		rgb(153pt)=(0.011773,0.420532,0.880435);
		rgb(154pt)=(0.0125898,0.42163,0.880129);
		rgb(155pt)=(0.0134066,0.422728,0.879823);
		rgb(156pt)=(0.0142234,0.423825,0.879516);
		rgb(157pt)=(0.0150703,0.424915,0.879195);
		rgb(158pt)=(0.0159892,0.425987,0.878838);
		rgb(159pt)=(0.0169081,0.427059,0.87848);
		rgb(160pt)=(0.017827,0.428132,0.878123);
		rgb(161pt)=(0.0187748,0.429194,0.877746);
		rgb(162pt)=(0.0197703,0.430241,0.877338);
		rgb(163pt)=(0.0207658,0.431287,0.876929);
		rgb(164pt)=(0.0217613,0.432334,0.876521);
		rgb(165pt)=(0.0227802,0.43338,0.876113);
		rgb(166pt)=(0.0238267,0.434427,0.875704);
		rgb(167pt)=(0.0248733,0.435473,0.875296);
		rgb(168pt)=(0.0259198,0.43652,0.874887);
		rgb(169pt)=(0.0269802,0.437553,0.874451);
		rgb(170pt)=(0.0280523,0.438574,0.873992);
		rgb(171pt)=(0.0291243,0.439595,0.873532);
		rgb(172pt)=(0.0301964,0.440616,0.873073);
		rgb(173pt)=(0.0312844,0.441621,0.872614);
		rgb(174pt)=(0.032382,0.442616,0.872154);
		rgb(175pt)=(0.0334796,0.443612,0.871695);
		rgb(176pt)=(0.0345772,0.444607,0.871235);
		rgb(177pt)=(0.0357108,0.445603,0.870758);
		rgb(178pt)=(0.0368595,0.446598,0.870273);
		rgb(179pt)=(0.0380081,0.447594,0.869788);
		rgb(180pt)=(0.0391568,0.448589,0.869303);
		rgb(181pt)=(0.0402652,0.449565,0.868798);
		rgb(182pt)=(0.0413628,0.450535,0.868287);
		rgb(183pt)=(0.0424604,0.451505,0.867777);
		rgb(184pt)=(0.043558,0.452474,0.867266);
		rgb(185pt)=(0.0445889,0.453444,0.866756);
		rgb(186pt)=(0.0456099,0.454414,0.866245);
		rgb(187pt)=(0.0466309,0.455384,0.865735);
		rgb(188pt)=(0.047652,0.456354,0.865224);
		rgb(189pt)=(0.0486,0.457324,0.864714);
		rgb(190pt)=(0.0495444,0.458294,0.864203);
		rgb(191pt)=(0.0504889,0.459264,0.863692);
		rgb(192pt)=(0.0514315,0.460234,0.863181);
		rgb(193pt)=(0.0523249,0.461204,0.862645);
		rgb(194pt)=(0.0532183,0.462174,0.862109);
		rgb(195pt)=(0.0541117,0.463144,0.861573);
		rgb(196pt)=(0.0549991,0.464111,0.861034);
		rgb(197pt)=(0.0558414,0.465056,0.860472);
		rgb(198pt)=(0.0566838,0.466,0.859911);
		rgb(199pt)=(0.0575261,0.466944,0.859349);
		rgb(200pt)=(0.0583532,0.467889,0.858793);
		rgb(201pt)=(0.0591189,0.468833,0.858257);
		rgb(202pt)=(0.0598847,0.469778,0.857721);
		rgb(203pt)=(0.0606505,0.470722,0.857185);
		rgb(204pt)=(0.0614018,0.471667,0.856641);
		rgb(205pt)=(0.0621165,0.472611,0.85608);
		rgb(206pt)=(0.0628312,0.473556,0.855518);
		rgb(207pt)=(0.0635459,0.4745,0.854957);
		rgb(208pt)=(0.064242,0.475444,0.854405);
		rgb(209pt)=(0.0649057,0.476389,0.853868);
		rgb(210pt)=(0.0655694,0.477333,0.853332);
		rgb(211pt)=(0.066233,0.478278,0.852796);
		rgb(212pt)=(0.0668625,0.479222,0.852249);
		rgb(213pt)=(0.0674495,0.480167,0.851687);
		rgb(214pt)=(0.0680366,0.481111,0.851126);
		rgb(215pt)=(0.0686237,0.482056,0.850564);
		rgb(216pt)=(0.0691838,0.483,0.850003);
		rgb(217pt)=(0.0697198,0.483944,0.849441);
		rgb(218pt)=(0.0702559,0.484889,0.84888);
		rgb(219pt)=(0.0707919,0.485833,0.848318);
		rgb(220pt)=(0.0712967,0.486778,0.847772);
		rgb(221pt)=(0.0717817,0.487722,0.847236);
		rgb(222pt)=(0.0722667,0.488667,0.8467);
		rgb(223pt)=(0.0727517,0.489611,0.846164);
		rgb(224pt)=(0.0732012,0.490573,0.845628);
		rgb(225pt)=(0.0736351,0.491543,0.845092);
		rgb(226pt)=(0.0740691,0.492513,0.844556);
		rgb(227pt)=(0.074503,0.493483,0.84402);
		rgb(228pt)=(0.0748973,0.494433,0.843484);
		rgb(229pt)=(0.0752802,0.495378,0.842948);
		rgb(230pt)=(0.0756631,0.496322,0.842412);
		rgb(231pt)=(0.0760459,0.497267,0.841876);
		rgb(232pt)=(0.0763631,0.498233,0.841362);
		rgb(233pt)=(0.0766694,0.499203,0.840851);
		rgb(234pt)=(0.0769757,0.500173,0.840341);
		rgb(235pt)=(0.077282,0.501143,0.83983);
		rgb(236pt)=(0.0775162,0.502137,0.83932);
		rgb(237pt)=(0.0777459,0.503132,0.838809);
		rgb(238pt)=(0.0779757,0.504128,0.838298);
		rgb(239pt)=(0.0782042,0.505123,0.837789);
		rgb(240pt)=(0.0783829,0.506093,0.837304);
		rgb(241pt)=(0.0785616,0.507063,0.836819);
		rgb(242pt)=(0.0787402,0.508033,0.836334);
		rgb(243pt)=(0.0789135,0.509008,0.835851);
		rgb(244pt)=(0.0790411,0.510029,0.835392);
		rgb(245pt)=(0.0791688,0.51105,0.834932);
		rgb(246pt)=(0.0792964,0.512071,0.834473);
		rgb(247pt)=(0.0794048,0.513092,0.834018);
		rgb(248pt)=(0.0794303,0.514113,0.833584);
		rgb(249pt)=(0.0794559,0.515134,0.83315);
		rgb(250pt)=(0.0794814,0.516155,0.832717);
		rgb(251pt)=(0.0794862,0.517183,0.832289);
		rgb(252pt)=(0.0794351,0.51823,0.831881);
		rgb(253pt)=(0.0793841,0.519276,0.831473);
		rgb(254pt)=(0.079333,0.520323,0.831064);
		rgb(255pt)=(0.079255,0.521369,0.830665);
		rgb(256pt)=(0.0791273,0.522416,0.830282);
		rgb(257pt)=(0.0789997,0.523462,0.829899);
		rgb(258pt)=(0.0788721,0.524509,0.829516);
		rgb(259pt)=(0.0786889,0.525589,0.829156);
		rgb(260pt)=(0.0784336,0.526712,0.828824);
		rgb(261pt)=(0.0781784,0.527835,0.828492);
		rgb(262pt)=(0.0779231,0.528958,0.82816);
		rgb(263pt)=(0.077615,0.530081,0.827868);
		rgb(264pt)=(0.0772577,0.531205,0.827613);
		rgb(265pt)=(0.0769003,0.532328,0.827357);
		rgb(266pt)=(0.0765429,0.533451,0.827102);
		rgb(267pt)=(0.0761243,0.534589,0.826862);
		rgb(268pt)=(0.0756649,0.535738,0.826632);
		rgb(269pt)=(0.0752054,0.536886,0.826403);
		rgb(270pt)=(0.0747459,0.538035,0.826173);
		rgb(271pt)=(0.0742168,0.539219,0.825961);
		rgb(272pt)=(0.0736553,0.540418,0.825756);
		rgb(273pt)=(0.0730937,0.541618,0.825552);
		rgb(274pt)=(0.0725321,0.542818,0.825348);
		rgb(275pt)=(0.0718925,0.544037,0.825183);
		rgb(276pt)=(0.0712288,0.545262,0.82503);
		rgb(277pt)=(0.0705652,0.546487,0.824877);
		rgb(278pt)=(0.0699015,0.547713,0.824723);
		rgb(279pt)=(0.0691514,0.548938,0.824614);
		rgb(280pt)=(0.0683856,0.550163,0.824511);
		rgb(281pt)=(0.0676198,0.551388,0.824409);
		rgb(282pt)=(0.0668541,0.552614,0.824307);
		rgb(283pt)=(0.0660408,0.553886,0.824205);
		rgb(284pt)=(0.065224,0.555162,0.824103);
		rgb(285pt)=(0.0644072,0.556439,0.824001);
		rgb(286pt)=(0.0635892,0.557715,0.823899);
		rgb(287pt)=(0.0626703,0.558991,0.823848);
		rgb(288pt)=(0.0617514,0.560268,0.823797);
		rgb(289pt)=(0.0608324,0.561544,0.823746);
		rgb(290pt)=(0.0599087,0.56282,0.823693);
		rgb(291pt)=(0.0589387,0.564096,0.823616);
		rgb(292pt)=(0.0579688,0.565373,0.82354);
		rgb(293pt)=(0.0569988,0.566649,0.823463);
		rgb(294pt)=(0.0560243,0.567925,0.823386);
		rgb(295pt)=(0.0550288,0.569202,0.82331);
		rgb(296pt)=(0.0540333,0.570478,0.823233);
		rgb(297pt)=(0.0530378,0.571754,0.823157);
		rgb(298pt)=(0.0520423,0.57303,0.82308);
		rgb(299pt)=(0.0510468,0.574307,0.823004);
		rgb(300pt)=(0.0500514,0.575583,0.822927);
		rgb(301pt)=(0.0490559,0.576859,0.82285);
		rgb(302pt)=(0.0480604,0.578127,0.822756);
		rgb(303pt)=(0.0470649,0.579377,0.822629);
		rgb(304pt)=(0.0460694,0.580628,0.822501);
		rgb(305pt)=(0.0450739,0.581879,0.822374);
		rgb(306pt)=(0.0441,0.583119,0.822235);
		rgb(307pt)=(0.0431556,0.584344,0.822082);
		rgb(308pt)=(0.0422111,0.585569,0.821929);
		rgb(309pt)=(0.0412667,0.586795,0.821776);
		rgb(310pt)=(0.0403351,0.58802,0.821597);
		rgb(311pt)=(0.0394162,0.589245,0.821392);
		rgb(312pt)=(0.0384973,0.59047,0.821188);
		rgb(313pt)=(0.0375784,0.591695,0.820984);
		rgb(314pt)=(0.0367495,0.592891,0.820735);
		rgb(315pt)=(0.0359838,0.594065,0.820454);
		rgb(316pt)=(0.035218,0.595239,0.820173);
		rgb(317pt)=(0.0344523,0.596413,0.819892);
		rgb(318pt)=(0.0337721,0.597553,0.819595);
		rgb(319pt)=(0.0331339,0.598676,0.819288);
		rgb(320pt)=(0.0324958,0.599799,0.818982);
		rgb(321pt)=(0.0318577,0.600923,0.818676);
		rgb(322pt)=(0.0312964,0.602026,0.818312);
		rgb(323pt)=(0.0307604,0.603124,0.817929);
		rgb(324pt)=(0.0302243,0.604222,0.817546);
		rgb(325pt)=(0.0296883,0.605319,0.817163);
		rgb(326pt)=(0.0292375,0.606395,0.816738);
		rgb(327pt)=(0.0288036,0.607468,0.816304);
		rgb(328pt)=(0.0283697,0.60854,0.81587);
		rgb(329pt)=(0.0279357,0.609612,0.815436);
		rgb(330pt)=(0.0275721,0.610637,0.814955);
		rgb(331pt)=(0.0272147,0.611658,0.81447);
		rgb(332pt)=(0.0268574,0.612679,0.813985);
		rgb(333pt)=(0.0265,0.6137,0.8135);
		rgb(334pt)=(0.0262447,0.614695,0.812964);
		rgb(335pt)=(0.0259895,0.615691,0.812428);
		rgb(336pt)=(0.0257342,0.616686,0.811892);
		rgb(337pt)=(0.0254853,0.61768,0.811352);
		rgb(338pt)=(0.0253066,0.61865,0.810765);
		rgb(339pt)=(0.0251279,0.61962,0.810177);
		rgb(340pt)=(0.0249492,0.62059,0.80959);
		rgb(341pt)=(0.024779,0.621551,0.808995);
		rgb(342pt)=(0.0246514,0.62247,0.808357);
		rgb(343pt)=(0.0245237,0.623389,0.807719);
		rgb(344pt)=(0.0243961,0.624308,0.80708);
		rgb(345pt)=(0.0242748,0.625221,0.80643);
		rgb(346pt)=(0.0241727,0.626114,0.805741);
		rgb(347pt)=(0.0240706,0.627008,0.805051);
		rgb(348pt)=(0.0239685,0.627901,0.804362);
		rgb(349pt)=(0.0238832,0.628786,0.803656);
		rgb(350pt)=(0.0238321,0.629654,0.802916);
		rgb(351pt)=(0.0237811,0.630522,0.802176);
		rgb(352pt)=(0.02373,0.631389,0.801435);
		rgb(353pt)=(0.023679,0.632247,0.800685);
		rgb(354pt)=(0.0236279,0.633089,0.799919);
		rgb(355pt)=(0.0235769,0.633932,0.799153);
		rgb(356pt)=(0.0235258,0.634774,0.798387);
		rgb(357pt)=(0.0234748,0.635604,0.797596);
		rgb(358pt)=(0.0234237,0.63642,0.79678);
		rgb(359pt)=(0.0233727,0.637237,0.795963);
		rgb(360pt)=(0.0233216,0.638054,0.795146);
		rgb(361pt)=(0.0232706,0.638856,0.794329);
		rgb(362pt)=(0.0232195,0.639647,0.793512);
		rgb(363pt)=(0.0231685,0.640439,0.792695);
		rgb(364pt)=(0.0231174,0.64123,0.791879);
		rgb(365pt)=(0.0230832,0.642005,0.791011);
		rgb(366pt)=(0.0230577,0.64277,0.790118);
		rgb(367pt)=(0.0230321,0.643536,0.789225);
		rgb(368pt)=(0.0230066,0.644302,0.788331);
		rgb(369pt)=(0.0229811,0.645049,0.787438);
		rgb(370pt)=(0.0229556,0.645789,0.786544);
		rgb(371pt)=(0.02293,0.646529,0.785651);
		rgb(372pt)=(0.0229045,0.647269,0.784758);
		rgb(373pt)=(0.022858,0.64801,0.783843);
		rgb(374pt)=(0.0228069,0.64875,0.782924);
		rgb(375pt)=(0.0227559,0.64949,0.782005);
		rgb(376pt)=(0.0227048,0.65023,0.781086);
		rgb(377pt)=(0.0227,0.650947,0.780144);
		rgb(378pt)=(0.0227,0.651662,0.7792);
		rgb(379pt)=(0.0227,0.652377,0.778256);
		rgb(380pt)=(0.0227,0.653092,0.777311);
		rgb(381pt)=(0.0228261,0.653781,0.776341);
		rgb(382pt)=(0.0229538,0.65447,0.775371);
		rgb(383pt)=(0.0230814,0.655159,0.774402);
		rgb(384pt)=(0.0232108,0.655849,0.77343);
		rgb(385pt)=(0.023364,0.656538,0.772434);
		rgb(386pt)=(0.0235171,0.657227,0.771439);
		rgb(387pt)=(0.0236703,0.657916,0.770443);
		rgb(388pt)=(0.0238312,0.658602,0.769444);
		rgb(389pt)=(0.0240354,0.659265,0.768423);
		rgb(390pt)=(0.0242396,0.659929,0.767402);
		rgb(391pt)=(0.0244438,0.660592,0.766381);
		rgb(392pt)=(0.0247021,0.661256,0.765354);
		rgb(393pt)=(0.025136,0.66192,0.764307);
		rgb(394pt)=(0.02557,0.662583,0.763261);
		rgb(395pt)=(0.0260039,0.663247,0.762214);
		rgb(396pt)=(0.0264541,0.663911,0.761168);
		rgb(397pt)=(0.026939,0.664574,0.760121);
		rgb(398pt)=(0.027424,0.665238,0.759074);
		rgb(399pt)=(0.027909,0.665902,0.758028);
		rgb(400pt)=(0.028445,0.666555,0.756971);
		rgb(401pt)=(0.0290577,0.667193,0.755899);
		rgb(402pt)=(0.0296703,0.667832,0.754827);
		rgb(403pt)=(0.0302829,0.66847,0.753755);
		rgb(404pt)=(0.030994,0.669095,0.752683);
		rgb(405pt)=(0.0318108,0.669708,0.751611);
		rgb(406pt)=(0.0326276,0.670321,0.750539);
		rgb(407pt)=(0.0334444,0.670933,0.749467);
		rgb(408pt)=(0.0343045,0.67156,0.748366);
		rgb(409pt)=(0.0351979,0.672198,0.747243);
		rgb(410pt)=(0.0360913,0.672837,0.74612);
		rgb(411pt)=(0.0369847,0.673475,0.744996);
		rgb(412pt)=(0.0380432,0.674096,0.743873);
		rgb(413pt)=(0.0391919,0.674709,0.74275);
		rgb(414pt)=(0.0403405,0.675322,0.741627);
		rgb(415pt)=(0.0414892,0.675934,0.740504);
		rgb(416pt)=(0.0427123,0.676528,0.739381);
		rgb(417pt)=(0.0439631,0.677115,0.738258);
		rgb(418pt)=(0.0452138,0.677702,0.737135);
		rgb(419pt)=(0.0464646,0.678289,0.736011);
		rgb(420pt)=(0.0477153,0.678897,0.734868);
		rgb(421pt)=(0.0489661,0.67951,0.733719);
		rgb(422pt)=(0.0502168,0.680123,0.73257);
		rgb(423pt)=(0.0514676,0.680735,0.731422);
		rgb(424pt)=(0.0529237,0.681325,0.73025);
		rgb(425pt)=(0.0544042,0.681912,0.729076);
		rgb(426pt)=(0.0558847,0.682499,0.727902);
		rgb(427pt)=(0.0573652,0.683086,0.726728);
		rgb(428pt)=(0.0587709,0.683673,0.725553);
		rgb(429pt)=(0.0601748,0.68426,0.724379);
		rgb(430pt)=(0.0615787,0.684847,0.723205);
		rgb(431pt)=(0.0629946,0.685435,0.722028);
		rgb(432pt)=(0.0646027,0.686022,0.720803);
		rgb(433pt)=(0.0662108,0.686609,0.719577);
		rgb(434pt)=(0.0678189,0.687196,0.718352);
		rgb(435pt)=(0.069427,0.687779,0.717131);
		rgb(436pt)=(0.0710351,0.688341,0.715931);
		rgb(437pt)=(0.0726432,0.688902,0.714731);
		rgb(438pt)=(0.0742514,0.689464,0.713532);
		rgb(439pt)=(0.0758709,0.690026,0.712326);
		rgb(440pt)=(0.07753,0.690587,0.711101);
		rgb(441pt)=(0.0791892,0.691149,0.709876);
		rgb(442pt)=(0.0808483,0.69171,0.70865);
		rgb(443pt)=(0.0825387,0.692272,0.707417);
		rgb(444pt)=(0.0843,0.692833,0.706167);
		rgb(445pt)=(0.0860613,0.693395,0.704916);
		rgb(446pt)=(0.0878225,0.693956,0.703665);
		rgb(447pt)=(0.089564,0.694518,0.702405);
		rgb(448pt)=(0.0912742,0.69508,0.701128);
		rgb(449pt)=(0.0929844,0.695641,0.699852);
		rgb(450pt)=(0.0946946,0.696203,0.698576);
		rgb(451pt)=(0.0965009,0.696752,0.697299);
		rgb(452pt)=(0.0984153,0.697288,0.696023);
		rgb(453pt)=(0.10033,0.697824,0.694747);
		rgb(454pt)=(0.102244,0.69836,0.693471);
		rgb(455pt)=(0.10413,0.698896,0.69218);
		rgb(456pt)=(0.105994,0.699432,0.690878);
		rgb(457pt)=(0.107857,0.699968,0.689577);
		rgb(458pt)=(0.10972,0.700505,0.688275);
		rgb(459pt)=(0.111632,0.701041,0.686973);
		rgb(460pt)=(0.113572,0.701577,0.685671);
		rgb(461pt)=(0.115512,0.702113,0.684369);
		rgb(462pt)=(0.117452,0.702649,0.683068);
		rgb(463pt)=(0.119429,0.703185,0.681747);
		rgb(464pt)=(0.12142,0.703721,0.68042);
		rgb(465pt)=(0.123411,0.704257,0.679093);
		rgb(466pt)=(0.125402,0.704793,0.677765);
		rgb(467pt)=(0.127372,0.705308,0.676438);
		rgb(468pt)=(0.129338,0.705819,0.675111);
		rgb(469pt)=(0.131303,0.706329,0.673783);
		rgb(470pt)=(0.133269,0.70684,0.672456);
		rgb(471pt)=(0.135369,0.70735,0.671084);
		rgb(472pt)=(0.137488,0.707861,0.669705);
		rgb(473pt)=(0.139607,0.708371,0.668327);
		rgb(474pt)=(0.141725,0.708882,0.666949);
		rgb(475pt)=(0.143795,0.709392,0.665595);
		rgb(476pt)=(0.145862,0.709903,0.664242);
		rgb(477pt)=(0.14793,0.710414,0.662889);
		rgb(478pt)=(0.150003,0.710924,0.661534);
		rgb(479pt)=(0.152198,0.711435,0.66013);
		rgb(480pt)=(0.154394,0.711945,0.658726);
		rgb(481pt)=(0.156589,0.712456,0.657322);
		rgb(482pt)=(0.158784,0.712963,0.655922);
		rgb(483pt)=(0.160979,0.713448,0.654543);
		rgb(484pt)=(0.163174,0.713933,0.653165);
		rgb(485pt)=(0.16537,0.714418,0.651786);
		rgb(486pt)=(0.16757,0.714908,0.650397);
		rgb(487pt)=(0.169791,0.715419,0.648968);
		rgb(488pt)=(0.172012,0.715929,0.647538);
		rgb(489pt)=(0.174232,0.71644,0.646109);
		rgb(490pt)=(0.176483,0.716935,0.64468);
		rgb(491pt)=(0.178806,0.717395,0.64325);
		rgb(492pt)=(0.181129,0.717854,0.641821);
		rgb(493pt)=(0.183452,0.718314,0.640391);
		rgb(494pt)=(0.185755,0.718783,0.638952);
		rgb(495pt)=(0.188027,0.719268,0.637497);
		rgb(496pt)=(0.190299,0.719753,0.636042);
		rgb(497pt)=(0.192571,0.720238,0.634587);
		rgb(498pt)=(0.194913,0.720711,0.633132);
		rgb(499pt)=(0.197338,0.72117,0.631677);
		rgb(500pt)=(0.199762,0.72163,0.630223);
		rgb(501pt)=(0.202187,0.722089,0.628768);
		rgb(502pt)=(0.204612,0.722549,0.627299);
		rgb(503pt)=(0.207037,0.723008,0.625818);
		rgb(504pt)=(0.209462,0.723468,0.624338);
		rgb(505pt)=(0.211887,0.723927,0.622857);
		rgb(506pt)=(0.214328,0.724386,0.621377);
		rgb(507pt)=(0.216778,0.724846,0.619896);
		rgb(508pt)=(0.219229,0.725305,0.618416);
		rgb(509pt)=(0.221679,0.725765,0.616935);
		rgb(510pt)=(0.224202,0.726188,0.615455);
		rgb(511pt)=(0.226754,0.726597,0.613974);
		rgb(512pt)=(0.229307,0.727005,0.612494);
		rgb(513pt)=(0.231859,0.727414,0.611014);
		rgb(514pt)=(0.234392,0.727842,0.609513);
		rgb(515pt)=(0.236919,0.728276,0.608007);
		rgb(516pt)=(0.239446,0.72871,0.606501);
		rgb(517pt)=(0.241973,0.729144,0.604995);
		rgb(518pt)=(0.244611,0.729556,0.603467);
		rgb(519pt)=(0.247266,0.729964,0.601935);
		rgb(520pt)=(0.24992,0.730372,0.600404);
		rgb(521pt)=(0.252575,0.730781,0.598872);
		rgb(522pt)=(0.25523,0.731189,0.597365);
		rgb(523pt)=(0.257884,0.731598,0.595859);
		rgb(524pt)=(0.260539,0.732006,0.594353);
		rgb(525pt)=(0.263194,0.732414,0.592846);
		rgb(526pt)=(0.265848,0.732796,0.591314);
		rgb(527pt)=(0.268503,0.733179,0.589783);
		rgb(528pt)=(0.271158,0.733562,0.588251);
		rgb(529pt)=(0.27383,0.733945,0.58672);
		rgb(530pt)=(0.276638,0.734328,0.585188);
		rgb(531pt)=(0.279446,0.734711,0.583657);
		rgb(532pt)=(0.282254,0.735094,0.582125);
		rgb(533pt)=(0.285051,0.735471,0.580594);
		rgb(534pt)=(0.287808,0.735829,0.579062);
		rgb(535pt)=(0.290565,0.736186,0.577531);
		rgb(536pt)=(0.293322,0.736544,0.575999);
		rgb(537pt)=(0.2961,0.736894,0.574468);
		rgb(538pt)=(0.298933,0.737226,0.572936);
		rgb(539pt)=(0.301767,0.737557,0.571405);
		rgb(540pt)=(0.3046,0.737889,0.569873);
		rgb(541pt)=(0.307452,0.738221,0.568351);
		rgb(542pt)=(0.310336,0.738553,0.566845);
		rgb(543pt)=(0.313221,0.738885,0.565339);
		rgb(544pt)=(0.316105,0.739217,0.563833);
		rgb(545pt)=(0.318978,0.739537,0.562315);
		rgb(546pt)=(0.321837,0.739843,0.560784);
		rgb(547pt)=(0.324696,0.74015,0.559252);
		rgb(548pt)=(0.327555,0.740456,0.557721);
		rgb(549pt)=(0.330468,0.740749,0.556216);
		rgb(550pt)=(0.333429,0.741029,0.554736);
		rgb(551pt)=(0.336389,0.74131,0.553255);
		rgb(552pt)=(0.33935,0.741591,0.551775);
		rgb(553pt)=(0.342296,0.741872,0.550279);
		rgb(554pt)=(0.345231,0.742153,0.548773);
		rgb(555pt)=(0.348167,0.742433,0.547267);
		rgb(556pt)=(0.351102,0.742714,0.545761);
		rgb(557pt)=(0.354038,0.742977,0.54429);
		rgb(558pt)=(0.356973,0.743232,0.542835);
		rgb(559pt)=(0.359908,0.743488,0.54138);
		rgb(560pt)=(0.362844,0.743743,0.539925);
		rgb(561pt)=(0.365839,0.743959,0.53847);
		rgb(562pt)=(0.368851,0.744163,0.537015);
		rgb(563pt)=(0.371863,0.744367,0.53556);
		rgb(564pt)=(0.374875,0.744571,0.534105);
		rgb(565pt)=(0.377843,0.744775,0.532672);
		rgb(566pt)=(0.380804,0.74498,0.531243);
		rgb(567pt)=(0.383765,0.745184,0.529814);
		rgb(568pt)=(0.386726,0.745388,0.528384);
		rgb(569pt)=(0.389711,0.745568,0.527003);
		rgb(570pt)=(0.392697,0.745747,0.525624);
		rgb(571pt)=(0.395684,0.745926,0.524246);
		rgb(572pt)=(0.39867,0.746104,0.522868);
		rgb(573pt)=(0.401657,0.746257,0.521489);
		rgb(574pt)=(0.404643,0.74641,0.520111);
		rgb(575pt)=(0.40763,0.746563,0.518732);
		rgb(576pt)=(0.410611,0.746716,0.517359);
		rgb(577pt)=(0.413546,0.746869,0.516032);
		rgb(578pt)=(0.416482,0.747023,0.514705);
		rgb(579pt)=(0.419417,0.747176,0.513377);
		rgb(580pt)=(0.422357,0.747319,0.512055);
		rgb(581pt)=(0.425318,0.747421,0.510753);
		rgb(582pt)=(0.428279,0.747523,0.509451);
		rgb(583pt)=(0.43124,0.747626,0.50815);
		rgb(584pt)=(0.43418,0.747735,0.506848);
		rgb(585pt)=(0.437065,0.747862,0.505546);
		rgb(586pt)=(0.439949,0.74799,0.504244);
		rgb(587pt)=(0.442834,0.748117,0.502942);
		rgb(588pt)=(0.445727,0.748227,0.501659);
		rgb(589pt)=(0.448637,0.748304,0.500408);
		rgb(590pt)=(0.451547,0.74838,0.499157);
		rgb(591pt)=(0.454457,0.748457,0.497906);
		rgb(592pt)=(0.457333,0.748522,0.496667);
		rgb(593pt)=(0.460167,0.748573,0.495441);
		rgb(594pt)=(0.463,0.748624,0.494216);
		rgb(595pt)=(0.465833,0.748675,0.492991);
		rgb(596pt)=(0.468667,0.748726,0.491779);
		rgb(597pt)=(0.4715,0.748777,0.490579);
		rgb(598pt)=(0.474333,0.748829,0.48938);
		rgb(599pt)=(0.477167,0.74888,0.48818);
		rgb(600pt)=(0.479969,0.748931,0.486995);
		rgb(601pt)=(0.482752,0.748982,0.485821);
		rgb(602pt)=(0.485534,0.749033,0.484647);
		rgb(603pt)=(0.488316,0.749084,0.483473);
		rgb(604pt)=(0.491081,0.7491,0.482316);
		rgb(605pt)=(0.493838,0.7491,0.481168);
		rgb(606pt)=(0.496595,0.7491,0.480019);
		rgb(607pt)=(0.499351,0.7491,0.47887);
		rgb(608pt)=(0.502069,0.74912,0.477722);
		rgb(609pt)=(0.504775,0.749145,0.476573);
		rgb(610pt)=(0.50748,0.749171,0.475424);
		rgb(611pt)=(0.510186,0.749196,0.474276);
		rgb(612pt)=(0.512892,0.7492,0.47317);
		rgb(613pt)=(0.515598,0.7492,0.472073);
		rgb(614pt)=(0.518303,0.7492,0.470975);
		rgb(615pt)=(0.521009,0.7492,0.469877);
		rgb(616pt)=(0.523644,0.749176,0.46878);
		rgb(617pt)=(0.526273,0.749151,0.467682);
		rgb(618pt)=(0.528902,0.749125,0.466585);
		rgb(619pt)=(0.531531,0.7491,0.465487);
		rgb(620pt)=(0.53416,0.749074,0.464415);
		rgb(621pt)=(0.536789,0.749049,0.463343);
		rgb(622pt)=(0.539418,0.749023,0.462271);
		rgb(623pt)=(0.542043,0.748998,0.461199);
		rgb(624pt)=(0.544621,0.748972,0.460127);
		rgb(625pt)=(0.547199,0.748947,0.459055);
		rgb(626pt)=(0.549777,0.748921,0.457983);
		rgb(627pt)=(0.55235,0.748891,0.45692);
		rgb(628pt)=(0.554903,0.74884,0.455899);
		rgb(629pt)=(0.557456,0.748789,0.454878);
		rgb(630pt)=(0.560008,0.748738,0.453857);
		rgb(631pt)=(0.562554,0.748687,0.452829);
		rgb(632pt)=(0.565081,0.748636,0.451783);
		rgb(633pt)=(0.567608,0.748585,0.450736);
		rgb(634pt)=(0.570135,0.748534,0.449689);
		rgb(635pt)=(0.572653,0.748474,0.44866);
		rgb(636pt)=(0.575155,0.748397,0.447665);
		rgb(637pt)=(0.577656,0.748321,0.446669);
		rgb(638pt)=(0.580158,0.748244,0.445674);
		rgb(639pt)=(0.582649,0.748168,0.444678);
		rgb(640pt)=(0.585125,0.748091,0.443683);
		rgb(641pt)=(0.587601,0.748014,0.442687);
		rgb(642pt)=(0.590077,0.747938,0.441692);
		rgb(643pt)=(0.59254,0.747861,0.440709);
		rgb(644pt)=(0.59499,0.747785,0.439739);
		rgb(645pt)=(0.597441,0.747708,0.438769);
		rgb(646pt)=(0.599891,0.747632,0.437799);
		rgb(647pt)=(0.602311,0.747555,0.436814);
		rgb(648pt)=(0.604711,0.747478,0.435819);
		rgb(649pt)=(0.60711,0.747402,0.434823);
		rgb(650pt)=(0.60951,0.747325,0.433828);
		rgb(651pt)=(0.611909,0.747232,0.432867);
		rgb(652pt)=(0.614308,0.747129,0.431922);
		rgb(653pt)=(0.616708,0.747027,0.430978);
		rgb(654pt)=(0.619107,0.746925,0.430033);
		rgb(655pt)=(0.621487,0.746823,0.429089);
		rgb(656pt)=(0.623861,0.746721,0.428144);
		rgb(657pt)=(0.626235,0.746619,0.4272);
		rgb(658pt)=(0.628609,0.746517,0.426256);
		rgb(659pt)=(0.630962,0.746393,0.425311);
		rgb(660pt)=(0.63331,0.746266,0.424367);
		rgb(661pt)=(0.635658,0.746138,0.423422);
		rgb(662pt)=(0.638007,0.746011,0.422478);
		rgb(663pt)=(0.640332,0.745906,0.421557);
		rgb(664pt)=(0.642654,0.745804,0.420638);
		rgb(665pt)=(0.644977,0.745702,0.419719);
		rgb(666pt)=(0.6473,0.7456,0.4188);
		rgb(667pt)=(0.649623,0.745472,0.417881);
		rgb(668pt)=(0.651946,0.745345,0.416962);
		rgb(669pt)=(0.654268,0.745217,0.416043);
		rgb(670pt)=(0.656587,0.745089,0.415124);
		rgb(671pt)=(0.658859,0.744962,0.414205);
		rgb(672pt)=(0.661131,0.744834,0.413286);
		rgb(673pt)=(0.663402,0.744707,0.412368);
		rgb(674pt)=(0.665674,0.744579,0.411453);
		rgb(675pt)=(0.667946,0.744451,0.410559);
		rgb(676pt)=(0.670218,0.744324,0.409666);
		rgb(677pt)=(0.672489,0.744196,0.408773);
		rgb(678pt)=(0.674755,0.744062,0.407879);
		rgb(679pt)=(0.677001,0.743909,0.406986);
		rgb(680pt)=(0.679247,0.743756,0.406092);
		rgb(681pt)=(0.681494,0.743603,0.405199);
		rgb(682pt)=(0.68374,0.743458,0.404306);
		rgb(683pt)=(0.685986,0.74333,0.403412);
		rgb(684pt)=(0.688232,0.743203,0.402519);
		rgb(685pt)=(0.690479,0.743075,0.401626);
		rgb(686pt)=(0.692704,0.742937,0.400732);
		rgb(687pt)=(0.694899,0.742784,0.399839);
		rgb(688pt)=(0.697094,0.742631,0.398945);
		rgb(689pt)=(0.699289,0.742477,0.398052);
		rgb(690pt)=(0.701497,0.742324,0.397171);
		rgb(691pt)=(0.703718,0.742171,0.396303);
		rgb(692pt)=(0.705939,0.742018,0.395435);
		rgb(693pt)=(0.708159,0.741865,0.394568);
		rgb(694pt)=(0.710351,0.741712,0.3937);
		rgb(695pt)=(0.71252,0.741559,0.392832);
		rgb(696pt)=(0.71469,0.741405,0.391964);
		rgb(697pt)=(0.71686,0.741252,0.391096);
		rgb(698pt)=(0.719029,0.741082,0.390228);
		rgb(699pt)=(0.721199,0.740904,0.38936);
		rgb(700pt)=(0.723369,0.740725,0.388492);
		rgb(701pt)=(0.725538,0.740546,0.387625);
		rgb(702pt)=(0.727708,0.740386,0.386757);
		rgb(703pt)=(0.729878,0.740233,0.385889);
		rgb(704pt)=(0.732047,0.74008,0.385021);
		rgb(705pt)=(0.734217,0.739927,0.384153);
		rgb(706pt)=(0.736366,0.739753,0.383285);
		rgb(707pt)=(0.73851,0.739574,0.382417);
		rgb(708pt)=(0.740654,0.739395,0.38155);
		rgb(709pt)=(0.742798,0.739217,0.380682);
		rgb(710pt)=(0.744919,0.739038,0.379837);
		rgb(711pt)=(0.747038,0.738859,0.378995);
		rgb(712pt)=(0.749156,0.738681,0.378152);
		rgb(713pt)=(0.751275,0.738502,0.37731);
		rgb(714pt)=(0.753394,0.738323,0.376442);
		rgb(715pt)=(0.755512,0.738145,0.375574);
		rgb(716pt)=(0.757631,0.737966,0.374707);
		rgb(717pt)=(0.75975,0.737789,0.373841);
		rgb(718pt)=(0.761868,0.737636,0.372998);
		rgb(719pt)=(0.763987,0.737483,0.372156);
		rgb(720pt)=(0.766105,0.73733,0.371314);
		rgb(721pt)=(0.76822,0.737169,0.370471);
		rgb(722pt)=(0.770313,0.736965,0.369629);
		rgb(723pt)=(0.772406,0.73676,0.368786);
		rgb(724pt)=(0.774499,0.736556,0.367944);
		rgb(725pt)=(0.776592,0.736358,0.367096);
		rgb(726pt)=(0.778686,0.736179,0.366228);
		rgb(727pt)=(0.780779,0.736001,0.36536);
		rgb(728pt)=(0.782872,0.735822,0.364492);
		rgb(729pt)=(0.784957,0.735643,0.363632);
		rgb(730pt)=(0.787024,0.735465,0.36279);
		rgb(731pt)=(0.789092,0.735286,0.361948);
		rgb(732pt)=(0.791159,0.735107,0.361105);
		rgb(733pt)=(0.793227,0.734929,0.360263);
		rgb(734pt)=(0.795295,0.73475,0.359421);
		rgb(735pt)=(0.797362,0.734571,0.358578);
		rgb(736pt)=(0.79943,0.734392,0.357736);
		rgb(737pt)=(0.801485,0.734214,0.356881);
		rgb(738pt)=(0.803527,0.734035,0.356014);
		rgb(739pt)=(0.805569,0.733856,0.355146);
		rgb(740pt)=(0.807611,0.733678,0.354278);
		rgb(741pt)=(0.809668,0.733499,0.353424);
		rgb(742pt)=(0.811735,0.73332,0.352582);
		rgb(743pt)=(0.813803,0.733142,0.35174);
		rgb(744pt)=(0.81587,0.732963,0.350897);
		rgb(745pt)=(0.817921,0.732784,0.350038);
		rgb(746pt)=(0.819963,0.732606,0.349171);
		rgb(747pt)=(0.822005,0.732427,0.348303);
		rgb(748pt)=(0.824047,0.732248,0.347435);
		rgb(749pt)=(0.826071,0.73207,0.346567);
		rgb(750pt)=(0.828087,0.731891,0.345699);
		rgb(751pt)=(0.830104,0.731712,0.344831);
		rgb(752pt)=(0.83212,0.731534,0.343963);
		rgb(753pt)=(0.834158,0.731355,0.343095);
		rgb(754pt)=(0.8362,0.731176,0.342228);
		rgb(755pt)=(0.838242,0.730998,0.34136);
		rgb(756pt)=(0.840284,0.730819,0.340492);
		rgb(757pt)=(0.842303,0.73064,0.339624);
		rgb(758pt)=(0.84432,0.730462,0.338756);
		rgb(759pt)=(0.846336,0.730283,0.337888);
		rgb(760pt)=(0.848353,0.730104,0.33702);
		rgb(761pt)=(0.850369,0.729926,0.336153);
		rgb(762pt)=(0.852386,0.729747,0.335285);
		rgb(763pt)=(0.854402,0.729568,0.334417);
		rgb(764pt)=(0.856419,0.729391,0.333546);
		rgb(765pt)=(0.858435,0.729238,0.332627);
		rgb(766pt)=(0.860452,0.729085,0.331708);
		rgb(767pt)=(0.862468,0.728932,0.330789);
		rgb(768pt)=(0.864481,0.728778,0.329874);
		rgb(769pt)=(0.866472,0.728625,0.32898);
		rgb(770pt)=(0.868463,0.728472,0.328087);
		rgb(771pt)=(0.870454,0.728319,0.327194);
		rgb(772pt)=(0.872445,0.728166,0.326295);
		rgb(773pt)=(0.874436,0.728013,0.325376);
		rgb(774pt)=(0.876427,0.727859,0.324457);
		rgb(775pt)=(0.878418,0.727706,0.323538);
		rgb(776pt)=(0.880417,0.727561,0.322619);
		rgb(777pt)=(0.882433,0.727433,0.3217);
		rgb(778pt)=(0.88445,0.727306,0.320781);
		rgb(779pt)=(0.886466,0.727178,0.319862);
		rgb(780pt)=(0.888463,0.72705,0.318933);
		rgb(781pt)=(0.890429,0.726923,0.317989);
		rgb(782pt)=(0.892394,0.726795,0.317044);
		rgb(783pt)=(0.894359,0.726668,0.3161);
		rgb(784pt)=(0.896337,0.726552,0.315132);
		rgb(785pt)=(0.898328,0.72645,0.314136);
		rgb(786pt)=(0.900319,0.726348,0.313141);
		rgb(787pt)=(0.90231,0.726246,0.312145);
		rgb(788pt)=(0.904301,0.726158,0.31115);
		rgb(789pt)=(0.906292,0.726081,0.310154);
		rgb(790pt)=(0.908283,0.726005,0.309159);
		rgb(791pt)=(0.910274,0.725928,0.308163);
		rgb(792pt)=(0.912249,0.725851,0.307151);
		rgb(793pt)=(0.914214,0.725775,0.30613);
		rgb(794pt)=(0.91618,0.725698,0.305109);
		rgb(795pt)=(0.918145,0.725622,0.304088);
		rgb(796pt)=(0.920111,0.7256,0.303031);
		rgb(797pt)=(0.922076,0.7256,0.301959);
		rgb(798pt)=(0.924041,0.7256,0.300886);
		rgb(799pt)=(0.926007,0.7256,0.299814);
		rgb(800pt)=(0.927972,0.7256,0.298722);
		rgb(801pt)=(0.929938,0.7256,0.297624);
		rgb(802pt)=(0.931903,0.7256,0.296527);
		rgb(803pt)=(0.933869,0.7256,0.295429);
		rgb(804pt)=(0.935812,0.725668,0.294264);
		rgb(805pt)=(0.937752,0.725744,0.29309);
		rgb(806pt)=(0.939692,0.725821,0.291916);
		rgb(807pt)=(0.941632,0.725897,0.290741);
		rgb(808pt)=(0.943571,0.726023,0.289518);
		rgb(809pt)=(0.945511,0.726151,0.288293);
		rgb(810pt)=(0.947451,0.726278,0.287068);
		rgb(811pt)=(0.949389,0.726411,0.285839);
		rgb(812pt)=(0.951278,0.726641,0.284537);
		rgb(813pt)=(0.953167,0.72687,0.283235);
		rgb(814pt)=(0.955056,0.7271,0.281933);
		rgb(815pt)=(0.956938,0.72734,0.280622);
		rgb(816pt)=(0.958776,0.727646,0.279243);
		rgb(817pt)=(0.960614,0.727952,0.277865);
		rgb(818pt)=(0.962451,0.728259,0.276486);
		rgb(819pt)=(0.964273,0.728597,0.275086);
		rgb(820pt)=(0.966034,0.729057,0.273606);
		rgb(821pt)=(0.967795,0.729516,0.272126);
		rgb(822pt)=(0.969557,0.729976,0.270645);
		rgb(823pt)=(0.971288,0.730473,0.269135);
		rgb(824pt)=(0.972947,0.73106,0.267552);
		rgb(825pt)=(0.974606,0.731647,0.265969);
		rgb(826pt)=(0.976265,0.732234,0.264387);
		rgb(827pt)=(0.977857,0.732879,0.262785);
		rgb(828pt)=(0.979338,0.733619,0.261151);
		rgb(829pt)=(0.980818,0.734359,0.259518);
		rgb(830pt)=(0.982299,0.735099,0.257884);
		rgb(831pt)=(0.983697,0.73591,0.256227);
		rgb(832pt)=(0.984999,0.736803,0.254542);
		rgb(833pt)=(0.986301,0.737697,0.252858);
		rgb(834pt)=(0.987603,0.73859,0.251173);
		rgb(835pt)=(0.988753,0.739566,0.249474);
		rgb(836pt)=(0.989774,0.740613,0.247764);
		rgb(837pt)=(0.990795,0.741659,0.246054);
		rgb(838pt)=(0.991816,0.742706,0.244344);
		rgb(839pt)=(0.992677,0.743816,0.242681);
		rgb(840pt)=(0.993443,0.744965,0.241048);
		rgb(841pt)=(0.994209,0.746114,0.239414);
		rgb(842pt)=(0.994975,0.747262,0.23778);
		rgb(843pt)=(0.995578,0.748465,0.236165);
		rgb(844pt)=(0.996114,0.74969,0.234557);
		rgb(845pt)=(0.99665,0.750915,0.232949);
		rgb(846pt)=(0.997186,0.752141,0.231341);
		rgb(847pt)=(0.997562,0.753386,0.229813);
		rgb(848pt)=(0.997893,0.754637,0.228307);
		rgb(849pt)=(0.998225,0.755887,0.226801);
		rgb(850pt)=(0.998557,0.757138,0.225295);
		rgb(851pt)=(0.998711,0.758433,0.223856);
		rgb(852pt)=(0.998839,0.759735,0.222426);
		rgb(853pt)=(0.998966,0.761037,0.220997);
		rgb(854pt)=(0.999094,0.762339,0.219567);
		rgb(855pt)=(0.999076,0.763641,0.218186);
		rgb(856pt)=(0.99905,0.764942,0.216808);
		rgb(857pt)=(0.999025,0.766244,0.21543);
		rgb(858pt)=(0.998995,0.767546,0.214054);
		rgb(859pt)=(0.998868,0.768848,0.212752);
		rgb(860pt)=(0.99874,0.77015,0.21145);
		rgb(861pt)=(0.998613,0.771451,0.210149);
		rgb(862pt)=(0.998473,0.772756,0.208856);
		rgb(863pt)=(0.998243,0.774083,0.207631);
		rgb(864pt)=(0.998014,0.775411,0.206405);
		rgb(865pt)=(0.997784,0.776738,0.20518);
		rgb(866pt)=(0.997539,0.77806,0.20396);
		rgb(867pt)=(0.997232,0.779362,0.20276);
		rgb(868pt)=(0.996926,0.780664,0.201561);
		rgb(869pt)=(0.99662,0.781966,0.200361);
		rgb(870pt)=(0.996299,0.783268,0.199168);
		rgb(871pt)=(0.995942,0.784569,0.197994);
		rgb(872pt)=(0.995584,0.785871,0.19682);
		rgb(873pt)=(0.995227,0.787173,0.195646);
		rgb(874pt)=(0.994842,0.788475,0.19449);
		rgb(875pt)=(0.994408,0.789777,0.193367);
		rgb(876pt)=(0.993974,0.791078,0.192244);
		rgb(877pt)=(0.99354,0.79238,0.191121);
		rgb(878pt)=(0.993083,0.793671,0.190021);
		rgb(879pt)=(0.992598,0.794947,0.188949);
		rgb(880pt)=(0.992113,0.796223,0.187877);
		rgb(881pt)=(0.991628,0.797499,0.186805);
		rgb(882pt)=(0.99113,0.798789,0.185732);
		rgb(883pt)=(0.990619,0.800091,0.18466);
		rgb(884pt)=(0.990109,0.801393,0.183588);
		rgb(885pt)=(0.989598,0.802695,0.182516);
		rgb(886pt)=(0.989072,0.803996,0.18146);
		rgb(887pt)=(0.988536,0.805298,0.180413);
		rgb(888pt)=(0.988,0.8066,0.179367);
		rgb(889pt)=(0.987464,0.807902,0.17832);
		rgb(890pt)=(0.98691,0.809186,0.177291);
		rgb(891pt)=(0.986349,0.810462,0.17627);
		rgb(892pt)=(0.985787,0.811738,0.175249);
		rgb(893pt)=(0.985226,0.813015,0.174228);
		rgb(894pt)=(0.984644,0.814311,0.173207);
		rgb(895pt)=(0.984057,0.815613,0.172186);
		rgb(896pt)=(0.98347,0.816914,0.171165);
		rgb(897pt)=(0.982883,0.818216,0.170144);
		rgb(898pt)=(0.982296,0.819518,0.169145);
		rgb(899pt)=(0.981709,0.82082,0.16815);
		rgb(900pt)=(0.981122,0.822122,0.167154);
		rgb(901pt)=(0.980535,0.823423,0.166159);
		rgb(902pt)=(0.979947,0.824725,0.165163);
		rgb(903pt)=(0.97936,0.826027,0.164168);
		rgb(904pt)=(0.978773,0.827329,0.163172);
		rgb(905pt)=(0.978186,0.828631,0.162177);
		rgb(906pt)=(0.977599,0.829932,0.161181);
		rgb(907pt)=(0.977012,0.831234,0.160186);
		rgb(908pt)=(0.976425,0.832536,0.15919);
		rgb(909pt)=(0.975838,0.833841,0.158195);
		rgb(910pt)=(0.975251,0.835168,0.157199);
		rgb(911pt)=(0.974664,0.836495,0.156204);
		rgb(912pt)=(0.974077,0.837823,0.155208);
		rgb(913pt)=(0.973489,0.83915,0.154213);
		rgb(914pt)=(0.972902,0.840477,0.153217);
		rgb(915pt)=(0.972315,0.841805,0.152222);
		rgb(916pt)=(0.971728,0.843132,0.151226);
		rgb(917pt)=(0.971155,0.844466,0.150224);
		rgb(918pt)=(0.970619,0.845819,0.149203);
		rgb(919pt)=(0.970083,0.847172,0.148182);
		rgb(920pt)=(0.969547,0.848525,0.147161);
		rgb(921pt)=(0.96902,0.849886,0.14614);
		rgb(922pt)=(0.968509,0.851265,0.145119);
		rgb(923pt)=(0.967999,0.852643,0.144098);
		rgb(924pt)=(0.967488,0.854022,0.143077);
		rgb(925pt)=(0.967,0.855411,0.142056);
		rgb(926pt)=(0.966541,0.856815,0.141035);
		rgb(927pt)=(0.966081,0.858219,0.140014);
		rgb(928pt)=(0.965622,0.859623,0.138992);
		rgb(929pt)=(0.965189,0.86104,0.137945);
		rgb(930pt)=(0.96478,0.862469,0.136873);
		rgb(931pt)=(0.964372,0.863899,0.135801);
		rgb(932pt)=(0.963963,0.865328,0.134729);
		rgb(933pt)=(0.96357,0.866773,0.133657);
		rgb(934pt)=(0.963187,0.868228,0.132585);
		rgb(935pt)=(0.962805,0.869683,0.131513);
		rgb(936pt)=(0.962422,0.871138,0.130441);
		rgb(937pt)=(0.962091,0.87261,0.129351);
		rgb(938pt)=(0.961785,0.874091,0.128253);
		rgb(939pt)=(0.961478,0.875571,0.127156);
		rgb(940pt)=(0.961172,0.877052,0.126058);
		rgb(941pt)=(0.960885,0.878571,0.124961);
		rgb(942pt)=(0.960605,0.880103,0.123863);
		rgb(943pt)=(0.960324,0.881634,0.122765);
		rgb(944pt)=(0.960043,0.883166,0.121668);
		rgb(945pt)=(0.959849,0.884719,0.120549);
		rgb(946pt)=(0.95967,0.886276,0.119426);
		rgb(947pt)=(0.959491,0.887833,0.118302);
		rgb(948pt)=(0.959313,0.88939,0.117179);
		rgb(949pt)=(0.959181,0.890995,0.116032);
		rgb(950pt)=(0.959054,0.892603,0.114884);
		rgb(951pt)=(0.958926,0.894211,0.113735);
		rgb(952pt)=(0.958799,0.895819,0.112587);
		rgb(953pt)=(0.958748,0.897453,0.111464);
		rgb(954pt)=(0.958697,0.899086,0.110341);
		rgb(955pt)=(0.958646,0.90072,0.109217);
		rgb(956pt)=(0.958602,0.902359,0.108089);
		rgb(957pt)=(0.958628,0.904043,0.106915);
		rgb(958pt)=(0.958653,0.905728,0.105741);
		rgb(959pt)=(0.958679,0.907413,0.104567);
		rgb(960pt)=(0.958718,0.909102,0.103393);
		rgb(961pt)=(0.95882,0.910812,0.102219);
		rgb(962pt)=(0.958922,0.912522,0.101044);
		rgb(963pt)=(0.959024,0.914232,0.0998703);
		rgb(964pt)=(0.959153,0.915962,0.0986961);
		rgb(965pt)=(0.959357,0.917749,0.0975219);
		rgb(966pt)=(0.959561,0.919536,0.0963477);
		rgb(967pt)=(0.959765,0.921323,0.0951736);
		rgb(968pt)=(0.959996,0.923118,0.0939907);
		rgb(969pt)=(0.960277,0.924931,0.092791);
		rgb(970pt)=(0.960557,0.926743,0.0915913);
		rgb(971pt)=(0.960838,0.928555,0.0903916);
		rgb(972pt)=(0.961151,0.930378,0.0891919);
		rgb(973pt)=(0.961509,0.932216,0.0879922);
		rgb(974pt)=(0.961866,0.934054,0.0867925);
		rgb(975pt)=(0.962223,0.935892,0.0855928);
		rgb(976pt)=(0.96262,0.937768,0.0843802);
		rgb(977pt)=(0.963053,0.939683,0.083155);
		rgb(978pt)=(0.963487,0.941597,0.0819297);
		rgb(979pt)=(0.963921,0.943512,0.0807045);
		rgb(980pt)=(0.964415,0.945426,0.0794643);
		rgb(981pt)=(0.964951,0.947341,0.0782135);
		rgb(982pt)=(0.965487,0.949255,0.0769628);
		rgb(983pt)=(0.966023,0.951169,0.075712);
		rgb(984pt)=(0.966594,0.953118,0.0744441);
		rgb(985pt)=(0.967181,0.955083,0.0731679);
		rgb(986pt)=(0.967768,0.957049,0.0718916);
		rgb(987pt)=(0.968355,0.959014,0.0706153);
		rgb(988pt)=(0.96898,0.960999,0.0693006);
		rgb(989pt)=(0.969619,0.96299,0.0679733);
		rgb(990pt)=(0.970257,0.964981,0.0666459);
		rgb(991pt)=(0.970895,0.966972,0.0653186);
		rgb(992pt)=(0.971554,0.968984,0.0639486);
		rgb(993pt)=(0.972218,0.971001,0.0625703);
		rgb(994pt)=(0.972882,0.973017,0.0611919);
		rgb(995pt)=(0.973545,0.975034,0.0598135);
		rgb(996pt)=(0.974232,0.97705,0.058318);
		rgb(997pt)=(0.974922,0.979067,0.056812);
		rgb(998pt)=(0.975611,0.981083,0.055306);
		rgb(999pt)=(0.9763,0.9831,0.0538)
	}
}

\usetikzlibrary{arrows,shapes,positioning}
\usetikzlibrary{decorations.markings}
\tikzstyle arrowstyle=[scale=1]
\tikzstyle directed=[postaction={decorate,decoration={markings,
		mark=at position .65 with {\arrow[arrowstyle]{stealth}}}}]
\tikzstyle reverse directed=[postaction={decorate,decoration={markings,
		mark=at position .65 with {\arrowreversed[arrowstyle]{stealth};}}}]

\newcommand{\secref}[1]{Section~\ref{#1}}

\newcommand{\rmref}[1]{Remark~\ref{#1}}

\newcommand{\figref}[1]{Figure~\ref{#1}}
\newcommand{\tabref}[1]{Table~\ref{#1}}
\newcommand{\algoref}[1]{Algorithm~\ref{#1}}

\usepackage{mathtools}
\DeclarePairedDelimiter{\ceil}{\lceil}{\rceil}

\newcommand{\R}{\mathbb{R}}

\newcommand{\Domain}{\ensuremath{X}} 
\newcommand{\dimension}{\ensuremath{d}} 
\newcommand{\timevar}{\ensuremath{t}} 

\newcommand{\Lp}[1]{\ensuremath{L_{#1}}}

\newcommand{\flux}{\ensuremath{\mathbf{f}}}

\newcommand{\x}{\ensuremath{x}}

\newcommand{\dx}{\partial_{\x}}

\newcommand{\quadxpoint}[1]{\hat{\x}_{#1}}
\newcommand{\quadxweight}[1]{\hat{w}_{#1}}
\newcommand{\quadRpoint}[1]{\hat{\uncertainty}_{#1}}
\newcommand{\quadRweight}[1]{\hat{\omega}_{#1}}

\newcommand{\dtstepsize}{\ensuremath{\Delta \timevar}}

\newcommand{\ncells}{\ensuremath{N_x}}

\newcommand{\timeind}{\ensuremath{n}}

\newcommand{\cellind}{\ensuremath{i}}
\newcommand{\cellindR}{\ensuremath{j}}

\newcommand{\cell}[1]{\ensuremath{{\Domain}_{#1}}}
\newcommand{\cellR}[1]{\ensuremath{\sampleSpace_{#1}}}

\newcommand{\interface}[1]{\ensuremath{\x_{#1}}}

\newcommand{\limitervariable}{\ensuremath{\theta}}

\newcommand{\Testfunction}[1][ ]{\ensuremath{v_{#1}}}

\newcommand{\numericalFlux}{\ensuremath{\widehat{\flux}}}
\newcommand{\viscosityconstant}{\ensuremath{c}}

\newcommand{\xiv}{\ensuremath{\tilde{\xi}}}

\newcommand{\density}{\ensuremath{\rho}}

\newcommand{\energy}{\ensuremath{E}}
\newcommand{\pressure}{\ensuremath{p}}
\newcommand{\eulergamma}{\ensuremath{\gamma}}

\newcommand{\momentum}{\ensuremath{m}}

\newcommand{\uncertainty}{\ensuremath{\xi}}
\newcommand{\xiPDF}{\ensuremath{f_\Xi}}
\newcommand{\xiBasisPoly}[1]{\ensuremath{\phi_{#1}}}
\newcommand{\SGsumIndex}{\ensuremath{k}}

\newcommand{\SGeqIndex}{\ensuremath{l}}
\newcommand{\SGapproach}{\ensuremath{\sum_{\SGsumIndex=0}^\SGtruncorder \solution_\SGsumIndex \xiBasisPoly{\SGsumIndex}}}
\newcommand{\xiPDFdxi}{\ensuremath{\xiPDF \mathrm{d}\uncertainty}}
\newcommand{\MExiPDFdxi}{\ensuremath{\xiPDF \mathrm{d}\uncertainty}}
\newcommand{\intRS}{\ensuremath{\int_{\randomSpace_\MEIndex}}}
\newcommand{\intcell}[2]{\ensuremath{\int_{\cell{#1}} #2\,\mathrm{d}\x}}

\newcommand{\MEintcellR}[2]{\ensuremath{\int_{\cellR{#1}} #2\,\MExiPDFdxi}}

\newcommand{\MEintxR}[1]{\intcell{\cellind}{\!\MEintcellR{\cellindR}{#1}}}
\newcommand{\SGtruncorder}{\ensuremath{{K_\randomSpace}}}

\newcommand{\weno}{\mathcal{W}}

\newcommand{\nbxnodes}{\ensuremath{Q_\Domain}}
\newcommand{\nbRnodes}{\ensuremath{Q_\randomSpace}}
\newcommand{\nbxiQuadNodes}{\ensuremath{Q_\Gamma}}

\newcommand{\xiQuadIndex}{\ensuremath{\rho}}
\newcommand{\xQuadIndex}{\ensuremath{q}}

\newcommand{\solutioncm}[2]{{\overline{\solution}}_{#2}^{(#1)}}

\newcommand{\sampleSpace}{{\ensuremath{\Xi}}}
\newcommand{\sigmaAlgebra}{\ensuremath{\mathcal{F}}}
\newcommand{\probabilityMeasure}{\ensuremath{\mathcal{P}}}
\newcommand{\randomSpace}{{\ensuremath{\Xi}}}

\newcommand{\solution}{\ensuremath{\bu}}
\newcommand{\solOned}{\ensuremath{u}}

\newcommand{\rhs}{\ensuremath{L_h}}

\newcommand{\stageind}{s}

\newcommand{\randomElement}[1]{\ensuremath{\sampleSpace_{#1}}}
\newcommand{\indicatorVar}[2]{\ensuremath{\chi_{#1}(#2)}}

\newcommand{\MEIndex}{\ensuremath{j}}
\newcommand{\MEElements}{\ensuremath{{N_\randomSpace}}}
\newcommand{\localxiBasisPoly}[2]{\ensuremath{\phi_{#1,#2}}}

\newcommand{\landau}{\mathcal{O}}

\newcommand{\RKorder}{\ensuremath{R}}
\newcommand{\WENOdegree}{\ensuremath{K_\Domain}}
\newcommand{\WENOdegreeCoeff}{\ensuremath{r}}
\newcommand{\solutioncmRK}[2]{{\overline{\boldsymbol{v}}}^{(#1)}}
\newcommand{\solutionRK}[2]{{\boldsymbol{v}}^{(#1)}}

\newcommand{\stochsllimiter}[1]{\ensuremath{\Lambda\Pi_\uncertainty\!\left(#1\right)}}

\newcommand{\scheme}{Weighted Essentially Non-Oscillatory stochastic Galerkin}
\newcommand{\schemeshort}{WENOsG}
\newcommand{\ME}{Multielement}
\newcommand{\fullWENO}{2D WENO}

\newcommand{\SD}{\sampleSpace}
\newcommand{\PD}{X}
\newcommand{\RR}{\mathbb{R}}
\newcommand{\pd}[1]{\frac{\partial}{\partial #1}}
\newcommand{\dint}[1]{{\,\text{d}#1}}

\pgfplotscreateplotcyclelist{color parula}{%
	royalblue!70,every mark/.append style={solid,line width = 0pt,fill=royalblue!60!black},mark=ball\\%
	goldenrod!50!yelloworange,every mark/.append style={solid,fill=goldenrod!30!black},mark=*\\%
	green,every mark/.append style={black,fill=yellowgreen},mark=diamond*\\%
	bluegreen,every mark/.append style={fill=black},mark=triangle*\\%
	royalblue!70,densely dashed,every mark/.append style={solid,line width = 0pt,fill=royalblue!60!black},mark=ball\\%
	goldenrod!50!yelloworange,densely dashed,every mark/.append style={solid,black,fill=goldenrod},mark=diamond*\\%
	yellowgreen,densely dashed,every mark/.append style={solid,fill=yellowgreen!60!royalblue},mark=square*, mark size= 1pt\\%
	bluegreen,densely dashed,mark size = 1.5pt,every mark/.append style={solid,fill=white},mark=otimes*\\
	royalblue,densely dashed,mark size = 1.5pt,every mark/.append style={solid,fill=royalblue!30!black},mark=pentagon*\\
	goldenrod!40!yelloworange,every mark/.append style={solid,fill=goldenrod!30!black},mark=|\\%
}

\pgfplotscreateplotcyclelist{color louisa}{%
	blue,every mark/.append style={solid,line width = 0pt,fill=royalblue!60!black},mark=ball\\%
	red,every mark/.append style={solid,fill=orangered!30!black},mark=*\\%
	green,every mark/.append style={black,fill=yellowgreen},mark=diamond*\\%
	goldenrod!50!yelloworange,every mark/.append style={fill=black},mark=triangle*\\%
	royalblue!70,densely dashed,every mark/.append style={solid,line width = 0pt,fill=royalblue!60!black},mark=ball\\%
	goldenrod!50!yelloworange,densely dashed,every mark/.append style={solid,black,fill=goldenrod},mark=diamond*\\%
	yellowgreen,densely dashed,every mark/.append style={solid,fill=yellowgreen!60!royalblue},mark=square*, mark size= 1pt\\%
	bluegreen,densely dashed,mark size = 1.5pt,every mark/.append style={solid,fill=white},mark=otimes*\\
	royalblue,densely dashed,mark size = 1.5pt,every mark/.append style={solid,fill=royalblue!30!black},mark=pentagon*\\
	goldenrod!40!yelloworange,every mark/.append style={solid,fill=goldenrod!30!black},mark=|\\%
}

\pgfplotscreateplotcyclelist{color fabian}{%
	black,dashed\\%
	red,every mark/.append style={solid},mark=diamond*\\%
	blue,every mark/.append style={solid},mark=*\\%
	goldenrod!50!yelloworange,every mark/.append style={fill=black},mark=triangle*\\%
	royalblue!70,densely dashed,every mark/.append style={solid,line width = 0pt,fill=royalblue!60!black},mark=ball\\%
	goldenrod!50!yelloworange,densely dashed,every mark/.append style={solid,black,fill=goldenrod},mark=diamond*\\%
	yellowgreen,densely dashed,every mark/.append style={solid,fill=yellowgreen!60!royalblue},mark=square*, mark size= 1pt\\%
	bluegreen,densely dashed,mark size = 1.5pt,every mark/.append style={solid,fill=white},mark=otimes*\\
	royalblue,densely dashed,mark size = 1.5pt,every mark/.append style={solid,fill=royalblue!30!black},mark=pentagon*\\
	goldenrod!40!yelloworange,every mark/.append style={solid,fill=goldenrod!30!black},mark=|\\%
}

\usepackage[colorinlistoftodos,obeyFinal]{todonotes}

\usepackage{media9}

\begin{document}

\begin{abstract}
In this paper we extensively study the stochastic Galerkin scheme for uncertain systems of conservation laws, which appears to produce oscillations already for a simple example of the linear advection equation with Riemann initial data. Therefore, we introduce a modified scheme that we call the weighted essentially non-oscillatory (WENO) stochastic Galerkin scheme, which is constructed to prevent the propagation of Gibbs phenomenon into the stochastic domain by applying a slope limiter in the stochasticity. In order to achieve a high order method, we use a spatial WENO reconstruction and also compare the results to a scheme that uses WENO reconstruction in both the physical and the stochastic domain. We evaluate these methods by presenting various numerical test cases where we observe the reduction of the total variation compared to classical stochastic Galerkin.
\end{abstract}
\maketitle

\noindent


\section{Introduction}
Many physical problem settings can be described by hyperbolic systems of conservation laws, however,  crucial data or parameters might not be available exactly due to measurement errors and thus have non-deterministic effects on the approximation of the deterministic systems. The modeling of the propagation of this uncertainty into the solution is the topic of Uncertainty Quantification (UQ).
In this context, UQ methods \cite{abgrall2017uncertainty,Kolb2018,le2010spectral,pettersson2015polynomial,Smith2014,Kusch2018,Kusch2018a} are gaining more and more popularity, whereas we distinguish between two approaches, namely the so-called non-intrusive and intrusive schemes. \reffour{In addition to that, recent developments also introduce so-called semi-intrusive \cite{Abgrall2013} or weakly intrusive \cite{Petrella2019,Tokareva2014} methods, which only require small modifications in the deterministic solvers.}

The most widely known non-intrusive UQ method is (Multi-Level) Monte Carlo \cite{lye2016multilevel,Heinrich2001,Giles2008} which is based on statistical sampling methods and can be easily implemented and adopted to any type of uncertain conservation law, but comes with potentially high cost due to repeated application of e.g.\ finite volume methods (FVM). These so-called MC-FVM schemes have been studied in \cite{Mishra2016,mishra2013} showing a slow convergence rate that is improved by Multi-Level MC-FVM algorithms for conservation laws \cite{Mishra2014,mishra2012multi}. Another non-intrusive UQ scheme is stochastic collocation \cite{Xiu2005}, further developments of this method are described in \cite{Witteveen2009,Ma2009,Nobile2008}.

Intrusive UQ methods aim to increase the overall efficiency but require a modification of the underlying solver for the deterministic problem (as FVM). The most popular methods are the Intrusive Polynomial Moment (IPM) method \cite{Poette2009,Kusch2017} and the stochastic Galerkin (sG) scheme. They both rely on the generalized Polynomial Chaos (gPC) expansion \cite{Xiu2003,Wan2006a,Abgrall2007,Chen2005} which is theoretically based on the Polynomial Chaos expansion from \cite{Wiener1938}. IPM expands the stochastic solution in the so-called entropic variables, which results in a hyperbolic gPC system that yields a good approximation quality, however, it is necessary to know a strictly convex entropy solution beforehand. Stochastic Galerkin expands the solution in the conserved variables such that the gPC system results in a weak formulation of the equations with respect to the stochastic variable.

For many non-hyperbolic equations, the underlying random field is sufficiently smooth in the stochasticity such that the sG method is superior to Monte Carlo type methods since the gPC approach exhibits spectral convergence \cite{Babuska2004,ghanem2003stochastic,Yan2005}.
The biggest challenge of UQ methods for hyperbolic conservation laws lies in the fact that discontinuities in the physical space propagate into the solution manifold \cite{barth2013non}.
The naive usage of sG for nonlinear hyperbolic problems even typically fails \cite{abgrall2017uncertainty,Poette2009} since the polynomial expansion of discontinuous data yields huge oscillations that results in the loss of hyperbolicity. In order to resolve this problem, we apply the hyperbolicity limiter from \cite{Schlachter2017a,Meyer2019a} to the classical sG approach. 
In addition to that, the authors of \cite{Wan2006} introduced the so-called Multi-Element approach, where  the random space is divided into disjoint elements in order to define local gPC approximations. Further developments of this method can be found in \cite{Tryoen2010,Wan2005,Wan2009}.

Within this article, we show a simple example of the linear transport equation with uncertain wave speed, supplemented with Riemann initial data, that still produces oscillations in the gPC expansion \refone{\cite{Poette2009,Barth2012,Ma2009}}.
For this reason, we propose a robust numerical method that is able to deal with this kind of Gibbs oscillations. We combine the hyperbolicity limiter and Multi-Element ansatz with a weighted essentially non-oscillatory (WENO) reconstruction \cite{Osher1994,Toro2009,quarteroni2006advanced,JIANG1996,Cravero2017} in the physical space to deduce a high-order method and apply a slope limiter in the stochastic variable. Furthermore, we consider a full two-dimensional WENO reconstruction in both the physical and stochastic domain, motivated by the stochastic finite volume method from \cite{Tokareva2014}.  We compare the performance of both methods by considering the total variation for various numerical test cases. 

Further intrusive UQ methods that aim to damp oscillations induced by the Gibbs phenomenon are given for example in \cite{Kusch2018}, where filters are applied to the gPC coefficients of the stochastic Galerkin approximation. Filters are a common technique from kinetic theory \cite{McClarren2010}, allowing a numerically cheap reduction of oscillations. \refthree{Moreover, if the entropy of the underlying system is known, the IPM method \cite{Poette2009} may be used to control oscillations since they are bounded to a certain range though the entropy. 
Another development of this method can be found in \cite{Kusch2017}, proposing a second-order IPM scheme which fulfills the maximum principle. In this article, we want to explore additional strategies that aim on reducing oscillations by using slope limiting and WENO techniques in sG.}

The paper is structured as follows. In \secref{sec:model} we describe our problem setting, that is the system of uncertain conservation laws which we discretize in the stochastic domain by the stochastic Galerkin scheme. We then demonstrate the propagation of Gibbs phenomenon by an introductory example which yields to the definition of the stochastic slope limiter. \secref{sec:dg} formulates the weighted essentially non-oscillatory stochastic Galerkin scheme such as a full 2D WENO reconstruction of the conservation law. Finally, we show some numerical results in \secref{sec:results}, demonstrating the reduction of the total variation for our methods compared to classical stochastic Galerkin.

\section{Modeling Uncertainties}\label{sec:model}

We consider stochastic conservation laws of the form
\begin{subequations}\label{eq:spde}
\begin{equation}
 \pd{\timevar} \solution(\timevar,\x,\uncertainty) + \pd{\x} \flux\bigl(\solution(\timevar,\x,\uncertainty),\xi\bigr) = 0,\quad \quad \text{for} \ \ \x\in\PD,\, t>0,\, \uncertainty\in\SD,
\end{equation}
with physical domain $\PD \subset \RR$, stochastic domain $\SD \subset \RR$ and initial conditions given by
\begin{equation}
 \solution(0,\x,\uncertainty) = \solution^{(0)}(\x,\uncertainty), \quad \quad \text{for} \ \ \x\in\PD,\, \xi\in\SD.
\end{equation}
\end{subequations}
Depending on $\PD$, additional boundary conditions have to \refone{be} prescribed.
The solution $\solution\in\RR^\dimension$ is depending on a one-dimensional random variable $\uncertainty$ with probability space $(\sampleSpace,\sigmaAlgebra,\probabilityMeasure)$ and probability density function $\xiPDF(\uncertainty) : \randomSpace \rightarrow \R_+$.
We abuse notation and denote the random space of this uncertainty by $\randomSpace := \uncertainty(\sampleSpace)$ and write $\uncertainty$ also for the realizations of the random variable. 

\refone{
\begin{definition}\label{def:hyp}
	The system \eqref{eq:spde} is called hyperbolic, if for any $\solution\in  \R^{\dimension}$ the matrix
	$\frac{\partial\flux(\solution)}{\partial\solution}$
	has $\dimension$ real eigenvalues and is diagonalizable.
\end{definition}
\begin{remark}\label{rem:existunique}
	The existence and uniqueness of a so-called random entropy solution to  \eqref{eq:spde}  is proved in \cite{Mishra2014}. They correspond $\probabilityMeasure$ almost everywhere to entropy solutions in the deterministic case which is discussed in Kru\v{z}kov's theorem \cite{Kruzov1970}. Here, we assume that \eqref{eq:spde} has \refthree{an entropy} - entropy flux pair that satisfies the entropy inequality. For more information on this topic, we refer to \cite{Mishra2016,Meyer2019}. 
\end{remark}
}

\subsection{Stochastic Galerkin}
We seek for an approximate solution by a finite-term generalized Polynomial Chaos (gPC) expansion (see e.g. \cite{Gottlieb2008})
\begin{equation}\label{gPCApprox}
\solution(\timevar,\x,\uncertainty) \approx  \sum\limits_{\SGsumIndex=0}^{\SGtruncorder} \solution_\SGsumIndex(\timevar,\x) \, \xiBasisPoly{\SGsumIndex}(\uncertainty) ,
\end{equation}
where the polynomials $\xiBasisPoly{\SGsumIndex}$ of degree $\SGsumIndex$ are supposed to satisfy the orthogonality relation
\begin{equation}\label{eq:orthogonality}
\int_{\randomSpace} \xiBasisPoly{\SGsumIndex}(\uncertainty) \, \xiBasisPoly{\SGeqIndex}(\uncertainty) \, \xiPDF(\uncertainty) \dint{\uncertainty} = \delta_{\SGsumIndex\SGeqIndex},
\end{equation}
for $\SGsumIndex,\SGeqIndex = 0, \dots, \SGtruncorder$.
Inserting \eqref{gPCApprox} into \eqref{eq:spde} and applying a Galerkin projection in the stochastic space leads to the so called stochastic Galerkin system
\begin{equation}\label{SGsystem} \pd{t} \solution_\SGeqIndex +  \pd{\x} \intRS \flux\!\left(\SGapproach\right) \!\xiBasisPoly{\SGeqIndex}\xiPDFdxi = 0, 
\end{equation}
\refone{with $\SGeqIndex = 0,\ldots,\SGtruncorder.$ \refthree{For example,} the expected value and variance of $\solution$ are given by
	\begin{align}\label{eq:exp}
\mathbb{E}(\solution) &\approx  \int_{\randomSpace}  \SGapproach \, \xiPDFdxi 
= \sum_{\SGsumIndex=0}^\SGtruncorder \solution_\SGsumIndex  \int_{\randomSpace} \!\xiBasisPoly{\SGsumIndex}\xiBasisPoly{0}\, \xiPDFdxi
= \solution_0,\\[0.2cm]
\mathbb{V}(\solution) &\approx  \int_{\randomSpace} \! \left(  \SGapproach\right)^2  \! \!\xiPDFdxi - \solution_0^2
= \sum_{\SGsumIndex=1}^\SGtruncorder \solution_\SGsumIndex^2, \label{eq:var}
\end{align}
since \eqref{eq:orthogonality} yields $\xiBasisPoly{0}=1$. \refthree{Higher statistical moments and correlation functions can be obtained analogously. For a detailed uncertainty quantification study of the underlying random conservation law, the gPC expansion \eqref{gPCApprox} may additionally be used for a sensitivity analysis as discussed in \cite{Sudret2008,Smith2014}.}
}

\subsection{\ME{} Ansatz}\label{sec:ME}
For discontinuous solutions, the gPC approach may converge slowly or even fail to converge, cf. \cite{Wan2005, Poette2009}.  As presented in \cite{Wan2006,Meyer2019a}, we therefore apply the \ME{} approach, where $\randomSpace$ is divided into disjoint elements with local gPC approximations of  \eqref{eq:spde}.

We assume that $\randomSpace=(\uncertainty_L,\uncertainty_R)$ and define a decomposition of $\randomSpace$ into $\MEElements$ \ME{}s $\randomElement{\MEIndex} = (\xi_\MEIndex, \xi_{\MEIndex+1})$ of width $\Delta\uncertainty=\frac{\uncertainty_L-\uncertainty_R}{\MEElements}$ for every \refone{$j=1,\ldots,\MEElements$}. 
\refone{\begin{remark}
	If the random variable $\uncertainty$ is defined on an unbounded domain (with unbounded support), for example in the case of a normal distribution, we refer to the strategy proposed in \cite{Wan2006}. Here, the idea is to subdivide $\R$ into three elements $(-\infty,\,a)$, $[-a,\,a]$, $(a,\,\infty)$ and choose $a\in\R$ such that the tail probability satisfies $\probabilityMeasure(X\geq a)\leq\epsilon$ for some small $\epsilon > 0$. Due to the small probability of the tail elements, one performs the ME method only on $[-a,\,a]$. \refthree{This technique will introduce a small error by restricting the domain to $[-a,\,a]$. Another strategy is to map the random variable to the interval $[0,\,1]$ via the cumulative distribution function, then applying the \ME{} approach on $[0,\,1]$, and finally mapping it back using the quantile, i.e. the inverse cumulative distribution function.}
\end{remark}} 
Moreover, we  introduce an indicator variable $\chi_{\MEIndex}: \Omega \to  \{0,1\}$ on every random element
\begin{equation}\label{def:indicator}
\indicatorVar{\MEIndex}{\uncertainty} :=
\begin{cases}
1&\text{if } \uncertainty \in \randomElement{\MEIndex}, \\
0 &\text{else, }
\end{cases}
\end{equation}
\refone{with $\MEIndex=1,\ldots,\MEElements $.} 
If we let $\{\localxiBasisPoly{\SGsumIndex}{\MEIndex} (\uncertainty)\}_{\SGsumIndex=0}^\infty$ be  orthonormal polynomials with respect to a conditional probability density function on the \ME{} $\randomElement{\MEIndex}$, \refone{$j=1,\ldots,\MEElements$} , 
the global approximation \eqref{gPCApprox} can be written as in \cite{Meyer2019a}
\begin{equation}\label{def:globalSG}
\solution(t,x,\uncertainty) = \sum_{\MEIndex=1}^{\MEElements}  \solution_\MEIndex(t,x,\uncertainty) \indicatorVar{\MEIndex}{\uncertainty} \approx \sum_{\MEIndex=1}^{\MEElements}  \sum_{\SGsumIndex=0}^{\SGtruncorder}  \solution_{\SGsumIndex, \MEIndex}(\timevar,\x) \localxiBasisPoly{\SGsumIndex}{\MEIndex}(\uncertainty) \indicatorVar{\MEIndex}{\uncertainty}.
\end{equation}
As $\MEElements,\SGtruncorder\to\infty$, the local approximation converges to the global solution in $\Lp{2}(\sampleSpace)$, cf. \cite{Alpert1993}.

The calculation of the expected value and the variance can be found in \cite{Meyer2019a}.

\refone{The stochastic Galerkin scheme is then applied to every \ME{} separately due to the disjoint decomposition of the random space and we come up with the following ME stochastic Galerkin system
\begin{equation}\label{SGsystemME}\pd{t} \solution_{\SGeqIndex,\MEIndex} +  \pd{\x} \int_{\cellR{\cellindR}}\! \flux\Big(\!\sum_{\MEIndex=1}^{\MEElements}  \sum_{\SGsumIndex=0}^{\SGtruncorder}  \solution_{\SGsumIndex, \MEIndex}(\timevar,\x) \localxiBasisPoly{\SGsumIndex}{\MEIndex}(\uncertainty) \Big) \xiBasisPoly{\SGeqIndex,\MEIndex}\xiPDFdxi = 0,\end{equation}
for $\SGeqIndex=0,\ldots,\SGtruncorder$ and $\MEIndex=1,\ldots,\MEElements$.
}
\refthree{Comparing this ME stochastic Galerkin system to sG without the ME approach in \eqref{SGsystem}, we observe that the complexity increases since the sG system has now to be solved in each \ME{} $\cellR{\cellindR}$, for $\cellindR=1,\ldots,\MEElements$.}

\subsection{Propagation of the Gibbs phenomenon}\label{sec:LAex}
\refone{In the following section we want to illustrate oscillations that appear within a simple example of an uncertain scalar linear advection equation due to the presence of Gibbs phenomenon, which arises when discontinuous data is interpolated with orthogonal basis functions. Further examples to this topic can be found in \cite{Poette2009,Barth2012,Ma2009}. We additionally show that the multi-element method from \secref{sec:ME} is not sufficient to prevent oscillations in our test case.}

As introductory example, we consider the one-dimensional hyperbolic problem
\begin{equation}\label{stochasticWaveEquation}
 \frac{\partial}{\partial t}u + a(\xi) \, \frac{\partial}{\partial x} u = 0
\end{equation}
with $\dimension=1$, $x\in(0,\infty)$ and uncertain wave speed
\begin{equation}\label{stochasticWaveSpeed}
 a(\xi) = 1.5 + 0.5\, \xi,
\end{equation}
where $\xi \sim U(-1,1)$, i.e., $\xi$ is uniformly distributed in $[-1,1]$. Hence, the density function is given by \refone{$\xiPDF(\uncertainty)=\frac{1}{2}\indicatorVar{(-1,1)}{\uncertainty}$} and the basis functions $\xiBasisPoly{k}$, \refone{$k=0,\ldots,\SGtruncorder$,} by the Legendre polynomials orthonormalized with respect to \eqref{eq:orthogonality}.
We further use non-smooth initial data
\refone{\begin{equation}\label{sweInitialConditions}
 u(0,\x,\xi) =  u(0,\x) = \begin{cases}
           1 & \text{for } 0 \le x \le 0.5,\\
	   0 & \text{for } 0.5 < x,
          \end{cases}
\end{equation}}
and constant boundary data
\begin{equation}\label{sweBoundaryConditions}
 u(t,0,\xi) = 1.
\end{equation}

We approximate the solution by the Polynomial Chaos (gPC) expansion \eqref{gPCApprox}
\begin{equation}\label{LAgPC}
 u(t,x,\xi) \approx \sum\limits_{\SGsumIndex=0}^{\SGtruncorder} u_\SGsumIndex(\timevar,\x) \, \xiBasisPoly{\SGsumIndex}(\uncertainty).
\end{equation}
Inserting \eqref{LAgPC} in \eqref{stochasticWaveEquation} and applying a Galerkin projection in the stochastic space leads to
\begin{equation*}
 \pd{t} u_l(t,x) + \sum\limits_{k=0}^{\SGtruncorder} a_{l,k} \pd{x} u_k(t,x) = 0, \qquad \text{for~}l = 0, \dots, \SGtruncorder,
\end{equation*}
with
\begin{equation*}
 a_{l,k} = \int_{-1}^{1} a(\xi) \xiBasisPoly{l}(\xi) \xiBasisPoly{k}(\xi) \xiPDF(\xi) \dint{\xi},  \qquad \refone{\text{for~}l,k = 0, \dots, \SGtruncorder.}
\end{equation*}

Collecting all $u_l$, \refone{$l=0,\ldots,\SGtruncorder$,} into a vector $U = (u_0, \dots, u_\SGtruncorder)^T$ and all $a_{l,k}$, \refone{$l,k=0,\dots,\SGtruncorder$,} into a matrix $A$ yields the system
\begin{equation}\label{sweProjection}
 \pd{t} U(t,x) + A \pd{x} U(t,x) = 0.
\end{equation}
Note that this system is hyperbolic because $A = A^T$. \refone{Moreover, we want to point out that the sG method is closely related to the
well-known $P_N$ closure from kinetic theory, cf. \cite{Lewis-Miller-1984}.}

The Galerkin projection of the initial conditions \eqref{sweInitialConditions} gives
\begin{equation*}
u_l(0,x) = \int_{-1}^{1} \refone{u(0,x,\xi)} \xiBasisPoly{l}(\xi) \xiPDF(\xi) \dint{\xi} = \refone{u(0,x)} \int_{-1}^{1} \!\xiBasisPoly{l}(\xi) \xiPDF(\xi) \dint{\xi}, \quad \quad \text{for} \ \ l = 0, \dots, \SGtruncorder
\end{equation*}
and therewith
\begin{equation*}
 U(0,x) = \begin{pmatrix}
           \refone{u(0,x)}\\ 0\\ \vdots\\ 0
          \end{pmatrix}\,.
\end{equation*}
For the boundary conditions \eqref{sweBoundaryConditions}, the Galerkin projection leads to
\begin{equation*}
 u_l(t,0,\xi) = \int_{-1}^{1} u(t,0,\xi) \xiBasisPoly{l}(\xi) \xiPDF(\xi) \dint{\xi} = \int_{-1}^{1} \xiBasisPoly{l}(\xi) \xiPDF(\xi) \dint{\xi}, \quad \quad \text{for} \ \ l = 0, \dots, \SGtruncorder
\end{equation*}
hence
\begin{equation*}
 U(t,0) = \begin{pmatrix}
           1\\ 0\\ \vdots\\ 0
          \end{pmatrix}\,.
\end{equation*}

Since \eqref{sweProjection} is a linear hyperbolic system, it can be decoupled into $\SGtruncorder$ scalar transport equations and the analytic solution can be explicitly given. For this purpose, we introduce a new variable $W := T^{-1} U$, where $T$ is the matrix of eigenvectors of $A$ (columnwise). With $\lambda_j$\refone{, $j=0,\ldots,\SGtruncorder$,} being the corresponding (real) eigenvalues, the components of $W$ satisfy
\begin{equation*}
 \pd{t} W_j + \lambda_j \pd{x} W_j = 0, \quad \quad \text{for} \ \ j = 0, \dots, \SGtruncorder
\end{equation*}
and the solution of the Galerkin system \refone{can be calculated by the method of characteristics \cite{evans2010partial} to}
\begin{equation}\label{eq:solLASG}
 W_j(t,x) = \begin{cases}
             \big( T^{-1} U(0,x-\lambda_j t) \big)_j\quad & \text{for } x-\lambda_j t \ge 0,\\
	     \big( T^{-1} U(t-\frac{x}{\lambda_j},0) \big)_j & \text{for } x-\lambda_j t \le 0,
            \end{cases}
\end{equation}
with \refone{$j=0,\ldots,\SGtruncorder$.}

For the following considerations, we use $\SGtruncorder=2$. From \eqref{eq:orthogonality} we deduce the Legendre polynomials
\begin{equation}\label{eq:legendrepolys}
 \xiBasisPoly{0}(\xi) = 1\,, \quad \xiBasisPoly{1}(\xi) = \sqrt{3} \xi\,, \quad \xiBasisPoly{2}(\xi) = \tfrac{\sqrt{5}}{2} (3\xi^2-1)\,,
\end{equation}
and further
\begin{equation*}
 A = \begin{pmatrix}
      \frac32 & \frac{\sqrt{3}}6 & 0\\
      \frac{\sqrt{3}}6 & \frac32 & \frac{\sqrt{15}}{15}\\
      0 & \frac{\sqrt{15}}{15} & \frac32
     \end{pmatrix}\,.
\end{equation*}
The eigenvalues of $A$ are $\tfrac32 - \tfrac{\sqrt{15}}{10}$, $\tfrac32$ and $\tfrac32 + \tfrac{\sqrt{15}}{10}$. 
Interestingly, \refone{if we choose $\xi \in \big\{-\sqrt{\tfrac{3}{5}}, 0, \sqrt{\tfrac{3}{5}}\big\}$, we obtain via \eqref{stochasticWaveSpeed} the corresponding stochastic wave speed $a(\xi)\in \big\{\tfrac32 - \tfrac{\sqrt{15}}{10}, \tfrac32, \tfrac32 + \tfrac{\sqrt{15}}{10}\big\}$, i.e. the eigenvalues of $A$. Thus, the analytical solution of \eqref{eq:spde} would coincide with the stochastic Galerkin solution \eqref{eq:solLASG}. Since $\big\{-\sqrt{\tfrac{3}{5}}, 0, \sqrt{\tfrac{3}{5}}\big\}$ represent the quadrature nodes for a third order Gauß-Legendre quadrature, the gPC Galerkin approximation \eqref{gPCApprox} equals the quadratic interpolating polynomial in this case.}
Now, the problem is that the discontinuity with respect to $x$ in the initial conditions \eqref{sweInitialConditions} carries over to the stochastic domain. This leads to overshoots (Gibbs phenomena) in the gPC Galerkin approximation as illustrated in Figure~\ref{fig:gPCApprox} for $t=0.5$ and a choice of $\xi$s inside and outside of the convex hull of the sampling points.

\externaltikz{figureGPCx}{\def\xmin{0.8}
\def\xmax{1.7}
\def\ymin{-0.4}
\def\ymax{1.4}
\def\width{\textwidth}
\def\height{1.25\textwidth}
\def\xtick{1.0, 1.5}

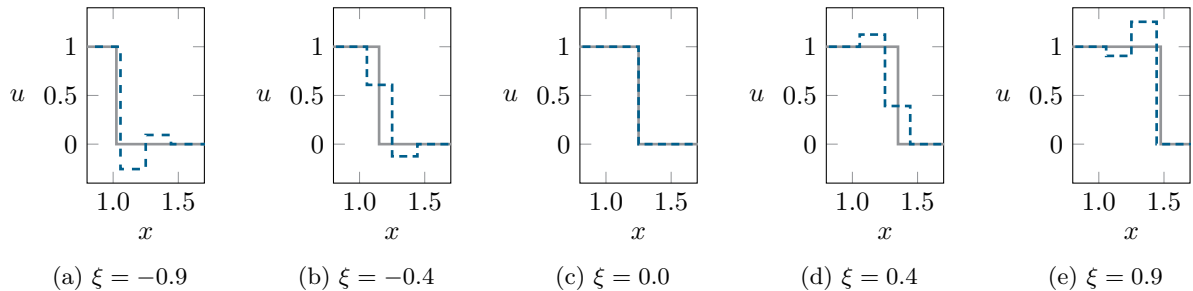
\begin{figure}[htb]
\begin{subfigure}[c]{0.19\textwidth}
 \begin{tikzpicture}
  \begin{axis}[width=\width, height=\height, xlabel=$x$,  ylabel=$u$, ylabel style = {rotate=-90}, xmin=\xmin, xmax=\xmax, ymin=\ymin, ymax=\ymax, xtick={\xtick}, xticklabels={{$1.0$},{$1.5$}}] 
   \addplot[mark=none, solid, color=gray, line width=1pt] table{data/gPCm09_exact.table};
   \addplot[mark=none, tuklblue, dashed, line width=1pt] table{data/gPCm09_approx.table};
  \end{axis}
 \end{tikzpicture}
 \subcaption{$\xi = -0.9$}
\end{subfigure}
\begin{subfigure}[c]{0.19\textwidth}
 \begin{tikzpicture}
  \begin{axis}[width=\width, height=\height, xlabel=$x$,  ylabel=$u$, ylabel style = {rotate=-90},
  xmin=\xmin, xmax=\xmax, ymin=\ymin, ymax=\ymax, xtick={\xtick}, xticklabels={{$1.0$},{$1.5$}}] 
   \addplot[mark=none, solid, color=gray, line width=1pt] table{data/gPCm04_exact.table};
   \addplot[mark=none, tuklblue, densely dashed, line width=1pt] table{data/gPCm04_approx.table};
  \end{axis}
 \end{tikzpicture}
 \subcaption{$\xi = -0.4$}
\end{subfigure}
\begin{subfigure}[c]{0.19\textwidth}
 \begin{tikzpicture}
  \begin{axis}[width=\width, height=\height, xlabel=$x$,  ylabel=$u$, ylabel style = {rotate=-90}, xmin=\xmin, xmax=\xmax, ymin=\ymin, ymax=\ymax, xtick={\xtick}, xticklabels={{$1.0$},{$1.5$}}] 
   \addplot[mark=none, solid, color=gray, line width=1pt] table{data/gPC00_exact.table};
   \addplot[mark=none, tuklblue, densely dashed, line width=1pt] table{data/gPC00_approx.table};
  \end{axis}
 \end{tikzpicture}
 \subcaption{$\xi = 0.0$}
\end{subfigure}
\begin{subfigure}[c]{0.19\textwidth}
 \begin{tikzpicture}
  \begin{axis}[width=\width, height=\height, xlabel=$x$,  ylabel=$u$, ylabel style = {rotate=-90}, xmin=\xmin, xmax=\xmax, ymin=\ymin, ymax=\ymax, xtick={\xtick}, xticklabels={{$1.0$},{$1.5$}}] 
   \addplot[mark=none, solid, color=gray, line width=1pt] table{data/gPC04_exact.table};
   \addplot[mark=none, tuklblue, densely dashed, line width=1pt] table{data/gPC04_approx.table};
  \end{axis}
 \end{tikzpicture}
 \subcaption{$\xi = 0.4$}
\end{subfigure}
\begin{subfigure}[c]{0.19\textwidth}
 \begin{tikzpicture}
  \begin{axis}[width=\width, height=\height, xlabel=$x$,  ylabel=$u$, ylabel style = {rotate=-90}, xmin=\xmin, xmax=\xmax, ymin=\ymin, ymax=\ymax, xtick={\xtick}, xticklabels={{$1.0$},{$1.5$}}] 
   \addplot[mark=none, solid, color=gray, line width=1pt] table{data/gPC09_exact.table};
   \addplot[mark=none, tuklblue, densely dashed, line width=1pt] table{data/gPC09_approx.table};
  \end{axis}
 \end{tikzpicture}
 \subcaption{$\xi = 0.9$}
\end{subfigure}

 \caption{Comparison of the gPC Galerkin approximation (dashed) and the exact solution (solid) for the given $\xi$s at $t=0.5$.}
 \label{fig:gPCApprox}
\end{figure}}

Note that this effect does not result from the fact that we made a global ansatz in the stochastic space but from the lack of regularity (see below). Moreover, if we want to refine the grid in the stochastic space close to the discontinuity, we are faced with the problem that the position of the jump in the stochastic space depends on $x$ and $t$ as illustrated in \figref{fig:discStochasticSpace}. A short computation shows that the discontinuity is located at $$\xi_{\text{jump}} = \frac{2x-1}{t} - 3,$$ which yields that for $t=0.5$, the position of the discontinuity in the stochastic space takes all values in $[-1,1]$ for $x \in [1,1.5]$.

\externaltikz{figureGPCxi}{\def\ximin{-1}
\def\ximax{1}
\def\xitick{-1,0,1}
\def\ymin{-0.4} 
\def\ymax{1.4} 
\def\width{\textwidth}
\def\height{1.25\textwidth}
\def\col{gray}

\begin{figure}[htb]
\begin{subfigure}[c]{0.19\textwidth}
 \begin{tikzpicture}
  \begin{axis}[width=\width, height=\height, xlabel=$\xi$, ylabel=$u$, ylabel style = {rotate=-90}, xmin=\ximin, xmax=\ximax, ymin=\ymin, ymax=\ymax, xtick=\xitick]
   \addplot[mark=none, solid, color=\col, line width=1pt] table{data/jumpPositionX095T05.table};
  \end{axis}
 \end{tikzpicture}
 \subcaption{$x = 0.95$}
\end{subfigure}
\begin{subfigure}[c]{0.19\textwidth}
 \begin{tikzpicture}
  \begin{axis}[width=\width, height=\height, xlabel=$\xi$, ylabel=$u$, ylabel style = {rotate=-90}, xmin=\ximin, xmax=\ximax, ymin=\ymin, ymax=\ymax, xtick=\xitick]
   \addplot[mark=none, solid, color=\col, line width=1pt] table{data/jumpPositionX110T05.table};
  \end{axis}
 \end{tikzpicture}
 \subcaption{$x = 1.1$}
\end{subfigure}
\begin{subfigure}[c]{0.19\textwidth}
 \begin{tikzpicture}
  \begin{axis}[width=\width, height=\height, xlabel=$\xi$, ylabel=$u$, ylabel style = {rotate=-90}, xmin=\ximin, xmax=\ximax, ymin=\ymin, ymax=\ymax, xtick=\xitick]
   \addplot[mark=none, solid, color=\col, line width=1pt] table{data/jumpPositionX125T05.table};
  \end{axis}
 \end{tikzpicture}
 \subcaption{$x = 1.25$}
\end{subfigure}
\begin{subfigure}[c]{0.19\textwidth}
 \begin{tikzpicture}
  \begin{axis}[width=\width, height=\height, xlabel=$\xi$, ylabel=$u$, ylabel style = {rotate=-90}, xmin=\ximin, xmax=\ximax, ymin=\ymin, ymax=\ymax, xtick=\xitick]
   \addplot[mark=none, solid, color=\col, line width=1pt] table{data/jumpPositionX140T05.table};
  \end{axis}
 \end{tikzpicture}
 \subcaption{$x = 1.4$}
\end{subfigure}
\begin{subfigure}[c]{0.19\textwidth}
 \begin{tikzpicture}
  \begin{axis}[width=\width, height=\height, xlabel=$\xi$, ylabel=$u$, ylabel style = {rotate=-90}, xmin=\ximin, xmax=\ximax, ymin=\ymin, ymax=\ymax, xtick=\xitick]
   \addplot[mark=none, solid, color=\col, line width=1pt] table{data/jumpPositionX155T05.table};
  \end{axis}
 \end{tikzpicture}
 \subcaption{$x = 1.55$}
\end{subfigure}

 \caption{Plots of the exact solution depending on $\xi$ and the given $x$ coordinates at $t=0.5$.}
 \label{fig:discStochasticSpace}
\end{figure}
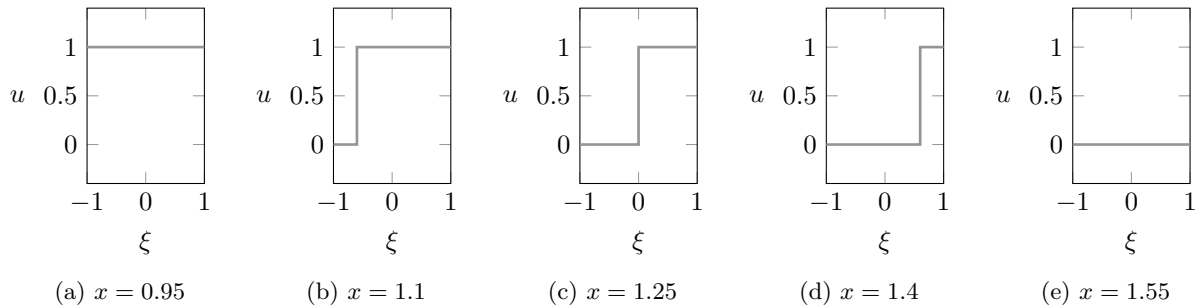}

Now, as noted above, we demonstrate that a grid refinement in the sense of a \ME{} approach does only marginally improve the so-found approximation. For this, we briefly consider the case \refone{of 3 \ME{}s,} where the stochastic domain is equally partitioned into \refone{the three intervals $\randomElement{1}=[-1,-\frac{1}{3}]$, $\randomElement{2}=[-\frac{1}{3},\frac{1}{3}]$, $\randomElement{3}=[\frac{1}{3},1]$} and we apply the gPC Galerkin approach in each of them. 
We look at the so-found \refone{\ME{}} approximations \refone{in $\randomElement{2}= [-\tfrac13, \tfrac13]$} and compare the results with the global gPC Galerkin approximations of $[-1,1]$, computed without the \ME{} approach. 
\refthree{For the local random variable $\xiv\in[-\tfrac13, \tfrac13]$, we now consider relative positions within the reference \ME{} $\widetilde{\Xi}=\frac{1}{\MEElements}[-1,1]$ compared to the positions in \figref{fig:gPCApprox}. We therefore introduce the following notation and write $3\xiv=\pm0.9$, $3\xiv=\pm0.4$ and $3\xiv=0.0$ for the point values $\xiv=\pm0.3$, $\xiv=\pm0.1\bar{3}$ and $\xiv=0$ in combination with the 3 \ME{} ansatz. At these relative positions we find again the overshoots observed for $\xi\in[-1,1]$ in \figref{fig:gPCApprox}. Note that the basis functions of the global gPC approximation $\xiBasisPoly{\SGsumIndex}(\xi)$, $\xi\in[-1,1]$, coincide with the ME basis polynomials $\xiBasisPoly{\SGsumIndex,2}(3\xiv)$, $\xiv\in[-\tfrac13, \tfrac13]$ for $\SGsumIndex=1,\ldots,\SGtruncorder$ and the ME approach describes a global gPC approximation restricted to $[-\tfrac13, \tfrac13]$, which is why we retrieve the oscillations at $3\xiv\in\{\pm0.9,\pm0.4,0\}$. }
\refone{As one can easily see from the plots in \figref{fig:gPCRefined}, the width of the overshoots in the compared approximations reduced but the height remained, which is a well-known issue for Gibbs phenomena. Thus, solely the refinement of the grid in the stochastic space cannot reduce the overshoots in the solution manifold.}

\externaltikz{figureGPCrefined}{\def\xmin{0.8}
\def\xmax{1.7}
\def\ymin{-0.4}
\def\ymax{1.4}
\def\width{\textwidth}
\def\height{1.25\textwidth}
\def\xtick{1.0, 1.5}

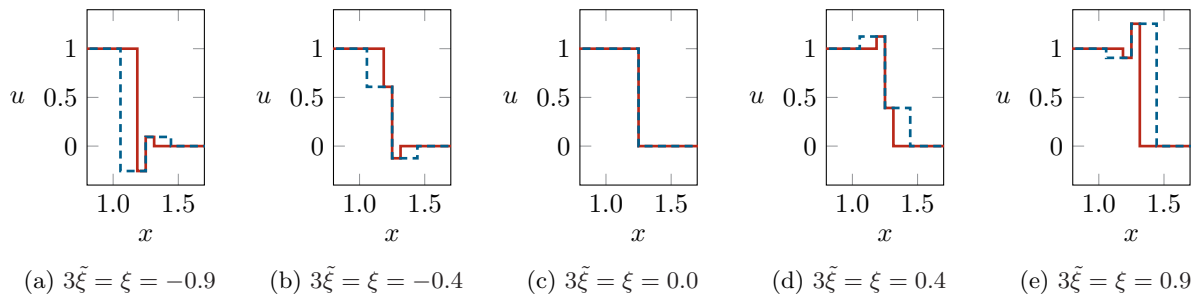
\begin{figure}[htb]
\begin{subfigure}[c]{0.19\textwidth}
 \begin{tikzpicture}
  \begin{axis}[width=\width, height=\height, xlabel=$x$, ylabel=$u$, ylabel style = {rotate=-90}, xmin=\xmin, xmax=\xmax, ymin=\ymin, ymax=\ymax, xtick={\xtick}, xticklabels={{$1.0$},{$1.5$}}] 
     \addplot[mark=none, tuklred, line width=1pt] table{data/gPCm09ref_approx.table};
   \addplot[mark=none, solid, tuklblue, densely dashed, line width=1pt] table{data/gPCm09_approx.table};
  \end{axis}
 \end{tikzpicture}
 \subcaption{\refthree{$3\xiv = \xi = -0.9$}}
\end{subfigure}
\begin{subfigure}[c]{0.19\textwidth}
 \begin{tikzpicture}
  \begin{axis}[width=\width, height=\height, xlabel=$x$, ylabel=$u$, ylabel style = {rotate=-90}, xmin=\xmin, xmax=\xmax, ymin=\ymin, ymax=\ymax, xtick={\xtick}, xticklabels={{$1.0$},{$1.5$}}] 
     \addplot[mark=none, tuklred, line width=1pt] table{data/gPCm04ref_approx.table};
   \addplot[mark=none, solid, tuklblue, densely dashed, line width=1pt] table{data/gPCm04_approx.table};
  \end{axis}
 \end{tikzpicture}
 \subcaption{\refthree{$3\xiv = \xi = -0.4$}}
\end{subfigure}
\begin{subfigure}[c]{0.19\textwidth}
 \begin{tikzpicture}
  \begin{axis}[width=\width, height=\height, xlabel=$x$, ylabel=$u$, ylabel style = {rotate=-90}, xmin=\xmin, xmax=\xmax, ymin=\ymin, ymax=\ymax, xtick={\xtick}, xticklabels={{$1.0$},{$1.5$}}] 
     \addplot[mark=none, tuklred, line width=1pt] table{data/gPC00ref_approx.table};
   \addplot[mark=none, solid, tuklblue, densely dashed, line width=1pt] table{data/gPC00_approx.table};
  \end{axis}
 \end{tikzpicture}
 \subcaption{\refthree{$3\xiv = \xi  = 0.0$}}
\end{subfigure}
\begin{subfigure}[c]{0.19\textwidth}
 \begin{tikzpicture}
  \begin{axis}[width=\width, height=\height, xlabel=$x$, ylabel=$u$, ylabel style = {rotate=-90}, xmin=\xmin, xmax=\xmax, ymin=\ymin, ymax=\ymax, xtick={\xtick}, xticklabels={{$1.0$},{$1.5$}}] 
     \addplot[mark=none, tuklred, line width=1pt] table{data/gPC04ref_approx.table};
   \addplot[mark=none, solid, tuklblue, densely dashed, line width=1pt] table{data/gPC04_approx.table};
  \end{axis}
 \end{tikzpicture}
 \subcaption{\refthree{$3\xiv = \xi = 0.4$}}
\end{subfigure}
\begin{subfigure}[c]{0.19\textwidth}
 \begin{tikzpicture}
  \begin{axis}[width=\width, height=\height, xlabel=$x$, ylabel=$u$, ylabel style = {rotate=-90}, xmin=\xmin, xmax=\xmax, ymin=\ymin, ymax=\ymax, xtick={\xtick}, xticklabels={{$1.0$},{$1.5$}}] 
     \addplot[mark=none, tuklred, line width=1pt] table{data/gPC09ref_approx.table};
   \addplot[mark=none, solid, tuklblue, densely dashed, line width=1pt] table{data/gPC09_approx.table};
  \end{axis}
 \end{tikzpicture}
 \subcaption{\refthree{$3\xiv = \xi = 0.9$}}
\end{subfigure}

 \caption{Comparison of the overshoots in the global gPC Galerkin approximation (dashed) and the approximations by the \ME{} approach (solid) for the given $\xiv$s and $\xi$s at $t=0.5$.}
 \label{fig:gPCRefined}
\end{figure}}

\subsection{Stochastic Slope Limiter} \label{sec:slopeLimiter}
We now consider a slope limiter within the stochastic space to reduce the oscillations described in the previous subsection. We present the minimod limiter \cite{Cockburn1989a,Cockburn1991}, given by the following troubled cell indicator for the $j$th \ME{} $\randomElement{\MEIndex}$
\begin{equation}\label{eq:TCminmod}
TC_j(\solution) = \begin{cases}
1 & \solution_{1,j} \neq \text{m}\big(\solution_{1,j},\, \solution_{0,j+1}-\solution_{0,j}, \, \solution_{0,j}-\solution_{0,j-1}\big),\\
0 & \text{else},
\end{cases}
\end{equation}
where $\text{m}(\cdot,\cdot,\cdot)$ is the minmod function, $\solution_{k,j}$, $k=0,1$, are the $k$th coefficients of the local gPC representations of the solution $\solution$ in the $j$th \ME{} as in \eqref{def:globalSG} \refone{and $j=1,\ldots,\MEElements$}. At the boundary, we copy the entries of the first (last) \ME{}. Finally, we replace each coefficient of $\solution\big|_{\randomElement{\MEIndex}}$ for $j=1,\ldots,\MEElements$ by
\begin{equation}\label{eq:slopeLimiter}
\refone{\stochsllimiter{\solution}}\big|_{\randomElement{\MEIndex}} =\stochsllimiter{\begin{matrix}
	\solution_{0,j}\\ \solution_{1,j}\\ \vdots \\ \solution_{\SGtruncorder,j}
	\end{matrix}}
= \begin{cases}
\begin{pmatrix}
(\solution_{0,j})^T\\[0.2cm]
 \text{m}\big(\solution_{1,j},\, \solution_{0,j+1}-\solution_{0,j}, \, \solution_{0,j}-\solution_{0,j-1}\big)\\[0.2cm]
 (0,\ldots,0)^T\\
 \vdots\\
 (0,\ldots,0)^T
\end{pmatrix} & \text{if } TC_j(\solution) = 1, \\[1.5cm]
\quad \solution\big|_{\randomElement{\MEIndex}} & \text{else}.
\end{cases}\end{equation}
Hence, we assume that the oscillations in the solution vector are mainly generated in the part with linear uncertainty \cite{Shu1987}.
\refthree{Note that the cell mean $ \solution_{0,j}$, $j=1,\ldots,\MEElements$, is not changed through the limiter \eqref{eq:slopeLimiter}.}

\begin{remark}
	The troubled cell indicator \eqref{eq:TCminmod} is defined for an expansion of $\solution$ in \refone{the space of monomials $\{1,x,x^2,\ldots,x^{\SGtruncorder}\}$. For other types of sG basis polynomials we need to adapt this definition using a basis transformation.
	In particular, if $\xiBasisPoly{k}$, $k=0,\ldots,\SGtruncorder$, are given by orthonormal basis polynomials in the spirit of \eqref{eq:orthogonality}, we obtain the corresponding coefficients for an expansion in monomials by multiplication of the matrix $V^{-1}$, \refthree{where} $V=(\int_\randomSpace x^m\xiBasisPoly{k}\xiPDFdxi)_{m,k=0:\SGtruncorder}$. After the application of the slope limiter, we transform the coefficients back by a multiplication of $V$.}
\end{remark}

If $\SGtruncorder\geq2$, we only apply the slope limiter if 
$|\solution_{1,j}| \geq M | \randomElement{\MEIndex} |^2$
in addition to $TC_j(\solution) = 1$ to achieve the TVBM property \refone{in every \ME{} $\randomElement{\MEIndex}$ for $\MEIndex=1,\ldots,\MEElements$}. The constant $M$ is chosen according to \cite{quarteroni2006advanced} as
$$M = \text{sup}\Big\{\big|\partial_\xi^2\refone{\solution(0,\x,\hat{\uncertainty})}\big| ~\Big|~ \hat{\uncertainty}\in\SD,\x\in\Domain,\,\partial_\xi\refone{\solution(0,\x,\hat{\uncertainty})}=0\Big\}.$$

The slope limited gPC Galerkin approximation can now be applied to the linear advection example from Section \ref{sec:LAex}. \figref{fig:gPClimiter}  compares the gPC approximation in 3 \ME{}s with its slope limited version which completely reduces the overshoots. Thus, the application of the minmod slope limiter in the stochastic space \eqref{eq:slopeLimiter} has been sufficient to prevent Gibbs phenomenon.

\externaltikz{figureGPClimiter}{\def\xmin{0.8}
\def\xmax{1.7}
\def\ymin{-0.4}
\def\ymax{1.4}
\def\width{\textwidth}
\def\height{1.25\textwidth}
\def\xtick{1.0, 1.5}

\begin{figure}[htb]
\begin{subfigure}[c]{0.19\textwidth}
 \begin{tikzpicture}
  \begin{axis}[width=\width, height=\height, xlabel=$x$, ylabel=$u$, ylabel style = {rotate=-90},
   xmin=\xmin, xmax=\xmax, ymin=\ymin, ymax=\ymax, xtick={\xtick}, xticklabels={{$1.0$},{$1.5$}}] 
   \addplot[mark=none, solid, tuklred, line width=1pt] table{data/gPCm09ref_approx.table};
   \addplot[mark=none, green, dashed, line width=1pt] file {Images/figureGPClimiter/GPClimiter_xim03.txt};
  \end{axis}
 \end{tikzpicture}
 \subcaption{\refthree{$3\xiv = -0.9$}}
\end{subfigure}
\begin{subfigure}[c]{0.19\textwidth}
 \begin{tikzpicture}
  \begin{axis}[width=\width, height=\height, xlabel=$x$, ylabel=$u$, ylabel style = {rotate=-90}, xmin=\xmin, xmax=\xmax, ymin=\ymin, ymax=\ymax, xtick={\xtick}, xticklabels={{$1.0$},{$1.5$}}] 
   \addplot[mark=none, solid, tuklred, line width=1pt] table{data/gPCm04ref_approx.table};
   \addplot[mark=none, green, densely dashed, line width=1pt] file {Images/figureGPClimiter/GPClimiter_xim013.txt};
  \end{axis}
 \end{tikzpicture}
 \subcaption{\refthree{$3\xiv = -0.4$}}
\end{subfigure}
\begin{subfigure}[c]{0.19\textwidth}
 \begin{tikzpicture}
  \begin{axis}[width=\width, height=\height, xlabel=$x$, ylabel=$u$, ylabel style = {rotate=-90}, xmin=\xmin, xmax=\xmax, ymin=\ymin, ymax=\ymax, xtick={\xtick}, xticklabels={{$1.0$},{$1.5$}}] 
   \addplot[mark=none, solid, tuklred, line width=1pt] table{data/gPC00ref_approx.table};
   \addplot[mark=none, green, densely dashed, line width=1pt] file {Images/figureGPClimiter/GPClimiter_xi0.txt};
  \end{axis}
 \end{tikzpicture}
 \subcaption{\refthree{$3\xiv =  0.0$}}
\end{subfigure}
\begin{subfigure}[c]{0.19\textwidth}
 \begin{tikzpicture}
  \begin{axis}[width=\width, height=\height, xlabel=$x$,ylabel=$u$, ylabel style = {rotate=-90},  xmin=\xmin, xmax=\xmax, ymin=\ymin, ymax=\ymax, xtick={\xtick}, xticklabels={{$1.0$},{$1.5$}}] 
   \addplot[mark=none, solid, tuklred, line width=1pt] table{data/gPC04ref_approx.table};
   \addplot[mark=none, green, densely dashed, line width=1pt] file {Images/figureGPClimiter/GPClimiter_xi013.txt};
  \end{axis}
 \end{tikzpicture}
 \subcaption{\refthree{$3\xiv =0.4$}}
\end{subfigure}
\begin{subfigure}[c]{0.19\textwidth}
 \begin{tikzpicture}
  \begin{axis}[width=\width, height=\height, xlabel=$x$, ylabel=$u$, ylabel style = {rotate=-90}, xmin=\xmin, xmax=\xmax, ymin=\ymin, ymax=\ymax, xtick={\xtick}, xticklabels={{$1.0$},{$1.5$}}] 
   \addplot[mark=none, solid, tuklred, line width=1pt] table{data/gPC09ref_approx.table};
   \addplot[mark=none, green, densely dashed, line width=1pt] file {Images/figureGPClimiter/GPClimiter_xi013.txt};
  \end{axis}
 \end{tikzpicture}
 \subcaption{\refthree{$3\xiv = 0.9$}}
\end{subfigure}

 \caption{Comparison of the overshoots in the gPC Galerkin approximation (solid) with 3 \ME{s} and the approximations by its slope limited approach (dashed) for the given $\xiv$s at $t=0.5$.}
 \label{fig:gPClimiter}
\end{figure}
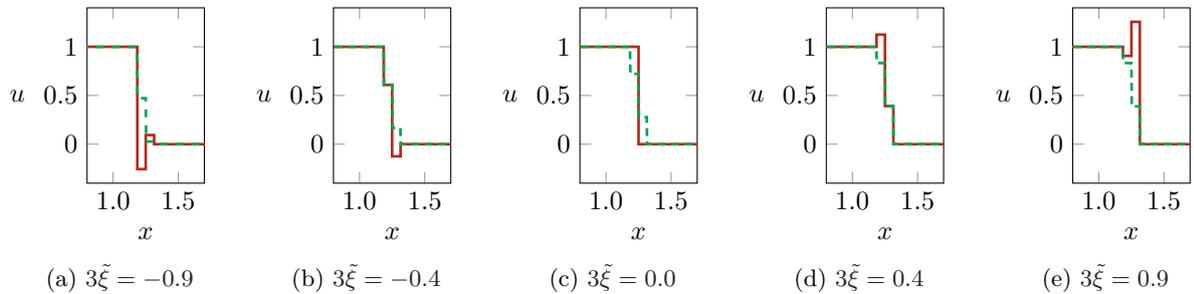}

In \figref{fig:gPClimiter_refined} we kept refining the number of \ME{}s to 10, where Gibbs phenomenon still causes overshoots of the same height as before. The slope limiter is able to eliminate these oscillations.

\externaltikz{figureGPClimiter_refined}{\def\xmin{0.8}
\def\xmax{1.7}
\def\ymin{-0.4}
\def\ymax{1.4}
\def\width{\textwidth}
\def\height{1.25\textwidth}
\def\xtick{1.0, 1.5}

\begin{figure}[htb]
\begin{subfigure}[c]{0.19\textwidth}
 \begin{tikzpicture}
  \begin{axis}[width=\width, height=\height, xlabel=$x$, ylabel=$u$, ylabel style = {rotate=-90}, xmin=\xmin, xmax=\xmax, ymin=\ymin, ymax=\ymax, xtick={\xtick}, xticklabels={{$1.0$},{$1.5$}}] 
   \addplot[mark=none, solid, tuklred, line width=1pt] file {Images/figureGPClimiter_refined/GPClimiter_xim09.txt};
   \addplot[mark=none, green, dashed, line width=1pt] file {Images/figureGPClimiter_refined/GPClimiter_xim09_SL.txt};
  \end{axis}
 \end{tikzpicture}
 \subcaption{\refthree{$10\xiv = -0.9$}}
\end{subfigure}
\begin{subfigure}[c]{0.19\textwidth}
 \begin{tikzpicture}
  \begin{axis}[width=\width, height=\height, xlabel=$x$, ylabel=$u$, ylabel style = {rotate=-90}, xmin=\xmin, xmax=\xmax, ymin=\ymin, ymax=\ymax, xtick={\xtick}, xticklabels={{$1.0$},{$1.5$}}] 
   \addplot[mark=none, solid, tuklred, line width=1pt] file {Images/figureGPClimiter_refined/GPClimiter_xim04.txt};
\addplot[mark=none, green, dashed, line width=1pt] file {Images/figureGPClimiter_refined/GPClimiter_xim04_SL.txt};
  \end{axis}
 \end{tikzpicture}
 \subcaption{\refthree{$10\xiv =-0.4$}}
\end{subfigure}
\begin{subfigure}[c]{0.19\textwidth}
 \begin{tikzpicture}
  \begin{axis}[width=\width, height=\height, xlabel=$x$, ylabel=$u$, ylabel style = {rotate=-90}, xmin=\xmin, xmax=\xmax, ymin=\ymin, ymax=\ymax, xtick={\xtick}, xticklabels={{$1.0$},{$1.5$}}] 
   \addplot[mark=none, solid, tuklred, line width=1pt] file {Images/figureGPClimiter_refined/GPClimiter_xi0.txt};
\addplot[mark=none, green, dashed, line width=1pt] file {Images/figureGPClimiter_refined/GPClimiter_xi0_SL.txt};
  \end{axis}
 \end{tikzpicture}
 \subcaption{\refthree{$10\xiv =0.0$}}
\end{subfigure}
\begin{subfigure}[c]{0.19\textwidth}
 \begin{tikzpicture}
  \begin{axis}[width=\width, height=\height, xlabel=$x$, ylabel=$u$, ylabel style = {rotate=-90}, xmin=\xmin, xmax=\xmax, ymin=\ymin, ymax=\ymax, xtick={\xtick}, xticklabels={{$1.0$},{$1.5$}}] 
   \addplot[mark=none, solid, tuklred, line width=1pt] file {Images/figureGPClimiter_refined/GPClimiter_xi04.txt};
\addplot[mark=none, green, dashed, line width=1pt] file {Images/figureGPClimiter_refined/GPClimiter_xi04_SL.txt};
  \end{axis}
 \end{tikzpicture}
 \subcaption{\refthree{$10\xiv=0.4$}}
\end{subfigure}
\begin{subfigure}[c]{0.19\textwidth}
 \begin{tikzpicture}
  \begin{axis}[width=\width, height=\height, xlabel=$x$, ylabel=$u$, ylabel style = {rotate=-90}, xmin=\xmin, xmax=\xmax, ymin=\ymin, ymax=\ymax, xtick={\xtick}, xticklabels={{$1.0$},{$1.5$}}] 
   \addplot[mark=none, solid, tuklred, line width=1pt] file {Images/figureGPClimiter_refined/GPClimiter_xi09.txt};
\addplot[mark=none, green, dashed, line width=1pt] file {Images/figureGPClimiter_refined/GPClimiter_xi09_SL.txt};
  \end{axis}
 \end{tikzpicture}
 \subcaption{\refthree{$10\xiv = 0.9$}}
\end{subfigure}

 \caption{Comparison of the overshoots in the gPC Galerkin approximation (solid) with 10 \ME{s} and the approximations by its slope limited approach (dashed) for the given $\xiv$s at $t=0.5$.}
 \label{fig:gPClimiter_refined}
\end{figure}
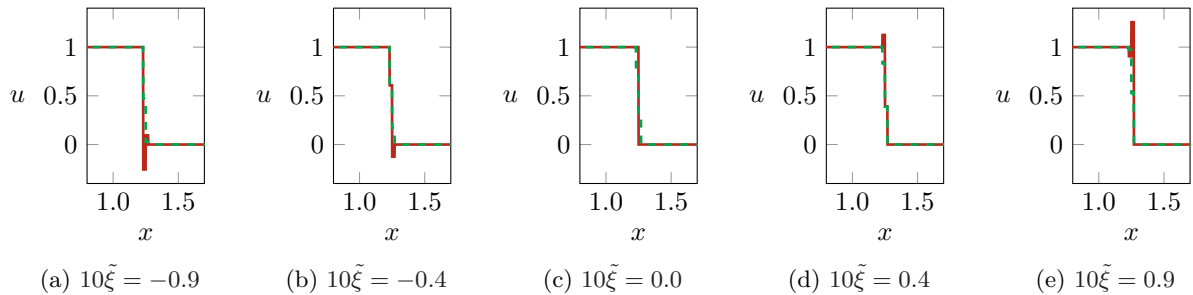}

\refthree{Moreover, in \figref{fig:gPCfixedx}, we show a comparison of the standard gPC approximation to its 3 \ME{} approach as well as the limited solution with 3 \ME{}s for some fixed values of $\x$ over the $\uncertainty$ domain. We observe oscillations in the standard stochastic Galerkin solution that have propagated from the spatial domain into the uncertainty, whereas the \ME{} approach describes a better approximation of the exact solution that tends to overshoots at the boundaries of the \ME{}s. The slope limiter is applied in the $\uncertainty$ manifold where it is able to reduce the overshoots shown in \figref{fig:gPCfixedx}, which also helps to deal with oscillations in $\x$ as presented in \figref{fig:gPClimiter}- \ref{fig:gPClimiter_refined}.}

\externaltikz{figureGPCxireview2}{\def\ximin{-1}
\def\ximax{1}
\def\xitick{-1,0,1}
\def\ymin{-0.4} 
\def\ymax{1.4} 
\def\width{\textwidth}
\def\height{1.25\textwidth}
\def\col{gray}

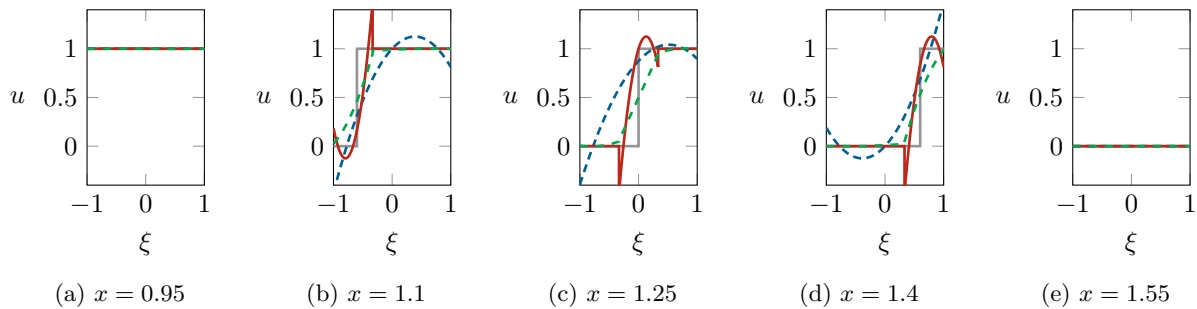
\begin{figure}[htb]
\begin{subfigure}[c]{0.19\textwidth}
 \begin{tikzpicture}
  \begin{axis}[width=\width, height=\height, xlabel=$\xi$, ylabel=$u$, ylabel style = {rotate=-90}, xmin=\ximin, xmax=\ximax, ymin=\ymin, ymax=\ymax, xtick=\xitick]
   \addplot[mark=none, solid, color=\col, line width=1pt] table{data/jumpPositionX095T05.table};
   \addplot[mark=none, solid, tuklblue, densely dashed, line width=1pt] table{Images/SGfixedx2/SGxa.txt};
   \addplot[mark=none, solid, tuklred, line width=1pt] table{Images/SGfixedx2/MESGxa.txt};
   \addplot[mark=none, solid, green, dashed, line width=1pt] table{Images/SGfixedx2/MESGlimxa.txt};
  \end{axis}
 \end{tikzpicture}
 \subcaption{$x = 0.95$}
\end{subfigure}
\begin{subfigure}[c]{0.19\textwidth}
 \begin{tikzpicture}
  \begin{axis}[width=\width, height=\height, xlabel=$\xi$, ylabel=$u$, ylabel style = {rotate=-90}, xmin=\ximin, xmax=\ximax, ymin=\ymin, ymax=\ymax, xtick=\xitick]
   \addplot[mark=none, solid, color=\col, line width=1pt] table{data/jumpPositionX110T05.table};
   \addplot[mark=none, solid, tuklblue, densely dashed, line width=1pt] table{Images/SGfixedx2/SGxb.txt};
  \addplot[mark=none, solid, tuklred, line width=1pt] table{Images/SGfixedx2/MESGxb.txt};
  \addplot[mark=none, solid, green, dashed, line width=1pt] table{Images/SGfixedx/MESGlimxb.txt};
  \end{axis}
 \end{tikzpicture}
 \subcaption{$x = 1.1$}
\end{subfigure}
\begin{subfigure}[c]{0.19\textwidth}
 \begin{tikzpicture}
  \begin{axis}[width=\width, height=\height, xlabel=$\xi$, ylabel=$u$, ylabel style = {rotate=-90}, xmin=\ximin, xmax=\ximax, ymin=\ymin, ymax=\ymax, xtick=\xitick]
   \addplot[mark=none, solid, color=\col, line width=1pt] table{data/jumpPositionX125T05.table};
   \addplot[mark=none, solid, tuklblue, densely dashed, line width=1pt] table{Images/SGfixedx2/SGxc.txt};
   \addplot[mark=none, solid, tuklred, line width=1pt] table{Images/SGfixedx2/MESGxc.txt};
   \addplot[mark=none, solid, green, dashed, line width=1pt] table{Images/SGfixedx/MESGlimxc.txt};
  \end{axis}
 \end{tikzpicture}
 \subcaption{$x = 1.25$}
\end{subfigure}
\begin{subfigure}[c]{0.19\textwidth}
 \begin{tikzpicture}
  \begin{axis}[width=\width, height=\height, xlabel=$\xi$, ylabel=$u$, ylabel style = {rotate=-90}, xmin=\ximin, xmax=\ximax, ymin=\ymin, ymax=\ymax, xtick=\xitick]
   \addplot[mark=none, solid, color=\col, line width=1pt] table{data/jumpPositionX140T05.table};
   \addplot[mark=none, solid, tuklblue, densely dashed, line width=1pt] table{Images/SGfixedx2/SGxd.txt};
   \addplot[mark=none, solid, tuklred, line width=1pt] table{Images/SGfixedx2/MESGxd.txt};
   \addplot[mark=none, solid, green, dashed, line width=1pt] table{Images/SGfixedx/MESGlimxd.txt};
  \end{axis}
 \end{tikzpicture}
 \subcaption{$x = 1.4$}
\end{subfigure}
\begin{subfigure}[c]{0.19\textwidth}
 \begin{tikzpicture}
  \begin{axis}[width=\width, height=\height, xlabel=$\xi$, ylabel=$u$, ylabel style = {rotate=-90}, xmin=\ximin, xmax=\ximax, ymin=\ymin, ymax=\ymax, xtick=\xitick]
   \addplot[mark=none, solid, color=\col, line width=1pt] table{data/jumpPositionX155T05.table};
   \addplot[mark=none, solid, tuklblue, densely dashed, line width=1pt] table{Images/SGfixedx2/SGxe.txt};
   \addplot[mark=none, solid, tuklred, line width=1pt] table{Images/SGfixedx2/MESGxe.txt};
   \addplot[mark=none, solid, green, dashed, line width=1pt] table{Images/SGfixedx2/MESGlimxe.txt};
  \end{axis}
 \end{tikzpicture}
 \subcaption{$x = 1.55$}
\end{subfigure}

 \caption{\refthree{Comparison of the exact solution (solid, gray), the gPC Galerkin approximation (dashed, blue), the 3 \ME{} ansatz (solid, red) and  its slope limited approach (dashed, green) depending on $\xi$ and the given $x$ coordinates at $t=0.5$.}}
\label{fig:gPCfixedx}
\end{figure}}

\section{\scheme{} scheme (WENOsG)}\label{sec:dg}
The previous example demonstrated the significance of limiting techniques not only on the spatial but also on the stochastic grid. 

\refone{
In this section we formulate a stochastic Galerkin scheme, which is constructed to prevent the propagation of Gibbs phenomenon into the stochastic domain by applying a slope limiter in the stochasticity while preserving a high-order approximation in space and time.
If we consider systems of equations, we also have to apply the hyperbolicity-preserving limiter introduced in \cite{Meyer2019a,Schlachter2017a}, since the stochastic slope limiter is not intended to retain hyperbolicity of solutions, however, the stochastic Galerkin approach is well known to lose hyperbolicity \cite{abgrall2017uncertainty,Poette2009}. 
In \secref{sec:LAex}, we presented an example of an uncertain scalar conservation law, i.e. where the hyperbolicity-preserving limiter does not alter the solution due to the preset hyperbolicity, that still produces Gibbs oscillations which is why we additionally have to construct the stochastic slope limiter for stochastic Galerkin.
We therefore embed this slope limiter into a weighted essentially non-oscillatory stochastic Galerkin scheme, using a WENO reconstruction in the physical and stochastic domain to further dampen Gibbs oscillations within the solution manifold and to ensure a high-order resolution.
}

To this end, we subdivide the spatial domain $\Domain=[\x_L,\x_R]\subset\R$ into $\ncells$ cells $\cell{\cellind}= \big[\x_{\cellind-\frac{1}{2}},\x_{\cellind+\frac{1}{2}}\big]$, \refone{$\cellind=1,\ldots,\ncells$}, of width $\Delta\x = \frac{\x_R-\x_L}{\ncells}$.

We deploy a $\RKorder$-th order SSP Runge Kutta scheme in time and a finite volume method combined with the Weighted Essentially Non-Oscillatory (WENO) reconstruction in space (cf. \cite{Toro2009,quarteroni2006advanced,JIANG1996}) to obtain a high-order scheme for which we apply the stochastic slope limiter  \eqref{eq:slopeLimiter}. Note that in our numerical experiments, we use the so called CWENOZ method from \cite{Cravero2017}, which will be explained in \eqref{eq:WENOpoly}. 

For this sake, we consider a semi-discretized finite volume scheme for the solution $\solution$. 

\refone{We test \eqref{eq:spde} by a test function $v(\x,\uncertainty)$ with supp$(v)\subseteq\cell{\cellind}\times\cellR{\cellindR}$ and obtain, after one formal integration by parts in $\x$, the weak formulation
\begin{equation}\label{eq:weakformulation}
 \frac{\partial}{\partial t} \MEintxR{\solution(\timevar,\x,\uncertainty)\,\Testfunction(\x,\uncertainty)} = \MEintxR{\flux (\solution(\timevar,\x,\uncertainty))\,\dx\Testfunction(\x,\uncertainty)}-\MEintcellR{\cellindR}{\flux (\solution(\timevar,\x,\uncertainty))\,\Testfunction(\x,\uncertainty)}\,\bigg|_{\interface{\cellind-\frac12}}^{\interface{\cellind+\frac12}},
\end{equation}
for $\cellind=1,\ldots,\ncells$ and $\MEIndex=1,\ldots,\MEElements$.}

Since $\solution$ is discontinuous at $\interface{\cellind\pm\frac12}$, \refone{$i=1,\ldots,\ncells$,} we replace the evaluation of $\flux(\solution)$ at these points with a numerical flux function $\numericalFlux(\solution^-,\solution^+)$ that approximately solves the Riemann problem at the cell interface.
The values 
\begin{equation}\label{eq:limitinterface}
\solution^-(\refone{\timevar},\interface{\cellind+\frac12},\uncertainty) := \lim\limits_{\x\uparrow\interface{\cellind+\frac12}} \solution(\refone{\timevar},\x,\uncertainty), \qquad 
\solution^+(\refone{\timevar},\interface{\cellind+\frac12},\uncertainty) := \lim\limits_{\x\downarrow\interface{\cellind+\frac12}} \solution(\refone{\timevar},\x,\uncertainty)
\end{equation}
denote the left and right limits of the piecewise polynomial solution at the interface $\interface{\cellind+\frac12}$, \refone{for $i=1,\ldots,\ncells$}.
\refone{Note that this allows us to write down the following numerical schemes in conservative form and we choose the numerical flux to be consistent, i.e. we require $\numericalFlux(\solution,\solution)=\flux(\solution)$.}

\refone{
We choose $v=\frac{1}{\Delta\x\Delta\uncertainty}$ in \eqref{eq:weakformulation} and obtain for one time step of the forward Euler scheme
\begin{subequations} \label{eq:FVscheme}
	\begin{align} 
	\solutioncm{n+1}{\cellind,\cellindR}= \solutioncm{n}{\cellind,\cellindR}
	&- \frac{\Delta \timevar}{\Delta x}
	\MEintcellR{\cellindR}{\hat{\flux}\big(\solution(\timevar_\timeind,\x^-_{i+\frac{1}{2}},\uncertainty),\solution(\timevar_\timeind,\x^+_{i+\frac{1}{2}},\uncertainty)\big)}\\
	&+ \frac{\Delta \timevar}{\Delta x}
	\MEintcellR{\cellindR}{\hat{\flux}(\solution(\timevar_\timeind,\x^-_{i-\frac{1}{2}},\uncertainty),\solution(\timevar_\timeind,\x^+_{i-\frac{1}{2}},\uncertainty)\big)},
	\end{align}
\end{subequations}
where $\solutioncm{n}{\cellind,\cellindR}$, denotes the $\x-\uncertainty$-cell mean of $\solution$ in $\cell{\cellind}\times\cellR{\cellindR}$ at time $\timevar_\timeind$, namely
\begin{equation}\label{eq:xcellmean}
\solutioncm{n}{\cellind,\cellindR}:= \frac{1}{\Delta\x\Delta\uncertainty}\int_{\cell{\cellind}}\!\MEintcellR{\cellindR}{\solution(\timevar_\timeind,\x,\uncertainty)}\,\mathrm{d}\x,
\end{equation}
where $\cellind=1,\ldots,\ncells$ and $\MEIndex=1,\ldots,\MEElements$.
For the coefficients of the gPC polynomial \eqref{def:globalSG} in the ME stochastic Galerkin system \eqref{SGsystemME} the finite volume scheme \eqref{eq:FVscheme} reads
\begin{subequations}\label{eq:LaxFriedrichs}
	\begin{align}
	\solutioncm{n+1}{\SGsumIndex,\cellind,\cellindR}
	= \solutioncm{n}{\SGsumIndex,\cellind,\cellindR}
	&- \frac{\dtstepsize}{\Delta\x} \intRS \hat{\flux}\big(\solution(\timevar_\timeind,\x^-_{i+\frac{1}{2}},\uncertainty),\solution(\timevar_\timeind,\x^+_{i+\frac{1}{2}},\uncertainty)\big)\xiBasisPoly{\SGsumIndex,\cellindR}(\uncertainty) \xiPDF \mathrm{d}\uncertainty \\
	&+\frac{\dtstepsize}{\Delta\x} \intRS\hat{\flux}(\solution(\timevar_\timeind,\x^-_{i-\frac{1}{2}},\uncertainty),\solution(\timevar_\timeind,\x^+_{i-\frac{1}{2}},\uncertainty)\big) 
	\xiBasisPoly{\SGsumIndex,\cellindR}(\uncertainty) \xiPDF \mathrm{d}\uncertainty,
	\end{align}
\end{subequations}
where $\solutioncm{n}{\SGsumIndex,\cellind,\cellindR}$ describes the $\SGsumIndex$th coefficient of the gPC expansion \eqref{def:globalSG} at time $\timevar_\timeind$, in the \ME{} $\cellR{\cellindR}$ and taking the spatial cell mean over $\cell{\cellind}$, i.e.
$$ \solutioncm{n}{\SGsumIndex,\cellind,\cellindR}:= \frac{1}{\Delta\x\Delta\uncertainty}\int_{\cell{\cellind}}\!\MEintcellR{\cellindR}{\solution(\timevar_\timeind,\x,\uncertainty)\,\xiBasisPoly{\SGsumIndex,\cellindR}}\,\mathrm{d}\x,$$
for  $\SGsumIndex=0,\ldots,\SGtruncorder$, $\cellind=1,\ldots,\ncells$ and $\MEIndex=1,\ldots,\MEElements$.
}

In our numerical studies, we set the numerical flux function as the global Lax-Friedrichs flux
\begin{align}
\label{eq:globalLF}
\numericalFlux(\solution^-, \solution^+) = \dfrac{1}{2} \left( \flux(\solution^-) + \flux(\solution^+) - \viscosityconstant ( \solution^+ - 
\solution^-) \right).
\end{align}
The numerical viscosity constant $\viscosityconstant$ is taken as the global estimate of the  absolute value of the largest eigenvalue of the Jacobian
$\frac{\partial \flux(\solution)}{\partial \solution}$, \refone{namely
	$$\viscosityconstant=\max\bigg\{|\lambda_1|,\ldots,|\lambda_\dimension|~\bigg|~\lambda_j, j=1,\ldots,\dimension, \text{ are the eigenvalues of } \frac{\partial \flux(\solution)}{\partial \solution}\bigg\}.$$}

As a quadrature rule in the physical domain, we apply a Gau\ss{}-Lobatto rule on $\cell{\cellind}$, \refone{$i=1,\ldots,\ncells$,} with $\nbxnodes+1 $ points and weights $(\quadxpoint{\xQuadIndex},\quadxweight{\xQuadIndex})$, for $\xQuadIndex=0,\ldots,\nbxnodes$ and where $\nbxnodes$ is chosen such that the quadrature integrates the WENO polynomial (will be defined in \eqref{eq:WENOpoly}) exactly. Gau\ss{}-Lobatto includes the endpoints, i.e. cell interfaces which will be used in the numerical flux function. Later on, we also require quadrature in the stochastic space. If the random variable is uniformly distributed, we can use a Gau\ss{}-Lobatto on $\cellR{\cellindR}$, \refone{$\cellindR=1,\ldots,\MEElements$,} with order $\SGtruncorder$, i.e., $\nbRnodes+1$  points and weights $(\quadRpoint{\xiQuadIndex},\quadRweight{\xiQuadIndex})$, $\xiQuadIndex=0,\ldots,\nbRnodes$, where $\nbRnodes = \ceil*{\frac{\SGtruncorder+1}{2}}$. 
For other distributions we use the corresponding Gauß quadrature based on the orthogonal basis polynomials and weighted by the conditional probability density function. We scale the quadrature weights such that
\begin{equation}\label{eq:quad}
\sum_{\xQuadIndex=0}^{\nbxnodes} \sum_{\xiQuadIndex=0}^{\nbRnodes}\quadxweight{\xQuadIndex}\quadRweight{\xiQuadIndex}=1, \qquad
\int_{\randomElement{\MEIndex}}\!\bg\refone{(\uncertainty)}\xiPDF\mathrm{d}\uncertainty\approx  \sum_{\xiQuadIndex=0}^{\nbRnodes}\bg(\quadRpoint{\xiQuadIndex})\quadRweight{\xiQuadIndex}.\end{equation}

%
The time discretization of the semi-discrete system \eqref{eq:weakformulation} is performed using a $\RKorder$-th order SSP Runge-Kutta method, see \cite{Schneider2015,Schneider2016}. In each Runge Kutta stage we apply the stochastic slope limiter \eqref{eq:slopeLimiter} to the gPC polynomial.  Afterwards we perform a polynomial reconstruction in order to derive the values of the numerical flux at the left and right limits of cell interfaces \eqref{eq:limitinterface}.

We reconstruct the solution vector $\solution(\timevar,\x,\uncertainty)$ as a polynomial of order (at most) $\RKorder$ for each quadrature node of $\uncertainty$ and in all spacial cells $\cell{\cellind}$, \refone{$i=1,\ldots,\ncells$}. This can be done by the CWENOZ scheme explained in \cite{Cravero2017}, which combines the CWENO method from \cite{Levy1999}, ensuring uniform accuracies within the cells, and WENOZ \cite{Don2013} for an optimal choice of the nonlinear WENO weights.
For all $\cellind=1,\ldots,\ncells$, $\xiQuadIndex=1,\ldots,\nbRnodes$, given the \refone{$x$-cell means at every $\uncertainty$ quadrature point  
$	\sum_{\SGsumIndex=0}^{\SGtruncorder} \solutioncm{\stageind}{\SGsumIndex,\cellind,\cellindR}\xiBasisPoly{\SGsumIndex,\cellindR}(\quadRpoint{\xiQuadIndex})$}, we represent the solution in each cell $\cell{\cellind}$ and for each $\quadRpoint{\xiQuadIndex}$ as a polynomial of degree $\WENOdegree := 2\WENOdegreeCoeff-2$ such that $\WENOdegree+1 \leq \RKorder$, more precisely
\begin{equation}\label{eq:WENOpoly} 
\weno_\Domain(\timevar, \x, \quadRpoint{\xiQuadIndex})\big|_{\cell{\cellind}} = \sum_{\kappa = 0}^{\WENOdegree} p^{(\kappa)}_{\cell{\cellind}}(\refone{\timevar},\quadRpoint{\xiQuadIndex})\,\refone{\varphi_{\kappa,i}(\x)}, 
\end{equation}
where $p^{(\kappa)}_{\cell{\cellind}}\refone{(\timevar,\quadRpoint{\xiQuadIndex})}$ are the coefficients obtained by the CWENOZ algorithm from \cite{Cravero2017} \refone{at time $\timevar$, in each cell $\cell{\cellind}$, $i=1,\ldots,\ncells$ and 
every quadrature point $\quadRpoint{\xiQuadIndex}$, $\xiQuadIndex=0,\ldots,\nbxiQuadNodes$ in $\cellR{\cellindR}$, $\cellindR=1,\ldots,\MEElements$. The basis polynomials on the physical cell $\cell{\cellind}$ are given by $(\varphi_\kappa)_{\kappa=0:\WENOdegree}$, we may set $$\varphi_{\kappa,i}(\x)=\frac{(\x-\x_{\cellind+1})^\kappa}{\kappa!},$$
This} polynomial will then be evaluated at the quadrature nodes $\x = \quadxpoint{0}$ and $\x = \quadxpoint{\nbxnodes}$ in $\cell{\cellind}$, \refone{$i=1,\ldots,\ncells$,} in order to obtain the limits \eqref{eq:limitinterface}.

\begin{remark}
	We can remap a solution vector \refone{$\solution\big(\timevar,\x,\quadRpoint{\xiQuadIndex}\big)$, $\xiQuadIndex=0,\ldots,\nbxiQuadNodes$} in each \ME{} $\randomElement{\MEIndex}$ to its gPC moments using  \begin{equation}\label{eq:SGremap}
	\refone{\solution_{\SGsumIndex,\cellindR}(\timevar,\x) } = \int_{\randomElement{\MEIndex}} \solution(\timevar,\refone{\x},\uncertainty) \xiBasisPoly{\SGsumIndex,j}(\uncertainty) \xiPDF \mathrm{d}\uncertainty\, \refone{\approx}  \sum_{\xiQuadIndex=0}^{\nbxiQuadNodes} \solution(\timevar,\refone{\x},\quadRpoint{\xiQuadIndex})\xiBasisPoly{\SGsumIndex,j}(\quadRpoint{\xiQuadIndex}) \quadRweight{\xiQuadIndex},
	\end{equation}
	for every $k=0,\ldots\SGtruncorder$ and $j=1,\ldots,\MEElements$.
\end{remark}

\subsection{Algorithm}
\refone{
The example from \secref{sec:LAex} demonstrated the significance of limiting techniques not only on the spatial but also on the stochastic grid. In this section we formulate an algorithm for the stochastic Galerkin scheme explained in the previous section.
} 

\refone{
The time discretization of the semi-discrete system \eqref{eq:weakformulation} is performed using a $\RKorder$-th order SSP Runge-Kutta method, see \cite{Schneider2016}. In each Runge-Kutta stage we apply the stochastic slope limiter \eqref{eq:slopeLimiter} to the gPC polynomial.  Afterwards we perform a polynomial reconstruction in order to derive the values of the numerical flux at the left and right limits of cell interfaces \eqref{eq:limitinterface}.
}
We summarize our results in the following algorithm of the \scheme{} scheme,
\refone{
where we denote by $(A,b)$, $A\in\R^{S\times S}$, $b\in\R^S$, the Butcher array from a $R$-th SSP Runge-Kutta scheme \cite{Gottlieb2003} with $S$ stages. We define the differential operator as the right hand side of \eqref{eq:LaxFriedrichs} to 
\begin{subequations}\label{eq:rhs}
 \begin{align}
 \rhs^{(n)}\big( \solutioncm{n}{\SGsumIndex,\cellind,\cellindR},\solution^{(n)}_{i,j}(\quadxpoint{0,\nbxnodes},\quadRpoint{0:\nbxiQuadNodes})\big) :=  \solutioncm{n}{\SGsumIndex,\cellind,\cellindR}
 &- \frac{\dtstepsize}{\Delta\x} \sum_{\xiQuadIndex=0}^{\nbxiQuadNodes} \hat{\flux}\big(\solution^{(n)}_{i,j}(\quadxpoint{\nbxnodes},\quadRpoint{\xiQuadIndex}),\solution^{(n)}_{i+1,j}(\quadxpoint{0},\quadRpoint{\xiQuadIndex})\big)\xiBasisPoly{\SGsumIndex,\cellindR}(\quadRpoint{\xiQuadIndex})\quadRweight{\xiQuadIndex} \\
 &+\frac{\dtstepsize}{\Delta\x} \sum_{\xiQuadIndex=0}^{\nbxiQuadNodes} \hat{\flux}\big(\solution^{(n)}_{i-1,j}(\quadxpoint{\nbxnodes},\quadRpoint{\xiQuadIndex}),\solution^{(n)}_{i,j}(\quadxpoint{0},\quadRpoint{\xiQuadIndex})\big) 
 \xiBasisPoly{\SGsumIndex,\cellindR}(\quadRpoint{\xiQuadIndex}) \quadRweight{\xiQuadIndex},\end{align}
\end{subequations}
for $k=0,\ldots,\SGtruncorder$ and where $\solution^{(n)}_{i,j}(\quadxpoint{\xQuadIndex},\quadRpoint{\xiQuadIndex})=\solution(t_n,\quadxpoint{\xQuadIndex},\quadRpoint{\xiQuadIndex})\big|_{\cell{i}\times\cellR{j}}$ denotes the solution $\solution$ at time $t_n$ in the cell $\cell{i}\times\cellR{j}$, $i=1,\ldots,\ncells$, $j=1,\ldots,\MEElements$ evaluated at the quadrature nodes $\quadxpoint{\xQuadIndex}$, $\xQuadIndex=0,\ldots,\nbxnodes$, and $\quadRpoint{\xiQuadIndex}$, $\xiQuadIndex=0,\ldots,\nbxiQuadNodes$.
}
\refone{
 \begin{algorithm}[H] 
	\caption{\refone{\scheme{} scheme (\schemeshort{})}}
	\label{algo:hSGWENO}
	\begin{algorithmic}[1]
	\State $\solutioncm{0}{}
	\leftarrow \text{vec}\Big(\sum_{\xQuadIndex=0}^{\nbxnodes}\sum_{\xiQuadIndex=0}^{\nbxiQuadNodes} \solution(0,\quadxpoint{\xQuadIndex},\quadRpoint{\xiQuadIndex})|_{\cell{\cellind}\times\cellR{\cellindR}}\xiBasisPoly{\SGsumIndex,j}(\quadRpoint{\xiQuadIndex})\quadxweight{\xQuadIndex} \quadRweight{\xiQuadIndex}\Big)_{k=0:\SGtruncorder, i=1:\ncells, j=1:\MEElements}$  \hfill \textit{\# initial state} \smallskip
	\For{$n=0$ to $N_t$} \hfill \textit{\# time loop} \smallskip
	\State Set $\solutioncmRK{0}{\SGsumIndex,\cellind,\cellindR} \leftarrow \solutioncm{\timeind}{}
	$ and $\rhs^{(1)}\leftarrow0$ \hfill \textit{\# initialization time step $\timeind$} \smallskip
	\For{$s=1$ to $S$} \hfill  $\#$ \textit{loop over RK stages} \smallskip
		\State $\solutioncmRK{\stageind}{\SGsumIndex,\cellind,\cellindR} \leftarrow \solutioncmRK{0}{\SGsumIndex,\cellind,\cellindR} + \Delta t_n \sum \limits_{\tilde{\stageind}=0}^{\stageind-1} A_{\tilde{\stageind},\stageind} \rhs^{(\tilde{\stageind})}$
		\hfill\textit{\# RK time update}
		\State $\solutioncmRK{\stageind}{\SGsumIndex,\cellind,\cellindR}\leftarrow \Lambda\Pi_\xi\big(\solutioncmRK{\stageind}{\SGsumIndex,\cellind,\cellindR}\big)$ \hfill \textit{\# call of slope limiter \eqref{eq:slopeLimiter}} \smallskip
	 \For{$\xiQuadIndex=0$ to $\nbxiQuadNodes$}\smallskip
		\State $\solutionRK{\stageind}{\cellind,\cellindR}(\x,\quadRpoint{\xiQuadIndex})\leftarrow
		\weno_\Domain\big(\solutioncmRK{\stageind}{\SGsumIndex,\cellind,\cellindR}
		\big)$  \hfill \textit{\# WENO reconstruction \eqref{eq:WENOpoly}} \smallskip
		\EndFor\smallskip
		\State $\rhs^{(\stageind+1)} \leftarrow \rhs\big(\solutioncmRK{\stageind}{\SGsumIndex,\cellind,\cellindR},\solutionRK{\stageind}{\cellind,\cellindR}(\quadxpoint{0,\nbxnodes},\quadRpoint{0:\nbxiQuadNodes})\big)$ \hfill \textit{\# update differential operator \eqref{eq:rhs}} \smallskip
	\EndFor 
		\State Set $\solutioncm{\timeind+1}{}
		 \leftarrow \solutioncmRK{0}{\SGsumIndex,\cellind,\cellindR} + \Delta t_n \sum\limits_{\stageind=1}^S b_\stageind\rhs^{(\stageind)}$ \hfill \textit{\# solution at new time step}
		\EndFor
\end{algorithmic}
\end{algorithm}
}
\begin{remark}\label{rem:hyplimiter}
	If we consider hyperbolic systems of equations, we may apply a hyperbolicity-preserving limiter in addition to the slope limiter after step 6 \& 8, \refone{since the stochastic slope limiter is not constructed to preserve admissible solutions within sG and the hyperbolicity limiter is not enough to satisfyingly dampen oscillations, as seen in \secref{sec:LAex}. We already mentioned this behavior in the beginning of \secref{sec:dg}}. For more details see \cite{Schlachter2017a}. \refone{The hyperbolicity limiter can be applied in the same manner for example to the filtered stochastic Galerkin scheme \cite{Kusch2018} such that the method can be used for any system of equations that might lose hyperbolicity. Note that the procedure of combining hyperbolicity- or positivity-preserving limiters with non-oscillatory limiters is performed in the theory of deterministic equations  for example in \cite{Schneider2015,Zhang2011}.}
\end{remark}

\begin{remark}
	\refthree{
	If the hyperbolicity limiter is applied in addition to the stochastic slope limiter in step 6, we can perform both limiters simultaneously by replacing \eqref{eq:slopeLimiter} with 
	$$\stochsllimiter{\solution}\big|_{\randomElement{\MEIndex}} 
	= \begin{cases}
	\begin{pmatrix}
	(\solution_{0,j})^T\\
	(1-\limitervariable)\text{m}\big(\solution_{1,j},\, \solution_{0,j+1}-\solution_{0,j}, \, \solution_{0,j}-\solution_{0,j-1}\big)\\[0.1cm]
	(0,\ldots,0)^T\\
	\vdots\\
	(0,\ldots,0)^T
	\end{pmatrix} & \text{if } TC_j(\solution) = 1, \\[1.5cm]
	\quad \solution\big|_{\randomElement{\MEIndex}} & \text{else},
	\end{cases}$$
	for $\MEIndex=1,\ldots,\MEElements$ and where the value $\limitervariable$ is derived as described in \cite{Schlachter2017a}.}
\end{remark}

\subsection{Full \fullWENO{} reconstruction of uncertain hyperbolic conservation laws}\label{sec:fullWENO}
In our numerical results, we compare \algoref{algo:hSGWENO} to a method using WENO reconstruction in both the physical and stochastic space. This is motivated by the Stochastic Finite Volume Method from \cite{Tokareva2014}, but now using a \fullWENO{} scheme instead of applying 1D reconstructions to the physical and stochastic space consecutively. 
\refone{In this method, we directly calculate the cell means $\overline{\solution}_{i,\MEIndex}$ based on \eqref{eq:FVscheme} instead of deriving the coefficients $\overline{\solution}_{k,i,\MEIndex}$ from \eqref{eq:LaxFriedrichs} as in the previously explained \schemeshort{} scheme.
In each Runge-Kutta stage, the numerical solution is then represented by a polynomial WENO reconstruction in $\x$ and $\uncertainty$.}
Given the $\x$-$\uncertainty$ cell means $\overline{\solution}_{i,\MEIndex}$, we obtain by the CWENOZ algorithm 
\begin{equation}\label{eq:WENO2Dpoly}
\weno_{\Domain\times\randomSpace}(\timevar, \x, \uncertainty)\big|_{\cell{\cellind}\times\randomElement{\MEIndex}}
= \sum_{\kappa = 0}^{\WENOdegree} \sum\limits_{\SGsumIndex=0}^{\SGtruncorder} p^{(\kappa,\SGsumIndex)}_{\cell{\cellind}\times\randomElement{\MEIndex}}\!\refone{(\timevar)} \varphi_{i,\kappa}(\x)\xiBasisPoly{j,\SGsumIndex}(\uncertainty)
\end{equation}
with basis polynomials $(\varphi_\kappa)_{\kappa=0:\WENOdegree}$ on the physical cells $\cell{\cellind}$, $i=1,\ldots,\ncells$, \refone{and basis polynomials $(\phi_k)_{k=0:\SGtruncorder}$ on the \ME{}s $\cellR{\cellindR}$, $j=1,\ldots,\MEElements$, given by \eqref{gPCApprox}.} We use Legendre polynomials for the basis functions w.r.t. $x$ in our numerical calculations in \secref{sec:results}.
Finally, we replace the evaluation of the flux at the interfaces in \eqref{eq:LaxFriedrichs} by the numerical Lax-Friedrichs flux \eqref{eq:globalLF}, deriving the left and right limits with help of the reconstructed polynomial \eqref{eq:WENO2Dpoly}.

\refone{The expansion of the solution into the WENO polynomial with respect to the uncertainty substitutes the application of the stochastic slope limiter \eqref{eq:slopeLimiter}. Similarly as if the troubled cell indicator would mark every cell as troubled, we use the WENO polynomial in order to reduce the oscillations that appear in the solution manifold. Therefore, the extension to a full WENO reconstruction within $\xi$ seems natural to be compared to \algoref{algo:hSGWENO}. To sum up, the first algorithm uses a WENO reconstruction only in $\x$ and applies the stochastic slope limiter in $\xi$, the second algorithm directly uses a two-dimensional WENO reconstruction in $\x$ and $\xi$. Both methods are supposed to reduce Gibbs oscillations but do not ensure hyperbolicity of the solution, thus}  
similar to \rmref{rem:hyplimiter}, we have to apply the hyperbolic slope limiter from \cite{Schlachter2017a} if we consider systems of conservation laws ($\dimension>1$).

The scheme is summarized in the following algorithm,
\refone{
where we now define the differential operator as the right hand side of \eqref{eq:FVscheme} to 
\begin{subequations}\label{eq:rhs2}
	\begin{align}
	\rhs^{(n)}\big( \solutioncm{n}{\cellind,\cellindR},\solution^{(n)}_{i,j}(\quadxpoint{0,\nbxnodes},\quadRpoint{0:\nbxiQuadNodes})\big) :=  \solutioncm{n}{\cellind,\cellindR}
	&- \frac{\dtstepsize}{\Delta\x} \sum_{\xiQuadIndex=0}^{\nbxiQuadNodes} \hat{\flux}\big(\solution^{(n)}_{i,j}(\quadxpoint{\nbxnodes},\quadRpoint{\xiQuadIndex}),\solution^{(n)}_{i+1,j}(\quadxpoint{0},\quadRpoint{\xiQuadIndex})\big)\quadRweight{\xiQuadIndex} \\
	&+\frac{\dtstepsize}{\Delta\x} \sum_{\xiQuadIndex=0}^{\nbxiQuadNodes} \hat{\flux}\big(\solution^{(n)}_{i-1,j}(\quadxpoint{\nbxnodes},\quadRpoint{\xiQuadIndex}),\solution^{(n)}_{i,j}(\quadxpoint{0},\quadRpoint{\xiQuadIndex})\big) 
	 \quadRweight{\xiQuadIndex},\end{align}
\end{subequations}
for $i=1,\ldots,\ncells$ and $j=1,\ldots,\MEElements$, thus we consider the $x-\xi$ cell means of $\solution$ instead of the gPC coefficients.
} 

\refone{
\begin{algorithm}[H] 
	\caption{\refone{Full \fullWENO{} reconstruction of uncertain hyperbolic conservation laws (\fullWENO{})}}
	\label{algo:fullWENO}
	\begin{algorithmic}[1]
	\State $\solutioncm{0}{}
		\leftarrow \text{vec}\Big(\sum_{\xQuadIndex=0}^{\nbxnodes}\sum_{\xiQuadIndex=0}^{\nbxiQuadNodes} \solution(0,\quadxpoint{\xQuadIndex},\quadRpoint{\xiQuadIndex})|_{\cell{\cellind}\times\cellR{\cellindR}}\quadxweight{\xQuadIndex} \quadRweight{\xiQuadIndex}\Big)_{i=1:\ncells, j=1:\MEElements}$  \hfill \textit{\# initial state} \smallskip
	\For{$n=0$ to $N_t$} \hfill \textit{\# time loop} \smallskip
	\State Set $\solutioncmRK{0}{\SGsumIndex,\cellind,\cellindR} \leftarrow \solutioncm{\timeind}{}
	$ and $\rhs^{(1)}\leftarrow0$ \hfill \textit{\# initialization time step $\timeind$} \smallskip
	\For{$s=1$ to $S$} \hfill  $\#$ \textit{loop over RK stages}
	\State $\solutioncmRK{\stageind}{\SGsumIndex,\cellind,\cellindR} \leftarrow \solutioncmRK{0}{\SGsumIndex,\cellind,\cellindR} + \Delta t_n \sum \limits_{\tilde{\stageind}=0}^{\stageind-1} A_{\tilde{\stageind},\stageind} \rhs^{(\tilde{\stageind})}$
	\hfill\textit{\# RK time update}
	\State $\solutionRK{\stageind}{\cellind,\cellindR}(\x,\uncertainty)
	\leftarrow\weno_{\Domain\times\randomSpace}\big(\solutioncmRK{\stageind}{\SGsumIndex,\cellind,\cellindR}\big)$  \hfill \textit{\# WENO reconstruction \eqref{eq:WENO2Dpoly}} \smallskip
	\State $\rhs^{(\stageind+1)} \leftarrow \rhs\big(\solutioncmRK{\stageind}{\SGsumIndex,\cellind,\cellindR},\solutionRK{\stageind}{\cellind,\cellindR}(\quadxpoint{0,\nbxnodes},\quadRpoint{0:\nbxiQuadNodes})\big)$ \hfill \textit{\# update differential operator \eqref{eq:rhs2}} \smallskip
	\EndFor \smallskip
	\State Set $\solutioncm{\timeind+1}{}
	\leftarrow \solutioncmRK{0}{\SGsumIndex,\cellind,\cellindR} + \Delta t_n \sum\limits_{\stageind=1}^S b_\stageind\rhs^{(\stageind)}$ \hfill \textit{\# solution at new time step}
	\EndFor
\end{algorithmic}
\end{algorithm}
}

\refone{The difference to \algoref{algo:hSGWENO} is that we perform the time stepping with the differential operator \eqref{eq:rhs2} based on the $\x$-$\xi$ cell means instead of the gPC coefficients as in \eqref{eq:rhs}. Moreover, in step 6, we use the two-dimensional WENO reconstruction \eqref{eq:WENO2Dpoly} with respect to $\x$ and $\xi$. In \algoref{algo:hSGWENO} this is substituted through steps 6-9 with the application of the stochastic slope limiter and the WENO reconstruction \eqref{eq:WENOpoly} with respect to $\x$ for every $\xi$ quadrature node.} 
\refthree{Hence, \algoref{algo:hSGWENO} uses a $\x$ WENO reconstruction and a $\uncertainty$ slope limiter, whereas \algoref{algo:fullWENO} is based on an $\x-\uncertainty$ WENO reconstruction. Both require the hyperbolicity limiter if systems of conservation laws are considered that can lose hyperbolicity.}
\section{Numerical Results}\label{sec:results}
\refone{
In this section we discuss extensive numerical experiments for the presented non-oscillation stochastic Galerkin methods that we described within \algoref{algo:hSGWENO} and \algoref{algo:fullWENO}. First of all, we present a convergence analysis for a refinement in the physical and stochastic domain. We then analyze the reduction of oscillations within the two algorithms by comparing it to classical stochastic Galerkin and a reference solution in terms of the $L_1$ error and the total variation.
\subsection{Convergence Tests}
For the convergence analysis, we illustrate two exemplary test cases for the Burgers' and Euler equations and we will observe that they verify the expected order of convergence of our numerical methods.
\subsubsection{Burgers' Equation}
We consider the Burgers' equation given by
\begin{equation}\label{eq:burgers}
 \frac{\partial}{\partial t} \solOned +  \frac{\partial}{\partial x} \Big(\frac{\solOned^2}{2}\Big)=0.
\end{equation}
We analyze the convergence of \algoref{algo:hSGWENO} and \algoref{algo:fullWENO} for this equation and a refinement within the physical space $\Domain$. In the following, we consider linear stochastic Galerkin polynomials with degree $\SGtruncorder = \reffour{2}$ and one random element $\MEElements = 1$. Moreover, we let $\uncertainty \sim \mathcal{U}(-1,1)$ and define the following stochastic modes
\begin{equation}\label{eq:BurgerExact}
\widehat{\solution}_0(\timevar,\x) = \frac{\x}{\timevar -1}, \qquad \widehat{\solution}_1(\timevar,\x) = \frac{\reffour{C_1}}{2\timevar-2},\qquad \reffour{\widehat{\solution}_2(\timevar,\x) = \frac{C_2}{2\timevar-2}}
\end{equation}
in $\Domain=[0,1]$, with non-zero \reffour{constants $C_1$, $C_2$} and $\timevar <2$. Hence, $\widehat{\solution}(\timevar,\x,\uncertainty)=\xiBasisPoly{0}(\uncertainty)\widehat{\solution}_0(\timevar,\x) + \xiBasisPoly{1}(\uncertainty)\widehat{\solution}_1(\timevar,\x)$ $\reffour{+\xiBasisPoly{2}(\uncertainty)\widehat{\solution}_2(\timevar,\x)}$ is an exact solution of the stochastic Galerkin system (note that \reffour{$\xiBasisPoly{0}$, $\xiBasisPoly{1}$ and $\xiBasisPoly{2}$ are given in \eqref{eq:legendrepolys}}) with expectation $\widehat{\solution}_0$ and variance $\widehat{\solution}_1^2$. Therefore, we derive the approximation of this analytical solution with our \schemeshort{} and \fullWENO{} schemes at $t=0.2$ with \reffour{$C_1=C_2=1$} and determine the error of the expectation \eqref{eq:exp} and variance \eqref{eq:var} in the $L_1$-norm. 
We denote 
$$
L_1\text{-Mean} = \int_\Domain \big|\mathbb{E}(\solution)(\x,\timevar) - \mathbb{E}(\widehat{\solution})(\x,\timevar)\big|\,\mathrm{d}\x,\\
$$
where $t=0.2$ and the error for the variance analogously with $L_1$-Var.
 The results for $\WENOdegree = 0$ (i.e. no WENO reconstruction) and $\WENOdegree = 2$ are shown in \figref{fig:BurgerConv} and \tabref{table:BurgerConv}, where every error is converging with the expected order $\WENOdegree +1$ of the WENO reconstruction, underlining the conservativity of our numerical scheme. 
\externaltikz{BburgerConv}{\begin{figure}[htb]
	\centering
	
		\begin{tikzpicture}
		\begin{groupplot}[
		group style={group size=2 by 2, horizontal sep = 2cm,  vertical sep = 2cm},
			]	
		\nextgroupplot[
		width=\figurewidth,
		height=\figureheight,
		scale only axis,
		xmode=log,
		xmin=0.001,
		xmax=0.2,
		xminorticks=true,
		xlabel style={font=\color{white!15!black}},
		xlabel={$\Delta\x$},
		ymode=log,
		ymin=5e-6,
		ymax = 1e-1,
		yminorticks=false,
		ylabel style={font=\color{white!15!black}},
		ylabel={error},
		ymajorgrids=true,
		axis background/.style={fill=white},
		legend style={at={(0.57,0.)},
			anchor=south west,
			, legend cell align=left, align=left, draw=white!15!black},
		legend columns=1
		]
		\addplot+[color=tuklblue, line width=1.0pt, mark size=2.5pt, mark=o, mark options={solid, tuklblue}]
		table[row sep=crcr]{%
			0.125	0.011661572639229\\
			0.0625	0.006024023957485\\
			0.03125	0.003062525104731\\
			0.015625	 0.001544135260667\\
			0.0078125	0.000775282027038\\
			0.00390625	0.000388491694569 \\
			0.001953125	0.0001945064600731\\
		};
		\addlegendentry{$\Lp{1}\text{-Mean}$}
		\addplot+[color=tuklred, line width=1.0pt, mark size=3.0pt, mark=triangle, mark options={solid, tuklred}]
		table[row sep=crcr]{%
			0.125	0.005957372504392\\
			0.0625	 0.003043328737822\\
			0.03125	0.001538073978334\\
			0.015625	0.000773217325290\\
			0.0078125	0.000387706953286\\
			0.00390625	0.000194076004132\\
			0.001953125	0.000097034582034\\
		};
		\addlegendentry{$\Lp{1}\text{-Var}$}
		\addplot+[color=tuklblue, densely dashed, line width=1.0pt, mark size=2.5pt, mark=o, mark options={solid, tuklblue}]
		table[row sep=crcr]{
			0.125	0.015490029298076\\
			0.0625	0.007847548715125\\
			0.03125	 0.003947740331894\\
			0.015625	 0.001979737610614\\
			0.0078125	0.000991460751382\\
			0.00390625	0.000496092495726 \\
			0.001953125	0.000248139784265\\
		};
		\addlegendentry{$\Lp{1}\text{-Mean}^*$}
				\addplot+[color=tuklred, line width=1.0pt, densely dashed, mark size=3.0pt, mark=triangle, mark options={solid, tuklred}]
		table[row sep=crcr]{%
			0.125	0.021944934114267\\			
			0.0625	0.011116625894822\\
			0.03125	0.005592021927508 \\
			0.015625	0.002804263253131\\
			0.0078125	0.001404364995972\\
			0.00390625	0.000702691732832\\
			0.001953125	0.000351477272962\\
		};
		\addlegendentry{$\Lp{1}\text{-Var}^*$}
		\addplot+[mark = none, color=black, line width=0.8pt]
		table[row sep=crcr]{%
			0.125	5e-02\\
			0.001953125	0.00078125\\
		};
		\addlegendentry{$\landau(\Delta\x)$}
\nextgroupplot[
width=1.\figurewidth,
height=\figureheight,
scale only axis,
xmode=log,
xmin=0.001,
xmax=0.2,
xminorticks=true,
xlabel style={font=\color{white!15!black}},
xlabel={$\Delta\x$},
ymode=log,
ymin=1e-15,
ymax=5e-7,
yminorticks=false,
ylabel style={font=\color{white!15!black}},
ylabel={error},
ymajorgrids=true,
axis background/.style={fill=white},
legend style={at={(0.57,0.)},
                anchor=south west,
, legend cell align=left, align=left, draw=white!15!black},
legend columns=1
]
\addplot+[color=tuklblue, line width=1.0pt, mark size=2.5pt, mark=o, mark options={solid, tuklblue}]
  table[row sep=crcr]{%
0.125	0.616394907095030e-07\\
0.0625	0.077519627357872e-07\\
0.03125	0.009705404087823e-07\\
0.015625	0.001212547671048e-07\\
0.0078125	0.000151801347778e-07\\
0.00390625	0.000019086651832e-07\\
0.001953125	0.000003222103103e-07\\
};
\addlegendentry{$\Lp{1}\text{-Mean}$}
\addplot+[color=tuklred, line width=1.0pt, mark size=3.0pt, mark=triangle, mark options={solid, tuklred}]
table[row sep=crcr]{%
	0.125	0.770493629315139e-07\\
	0.0625	0.096899606359668e-07\\
	0.03125	0.012131892677958e-07\\
	0.015625	0.001515969572097e-07\\
	0.0078125	0.000190384192786e-07\\
	0.00390625	0.000025019522044e-07\\
	0.001953125	0.000005996284198e-07\\
};
\addlegendentry{$\Lp{1}\text{-Var}$}
\addplot+[color=tuklblue, densely dashed, line width=1.0pt, mark size=2.5pt, mark=o, mark options={solid, tuklblue}]
table[row sep=crcr]{%
	0.125	0.496642082788851e-07\\
	0.0625	0.063041818596996e-07\\
	0.03125	0.007879672741119e-07\\
	0.015625	0.000984456100081e-07\\
	0.0078125	0.000123058079704e-07\\
	0.00390625	0.000015337893478e-07\\
	0.001953125	0.000001840998540e-07\\
};
\addlegendentry{$\Lp{1}\text{-Mean}^*$}
\addplot+[color=tuklred, densely dashed, line width=1.0pt, mark size=3.0pt, mark=triangle, mark options={solid, tuklred}]
table[row sep=crcr]{%
	0.125	0.809427563380643e-07\\
	0.0625	0.102896576691886e-07\\
	0.03125	0.012865883944314e-07\\
	0.015625	0.001607559886379e-07\\
	0.0078125	0.000200929196996e-07\\
	0.00390625	0.000025001444345e-07\\
	0.001953125	0.000002965572392e-07\\
};
\addlegendentry{$\Lp{1}\text{-Var}^*$}
\addplot+[mark = none, color=black, line width=0.8pt]
  table[row sep=crcr]{%
0.125	3e-07\\
0.001953125	1.144409179687500e-12\\ 
};
\addlegendentry{$\landau(\Delta\x^3)$}
\end{groupplot}
\end{tikzpicture}%
\caption{Error plot of Burgers' equation with $\MEElements=1$, $\SGtruncorder=\reffour{2}$, $\WENOdegree=0$ (left) and $\WENOdegree=2$ (right). $L_1$-Mean and $L_1$-Var denote the errors of the mean and variance in \schemeshort{} and the dashed lines of  $L_1$-Mean$^*$ and $L_1$-Var$^*$   correspond to \fullWENO{}. Example \eqref{eq:BurgerExact}.}
\label{fig:BurgerConv}
\end{figure}}
\begin{table}[htbp]
	\centering
	\begin{tabular}{ccccccccc}
		\toprule
&\multicolumn{4}{c}{$\WENOdegree = 0$ - $\Lp{1}$-Mean}& \multicolumn{4}{c}{$\WENOdegree = 2$ - $\Lp{1}$-Mean}\\
\cmidrule(r){2-5} \cmidrule(r){6-9}
$\ncells$ & \schemeshort{} & eoc & \fullWENO{} & eoc & \schemeshort{} & eoc & \fullWENO{}& eoc \\
\midrule
8   & 1.1661e-03 & -- & 1.5490e-03 & -- &
6.1639e-08  & --       & 4.9664e-08 & --   \\
16  & 0.6024e-03 & 0.95 & 0.7847e-03 & 0.98 &
0.7751e-08  & 2.9     & 0.6304e-08 & 2.9  \\
32  & 0.3062e-03 & 0.98 & 0.3948e-03 & 0.99 &
0.0970e-08  & 2.9     & 0.0788e-08 & 3.0  \\ 
64  & 0.1544e-03 & 0.99 & 0.1980e-03 & 0.99 &
0.0121e-08  & 3.0   & 0.0098e-08 & 3.0  \\
128 &  0.0775e-03 & 0.99 & 0.0991e-03 & 0.99 &
0.0015e-08  & 3.0      & 0.0012e-08 & 3.0  \\
256 & 0.0388e-03 & 0.99 & 0.0496e-03 & 0.99 &
0.0002e-08  & 3.0  & 0.0001e-08 & 3.0   \\
512 & 0.0194e-03 & 0.99 & 0.0248e-03 & 1.0 &
0.0003e-09  & 2.6     & 0.0002e-09 & 3.0 \\
		\bottomrule
	\end{tabular}
	\caption{$\Lp{1}$ errors and experimental order of convergence (eoc) of \schemeshort{} and \fullWENO{} for the Burgers' equation with $\SGtruncorder=\reffour{2}$ and $\MEElements=1$. Example \eqref{eq:BurgerExact}.}
	\label{table:BurgerConv}
\end{table}
\subsubsection{Compressible Euler Equations}\label{sec:compreuler}
The one-dimensional compressible Euler equations for the flow of an ideal gas are given by 
\begin{subequations}\label{eulereq}
\begin{align}
 \frac{\partial}{\partial t} \density &+  \frac{\partial}{\partial x} \momentum &&= 0 ,\\
 \frac{\partial}{\partial t} \momentum &+  \frac{\partial}{\partial x} \left(\frac{\momentum^2}{\density} + \pressure\right)  &&=0,\\
 \frac{\partial}{\partial t} \energy &+  \frac{\partial}{\partial x} \left( (\energy + \pressure) \,\frac{\momentum}{\density}\right) &&  = 0,
\end{align}
\end{subequations}
with pressure $$\pressure = (\eulergamma-1)\left(\energy - \frac{1}{2}\frac{\momentum^2}{\density}\right)$$
and adiabatic constant $\eulergamma>1$. }

\refone{
We apply \algoref{algo:hSGWENO} and \algoref{algo:fullWENO} to the compressible Euler equations with $\eulergamma = 1.4$ and $\x\in\Domain=[0,2]$. In order to obtain an analytical solution, we use the method of manufactured solutions and introduce an additional source term $S(\x,\timevar,\uncertainty)$. We choose the following analytical solution
\begin{align} \label{eq:eulerSmoothSolution}
\widehat{\solution}(\timevar,\x,\uncertainty) = 
\begin{pmatrix}
\widehat{\density}(\timevar,\x,\uncertainty) \\
\widehat{\momentum}(\timevar,\x,\uncertainty)  \\
\widehat{\energy}(\timevar,\x,\uncertainty)  
\end{pmatrix}=
\begin{pmatrix}
1+ 0.1\cos(\pi(\x-\uncertainty \timevar)) \\[0.1cm]
\Big(1+ 0.1\cos(\pi(\x-\uncertainty \timevar)) \Big)\Big(1+0.1\sin(\pi(\x-\uncertainty \timevar))\Big) \\
\Big(1+ 0.1\cos(\pi(\x-\uncertainty \timevar)) \Big)^2
\end{pmatrix}.
\end{align}
The source term is then given by inserting \eqref{eq:eulerSmoothSolution} into \eqref{eulereq}.
We use periodic boundary conditions, $\uncertainty\sim\mathcal{U}(-1,1)$ and set $t=0.5$. Moreover, we consider $\WENOdegree = 2$ and divide the spatial domain $\Domain$ into $\ncells=1000$ cells. We then compute the $L_1$ error of our numerical scheme to the analytical solution \eqref{eq:eulerSmoothSolution} for the stochastic Galerkin degrees $\SGtruncorder=0,1$ and $\MEElements=1,\ldots,10$ \ME{}s. The error of the mean and variance from the density $\density$ is shown in \figref{fig:EulerConv} and \tabref{table:EulerConv}. We observe that the rates of convergence agree with the theoretical rates which are given by $\landau(\MEElements^{-2(\SGtruncorder+1)})$ according to \cite{Wan2006}. For $\SGtruncorder=1$ and a large number of \ME{}s $\MEElements\geq8$, the error introduced by the spatial discretization dominates the error in the stochastic discretization, so no further convergence can be obtained. Both algorithms yield comparable results.
\externaltikz{EulerConv2}{\definecolor{mycolor1}{rgb}{0.91000,0.41000,0.17000}%
\begin{figure}[htb]
	\centering
\begin{tikzpicture}
\begin{groupplot}[
group style={group size=2 by 2, horizontal sep = 2cm,  vertical sep = 2cm},
]	
\nextgroupplot[
width=\figurewidth,
height=\figureheight,
scale only axis,
xmode=log,
xmin=0,
xmax=12,
xminorticks=true,
xtick={2,4,6,8,10},
xticklabels={2,4,6,8,10},
xlabel style={font=\color{white!15!black}},
xlabel={$\MEElements$},
ymode=log,
ymin=1e-7,
ymax=4e-2,
yminorticks=false,
ylabel style={font=\color{white!15!black}},
ylabel={error},
ymajorgrids=true,
axis background/.style={fill=white},
legend style={at={(0.,0.)},
	anchor=south west,
	, legend cell align=left, align=left, draw=white!15!black},
legend columns=1
]
\addplot+[color=tuklblue, line width=1.0pt, mark size=2.5pt, mark=o, mark options={solid, tuklblue}]
table[row sep=crcr]{%
1	0.003284301899236\\
2	0.000595902129711\\
3	0.000254706099314\\
4	0.000141439679910\\
5	0.000090000964795 \\
6	0.000062312159796\\
7	0.000045701464198\\
8	0.000034954304176\\
9	0.000027601380853 \\
10	0.000022349480894\\
};
\addlegendentry{$\Lp{1}\text{-Mean}$}
\addplot+[color=tuklred, line width=1.0pt, mark size=3.0pt, mark=triangle, mark options={solid, tuklred}]
table[row sep=crcr]{%
1	0.005947152654306\\
2	0.001892437867540\\
3	0.000880193957933\\
4	0.000503093374615\\
5	0.000324392137213 \\
6	0.000226200682734\\
7	0.000166610712137\\
8	0.000127779416497\\
9	0.000101086062991 \\
10	0.000081957094311\\
};
\addlegendentry{$\Lp{1}\text{-Var}$}
\addplot+[color=tuklblue, densely dashed, line width=1.0pt, mark size=2.5pt, mark=o, mark options={solid, tuklblue}]
table[row sep=crcr]{
1	0.001248597999166\\
2	0.000326885663218\\
3   0.000147987712388\\
4	 0.000084751028845\\
5	0.000054371571508\\
6	0.000037892495730  \\
7	0.000027913978426\\
8	0.000021413978426\\
9	0.000016931571508\\
10	0.000013721545270\\
};
\addlegendentry{$\Lp{1}\text{-Mean}^*$}
\addplot+[color=tuklred, line width=1.0pt, densely dashed, mark size=3.0pt, mark=triangle, mark options={solid, tuklred}]
table[row sep=crcr]{%
1	0.004273118694632\\
2	0.001130202047710\\
3   0.000510638395431\\
4	  0.000291922006511\\
5	0.000187612651680\\
6	0.000130532710135 \\
7	0.000095902182906\\
8	0.000073523036197\\
9	0.000058135352815\\
10	0.000047116998932\\
};
\addlegendentry{$\Lp{1}\text{-Var}^*$}
\addplot [mark=none,thick,color=black, line width=0.8pt]
table[row sep=crcr]{%
	1	5e-04 \\
	10 5e-06\\
};
\addlegendentry{$\landau(\MEElements^{-2})$}
\nextgroupplot[
width=1.\figurewidth,
height=\figureheight,
scale only axis,
xmode=log,
xmin=0,
xmax=12,
xminorticks=true,
xtick={2,4,6,8,10},
xticklabels={2,4,6,8,10},
xlabel style={font=\color{white!15!black}},
xlabel={$\MEElements$},
ymode=log,
ymin=1e-11,
ymax=1e-2,
yminorticks=false,
ylabel style={font=\color{white!15!black}},
ylabel={error},
ymajorgrids=true,
axis background/.style={fill=white},
legend style={at={(0.,0.)},
	anchor=south west,
	, legend cell align=left, align=left, draw=white!15!black},
legend columns=1
]
\addplot+[color=tuklblue, line width=1.0pt, mark size=2.5pt, mark=o, mark options={solid, tuklblue}]
table[row sep=crcr]{%
1	0.000233960442042\\
2	0.000007014218721\\
3	0.000001249843933\\
4	0.000000376389928\\
5	0.000000148193827\\
6	0.000000073486778 \\
7	0.000000038881785\\
8   0.000000022353195\\
9	0.000000013834966\\
10	0.000000009310806\\
};
\addlegendentry{$\Lp{1}\text{-Mean}$}
\addplot+[color=tuklred, line width=1.0pt, mark size=3.0pt, mark=triangle, mark options={solid, tuklred}]
table[row sep=crcr]{%
1	0.0010277426714000\\
2	0.000078231738275\\
3	0.000016083640664\\
4	0.000005172231493\\
5	0.000002143524719\\
6	0.000001047423623\\
7	0.000000575574155\\
8   0.000000345968641\\
9	0.000000223564421\\
10	0.000000153533478\\
};
\addlegendentry{$\Lp{1}\text{-Var}$}
\addplot+[color=tuklblue, densely dashed, line width=1.0pt, mark size=2.5pt, mark=o, mark options={solid, tuklblue}]
table[row sep=crcr]{%
	1	0.000097732184902\\
	2	0.000006251699893\\
	3	0.000001258211366\\
	4	0.000000401856921\\
	5	0.000000164845635\\
	6	0.000000079518652\\
	7	0.000000042937870\\
	8	0.000000025230374\\
	9	0.000000016319690\\
	10	0.000000011523787 \\
};
\addlegendentry{$\Lp{1}\text{-Mean}^*$}
\addplot+[color=tuklred, densely dashed, line width=1.0pt, mark size=3.0pt, mark=triangle, mark options={solid, tuklred}]
table[row sep=crcr]{%
1	0.213457983584576e-03\\
2	0.015345798358458e-03\\
3	0.002910846560000e-03\\
4	0.000943278771608e-03\\
5	0.000393243541026e-03\\
6	0.000191283440913e-03\\
7	0.000103800911135e-03\\
8	0.000061120952176e-03\\
9	0.000039586453788e-03\\
10	0.000027675341206e-03\\
};
\addlegendentry{$\Lp{1}\text{-Var}^*$}
\addplot [mark = none,thick,color=black, line width=1.0pt]
table[row sep=crcr]{%
	1	2e-05\\
	10	2e-09\\
};
\addlegendentry{$\landau(\MEElements^{-4})$}
\end{groupplot}
\end{tikzpicture}%
\caption{Error plot of Euler equations (density) with manufactured source term, where $\ncells=1000$, $\WENOdegree=2$ and $\SGtruncorder=0$ (left), $\SGtruncorder=1$ (right). The solid lines correspond to \schemeshort{} and the dashed lines to \fullWENO{}. Example \eqref{eq:eulerSmoothSolution}.}
\label{fig:EulerConv}
\end{figure}}
\begin{table}[htbp]
	\centering
	\begin{tabular}{ccccccccc}
		\toprule
		\multicolumn{1}{c}{}&\multicolumn{4}{c}{$\SGtruncorder = 0$ - $\Lp{1}$-Mean}& \multicolumn{4}{c}{$\SGtruncorder = 1$ - $\Lp{1}$-Mean} \\
\cmidrule(r){2-5} \cmidrule(r){6-9}
$\MEElements$ & \schemeshort {} & eoc & \fullWENO{} & eoc & \schemeshort{} & eoc & \fullWENO{} & eoc \\
\toprule
1   & 3.2844e-03 & -- &  1.2486e-03 & -- & 
2.33996e-04  & --    & 0.9773e-04 & --      \\
2  & 0.5959e-03 & 2.46 & 0.0625e-04  & 1.93 &
0.0703e-04  & 5.06 & 0.0782e-03 & 3.96     \\
3  & 0.2547e-03 & 2.10 & 0.0126e-04 & 1.95 &
0.0126e-04  & 4.24 & 0.0161e-03 & 3.95     \\ 
4  & 0.1414e-03 & 2.05 & 0.0040e-04 & 1.93 &
0.0039e-04  & 4.11 & 0.0052e-03 & 3.96    \\
5 &  0.0899e-03 & 2.03 & 0.0016e-04 & 1.98 &
0.0015e-04  & 4.09 & 0.0021e-03 &  4.00    \\
6 & 0.0623e-03 & 2.02 &  0.7951e-07& 1.98 &
0.7349e-07 & 4.10 & 1.0474e-06&  4.00 \\
7 & 0.0457e-03 & 2.01 &  0.4294e-07 &  1.98 &
0.3888e-07  & 4.13  & 0.5756e-06 &4.00 \\
8 & 0.0349e-03& 2.01 & 0.2523e-07 & 1.98 & 
0.2235e-07& 4.14 & 0.3460e-06& 3.98\\
9 & 0.0276e-03& 2.00 &0.1632e-07 & 2.00 & 
0.1383e-07& 4.07 & 0.2236e-06 & 3.69\\
10 & 0.0223e-03& 2.00 & 0.1152e-07 & 2.00 & 
0.0931e-07& 3.76 & 0.1535e-06 & 3.30\\
		\bottomrule
	\end{tabular}
	\caption{$\Lp{1}$ errors and experimental order of convergence (eoc) for the Euler equations (density) with $\ncells=1000$ and $\WENOdegree=2$. Example \eqref{eq:eulerSmoothSolution}.}
	\label{table:EulerConv}
\end{table}
}

\subsection{Numerical Results: Analysis of Oscillation Reduction}
\refone{We compare the \schemeshort{} method from \algoref{algo:hSGWENO}, the 2D WENO reconstruction from \algoref{algo:fullWENO} as well as the standard stochastic Galerkin scheme. Moreover, we derive the $L_1$ error to a reference solution and the total variation for multiple test cases in order to analyze whether the methods are able to reduce the Gibbs oscillations detected in \secref{sec:LAex}.}

\subsubsection{Linear Advection}
We apply the \scheme{} scheme from Algorithm \ref{algo:hSGWENO} to the linear advection problem
$$ \frac{\partial}{\partial t}\solOned + a(\uncertainty)  \frac{\partial}{\partial x}\solOned = 0$$
with uncertain wave speed $a(\xi) = 1.5 + 0.5\, \xi$, $\uncertainty\sim\mathcal{U}(-1,1)$ as in \secref{sec:LAex} and where $\x\in\Domain=[0.4,\,2].$

At first, we again consider the initial conditions \eqref{sweInitialConditions} at $t=0.5$ and compare the \schemeshort{} scheme to a stochastic Galerkin scheme with WENO reconstruction in $\x$ (cf. \eqref{eq:WENOpoly}) but without the stochastic slope limiter. We refer to this method as the standard stochastic Galerkin (sG) scheme in the following. The results are shown for 3 \ME{}s in \figref{fig:LA_hSG} and for 10 \ME{}s in \figref{fig:LA_hSG_refined}. \refone{For an explanation of the notation used within these figures we refer to \secref{sec:LAex} and \figref{fig:gPCRefined}.} We observe the Gibbs phenomenon at the boundaries of the stochastic domain, whereas a refinement to 10 \ME{}s reduced the width but not the height of the overshoots. This validates our theoretical results from \secref{sec:LAex}. 
The oscillations are completely eliminated through the application of the stochastic slope limiter. 

\externaltikz{LA_hSG}{\def\xmin{0.8}
\def\xmax{1.7}
\def\ymin{-0.4}
\def\ymax{1.4}
\def\width{\textwidth}
\def\height{1.25\textwidth}
\def\xtick{1.0, 1.5}

\begin{figure}[htb]
	\begin{subfigure}[c]{0.19\textwidth}
		\begin{tikzpicture}
		\begin{axis}[width=\width, height=\height, xlabel=$x$, ylabel=$u$, ylabel style = {rotate=-90}, xmin=\xmin, xmax=\xmax, ymin=\ymin, ymax=\ymax, xtick={\xtick}, xticklabels={{$1.0$},{$1.5$}}] 
		\addplot[mark=none, solid, tuklred, line width=1pt] file {Images/LA_hSG/GPClimiter_xim09.txt};
		\addplot[mark=none, densely dashed, green, line width=1pt] file {Images/LA_hSG/GPClimiter_xim09_SL.txt};
		\end{axis}
		\end{tikzpicture}
		\subcaption{\refthree{$3\xiv =  -0.9$}}
	\end{subfigure}
	\begin{subfigure}[c]{0.19\textwidth}
		\begin{tikzpicture}
		\begin{axis}[width=\width, height=\height, xlabel=$x$, ylabel=$u$, ylabel style = {rotate=-90}, xmin=\xmin, xmax=\xmax, ymin=\ymin, ymax=\ymax, xtick={\xtick}, xticklabels={{$1.0$},{$1.5$}}] 
		\addplot[mark=none, solid, tuklred, line width=1pt] file {Images/LA_hSG/GPClimiter_xim04.txt};	
		\addplot[mark=none, densely dashed, green, line width=1pt] file {Images/LA_hSG/GPClimiter_xim04_SL.txt};
		\end{axis}
		\end{tikzpicture}
		\subcaption{\refthree{$3\xiv = -0.4$}}
	\end{subfigure}
	\begin{subfigure}[c]{0.19\textwidth}
		\begin{tikzpicture}
		\begin{axis}[width=\width, height=\height, xlabel=$x$, ylabel=$u$, ylabel style = {rotate=-90}, xmin=\xmin, xmax=\xmax, ymin=\ymin, ymax=\ymax, xtick={\xtick}, xticklabels={{$1.0$},{$1.5$}}] 
		\addplot[mark=none, solid, tuklred, line width=1pt] file {Images/LA_hSG/GPClimiter_xi0.txt};		
		\addplot[mark=none, densely dashed, green, line width=1pt] file {Images/LA_hSG/GPClimiter_xi0_SL.txt};
		\end{axis}
		\end{tikzpicture}
		\subcaption{\refthree{$3\xiv =  0.0$}}
	\end{subfigure}
	\begin{subfigure}[c]{0.19\textwidth}
		\begin{tikzpicture}
		\begin{axis}[width=\width, height=\height, xlabel=$x$, ylabel=$u$, ylabel style = {rotate=-90}, xmin=\xmin, xmax=\xmax, ymin=\ymin, ymax=\ymax, xtick={\xtick}, xticklabels={{$1.0$},{$1.5$}}] 
		\addplot[mark=none, solid, tuklred, line width=1pt] file {Images/LA_hSG/GPClimiter_xi04.txt};		
		\addplot[mark=none, densely dashed, green, line width=1pt] file {Images/LA_hSG/GPClimiter_xi04_SL.txt};
		\end{axis}
		\end{tikzpicture}
		\subcaption{\refthree{$3\xiv =0.4$}}
	\end{subfigure}
	\begin{subfigure}[c]{0.19\textwidth}
		\begin{tikzpicture}
		\begin{axis}[width=\width, height=\height, xlabel=$x$, ylabel=$u$, ylabel style = {rotate=-90}, xmin=\xmin, xmax=\xmax, ymin=\ymin, ymax=\ymax, xtick={\xtick}, xticklabels={{$1.0$},{$1.5$}}] 
		\addplot[mark=none, solid, tuklred, line width=1pt] file {Images/LA_hSG/GPClimiter_xi09.txt};		
		\addplot[mark=none, densely dashed, green, line width=1pt] file {Images/LA_hSG/GPClimiter_xi09_SL.txt};
		\end{axis}
		\end{tikzpicture}
		\subcaption{\refthree{$3\xiv =0.9$}}
	\end{subfigure}
	
	\caption{Comparison of the \schemeshort{} approximation (dashed) with 3 \ME{}s for the linear advection problem and the sG approach (solid) for the given $\xiv$s at $t=0.5$, 2000 space cells, $\SGtruncorder=2$, $\WENOdegree=2$. Example \eqref{sweInitialConditions}.}
	\label{fig:LA_hSG}
\end{figure}}

\externaltikz{LA_hSG_refined}{\def\xmin{0.8}
\def\xmax{1.7}
\def\ymin{-0.4}
\def\ymax{1.4}
\def\width{\textwidth}
\def\height{1.25\textwidth}
\def\xtick{1.0, 1.5}

\begin{figure}[htb]
	\begin{subfigure}[c]{0.19\textwidth}
		\begin{tikzpicture}
		\begin{axis}[width=\width, height=\height, xlabel=$x$, ylabel=$u$, ylabel style = {rotate=-90}, xmin=\xmin, xmax=\xmax, ymin=\ymin, ymax=\ymax, xtick={\xtick}, xticklabels={{$1.0$},{$1.5$}}] 
		\addplot[mark=none, solid, tuklred, line width=1pt] file {Images/LA_hSG_refined/GPClimiter_xim09.txt};
		\addplot[mark=none, densely dashed, green, line width=1pt] file {Images/LA_hSG_refined/GPClimiter_xim09_SL.txt};
		\end{axis}
		\end{tikzpicture}
		\subcaption{\refthree{$10\xiv = -0.9$}}
	\end{subfigure}
	\begin{subfigure}[c]{0.19\textwidth}
		\begin{tikzpicture}
		\begin{axis}[width=\width, height=\height, xlabel=$x$, ylabel=$u$, ylabel style = {rotate=-90},, xmin=\xmin, xmax=\xmax, ymin=\ymin, ymax=\ymax, xtick={\xtick}, xticklabels={{$1.0$},{$1.5$}}] 
		\addplot[mark=none, solid, tuklred, line width=1pt] file {Images/LA_hSG_refined/GPClimiter_xim04.txt};	
		\addplot[mark=none, densely dashed, green, line width=1pt] file {Images/LA_hSG_refined/GPClimiter_xim04_SL.txt};
		\end{axis}
		\end{tikzpicture}
		\subcaption{\refthree{$10\xiv = -0.4$}}
	\end{subfigure}
	\begin{subfigure}[c]{0.19\textwidth}
		\begin{tikzpicture}
		\begin{axis}[width=\width, height=\height, xlabel=$x$, ylabel=$u$, ylabel style = {rotate=-90}, xmin=\xmin, xmax=\xmax, ymin=\ymin, ymax=\ymax, xtick={\xtick}, xticklabels={{$1.0$},{$1.5$}}] 
		\addplot[mark=none, solid, tuklred, line width=1pt] file {Images/LA_hSG_refined/GPClimiter_xi0.txt};		
		\addplot[mark=none, densely dashed, green, line width=1pt] file {Images/LA_hSG_refined/GPClimiter_xi0_SL.txt};
		\end{axis}
		\end{tikzpicture}
		\subcaption{\refthree{$10\xiv = 0.0$}}
	\end{subfigure}
	\begin{subfigure}[c]{0.19\textwidth}
		\begin{tikzpicture}
		\begin{axis}[width=\width, height=\height, xlabel=$x$, ylabel=$u$, ylabel style = {rotate=-90}, xmin=\xmin, xmax=\xmax, ymin=\ymin, ymax=\ymax, xtick={\xtick}, xticklabels={{$1.0$},{$1.5$}}] 
		\addplot[mark=none, solid, tuklred, line width=1pt] file {Images/LA_hSG_refined/GPClimiter_xi04.txt};		
		\addplot[mark=none, densely dashed, green, line width=1pt] file {Images/LA_hSG_refined/GPClimiter_xi04_SL.txt};
		\end{axis}
		\end{tikzpicture}
		\subcaption{\refthree{$10\xiv = 0.4$}}
	\end{subfigure}
	\begin{subfigure}[c]{0.19\textwidth}
		\begin{tikzpicture}
		\begin{axis}[width=\width, height=\height, xlabel=$x$, ylabel=$u$, ylabel style = {rotate=-90}, xmin=\xmin, xmax=\xmax, ymin=\ymin, ymax=\ymax, xtick={\xtick}, xticklabels={{$1.0$},{$1.5$}}] 
		\addplot[mark=none, solid, tuklred, line width=1pt] file {Images/LA_hSG_refined/GPClimiter_xi09.txt};		
		\addplot[mark=none, densely dashed, green, line width=1pt] file {Images/LA_hSG_refined/GPClimiter_xi09_SL.txt};
		\end{axis}
		\end{tikzpicture}
		\subcaption{\refthree{$10\xiv = 0.9$}}
	\end{subfigure}
	
	\caption{Comparison of the \schemeshort{} approximation (dashed) with 10 \ME{}s for the linear advection problem and the sG approach (solid) for the given $\xiv$s at $t=0.5$, 2000 space cells, $\SGtruncorder=2$, $\WENOdegree=2$. Example \eqref{sweInitialConditions}.}
	\label{fig:LA_hSG_refined}
\end{figure}}

We compare the different approaches by calculating the $L_1$ error to the analytical solution and the total variation of the numerical solution over the spatial domain $\Domain$ and one of the \ME{}s $\randomElement{\MEIndex}$ (due to the uniform decomposition of the random space). The total variations with respect to $\x$ and $\uncertainty$ are given by
\begin{align}\label{eq:TVx}
TV_\x(\solution) &= \int_{\widetilde{\randomElement{}}}\sum_{k=0}^{\nbxnodes{\cdot}\ncells}\Big| \solution(\timevar,\tilde{x}_{k},\uncertainty)-\solution(\timevar,\tilde{x}_{k-1},\uncertainty)\Big| \xiPDFdxi,\\
\label{eq:TVxi}
TV_\uncertainty(\solution) &=  \int_{\Domain} \sum_{\xiQuadIndex=0}^{\nbRnodes} \Big| \solution(\timevar,\x,\quadRpoint{\xiQuadIndex})-\solution(\timevar,\x,\quadRpoint{\xiQuadIndex-1})\Big| \mathrm{d}\x,
\end{align}
where $\quadxpoint{k}$ with $k=0,\ldots,\nbxnodes{\cdot}\ncells$ enumerates all quadrature points in $\x$ over each of the cells $\cell{1},\ldots,\cell{\ncells}$, \refone{namely $\tilde{x}_{0}= \quadxpoint{0}\big|_{\cell{1}}$, $\tilde{x}_{1}= \quadxpoint{1}\big|_{\cell{1}}, \ldots, \tilde{x}_{\nbxnodes+1}= \quadxpoint{0}\big|_{\cell{2}}, \ldots,\tilde{x}_{\nbxnodes{\cdot}\ncells}=\quadxpoint{\nbxnodes}\big|_{\cell{\ncells}}$. In addition to that, we consider the total variation within our reference \ME{} $\widetilde{\randomElement{}}=\frac{1}{\MEElements}[-1,\,1]$ and $\quadRpoint{\xiQuadIndex}$, $\xiQuadIndex=0,\ldots,\nbxiQuadNodes$ denotes the corresponding quadrature nodes within this interval. 
} 

The results are given in  \tabref{table:LAerror}. We choose a fine discretization of 2000 space cells in the physical domain and consider the convergence within the stochastic space. We additionally compare the \schemeshort{} scheme to the \fullWENO{} method described in \secref{sec:fullWENO} and \algoref{algo:fullWENO}. Since the stochastic Galerkin solution represents the best approximating polynomial due to the underlying Galerkin projection, we do not expect to improve the $L_2$ or $L_1$ error. However, the oscillations and therefore the total variations \eqref{eq:TVx} and \eqref{eq:TVxi} should be minimized by the application of our slope limiter. 

Indeed, we observe that the $L_1$ errors of the new methods are slightly higher as for standard stochastic Galerkin while the total variation and hence the oscillation are reduced. The total variation with respect to the uncertainty $TV_\uncertainty$ is even smaller in \schemeshort{} and \fullWENO{} as for the analytical solution. This situation arises due to the lack of information about $\uncertainty$ that is transported within the numerical schemes \refone{since there is no numerical flux with respect to $\uncertainty$. This can be observed in \figref{fig:LAsurf}, where the  approximations in \figref{fig:surfLAc} and \figref{fig:surfLAd} reveal only a small impact of the different values for $\uncertainty$ compared to analytic solution in \figref{fig:surfLAa}.} Especially the \schemeshort{} method was able to improve the total variation wrt. $\x$ tremendously, more precisely, it only increased the total variation of the analytic solution for 3 \ME{}s by 0.4\% while standard Galerkin yield an increase of almost 24\%.  In this test case, the \fullWENO{} reconstruction performed not as good as the \schemeshort{} scheme and sometimes even has a larger total variation in $\x$ than standard stochastic Galerkin. In the upcoming examples, the two schemes will show comparable results. Note that we have used $\nbRnodes=1000$ quadrature nodes in $\uncertainty$  and $\nbxnodes=4$ quadrature nodes in $\x$ for the calculations of the total variation.

\begin{table}[htbp]
\centering
\begin{tabular*}{0.47\textwidth}{C{2.2cm}ccc}
	\toprule
	\begin{tabular}{@{}c@{}} \mbox{}\\\mbox{}\end{tabular} $L_1$ error & $\MEElements =3$ & $\MEElements=10$ & $\MEElements =30$\\[0.2cm]
	\midrule
	analytic&--&--&--\\[0.1cm]
	\midrule 
	sG&  0.0265&  0.0076 & 0.0031 \\[0.1cm]
	\midrule
	\schemeshort{}&  0.0329 &  0.0101  & 0.0040  \\[0.1cm]
	\midrule
	\fullWENO{} &  0.0392 & 0.0140 & 0.0090 \\[0.1cm]
	\bottomrule
\end{tabular*}
\hspace*{0.35cm}
\begin{tabular*}{0.49\textwidth}{C{2.5cm}ccc}
	\toprule
	\begin{tabular}{@{}c@{}} \mbox{}\\\mbox{}\end{tabular} $TV_\uncertainty$ & $\MEElements =3$ & $\MEElements=10$ & $\MEElements =30$\\[0.2cm]
	\midrule
	analytic  &  0.1664 &  0.0503 & 0.0167\\[0.1cm]
	\midrule
	sG&  0.2422 & 0.0688 & 0.0200 \\[0.1cm]
	\midrule
	\schemeshort{} &  0.0572&  0.0188&  0.0086  \\[0.1cm]
	\midrule
	\fullWENO{} &  0.0186 &  0.0023 & 0.0008 \\[0.1cm]
	\bottomrule
\end{tabular*}
\begin{tabular*}{0.47\textwidth}{C{2.2cm}ccc}
	\multicolumn{4}{c}{} \\[0.2cm]
	\toprule
	\begin{tabular}{@{}c@{}} \mbox{}\\\mbox{}\end{tabular} $TV_\x$ & $\MEElements =3$ & $\MEElements=10$ & $\MEElements =30$\\[0.2cm]
	\midrule
	analytic &  1.0 &  1.0 & 1.0 \\[0.1cm]
	\midrule
	sG& 1.3118 & 1.2791 & 1.0795\\[0.1cm]
	\midrule
	\schemeshort{}&  1.0037 &  1.0044 &  1.0080\\[0.1cm]
	\midrule
	\fullWENO{} &  1.2316&  1.7983 & 1.3712\\[0.1cm]
	\bottomrule
\end{tabular*}
\hspace*{0.35cm}
\begin{tabular*}{0.49\textwidth}{cccc}
	\multicolumn{4}{c}{} \\[0.2cm]
	\toprule
	\begin{tabular}{@{}c@{}}percentage above \\ $TV_\x$ -- analytic\end{tabular} & $\MEElements =3$ & $\MEElements=10$ & $\MEElements=30$\\[0.2cm]
	\midrule
	analytic&--  &--   &--\\[0.1cm]
	\midrule
	sG& 23.8\%  & 21.8\% &  7.4\% \\[0.1cm]
	\midrule
	\schemeshort{}&  0.4\% & 0.4\% & 0.8\%\\[0.1cm]
	\midrule
	\fullWENO{} &  18.9\% & 44.4\% & 27.1\%\\[0.1cm]
	\bottomrule
\end{tabular*}
	\caption{$L_1$ error and total variation for linear advection with and without stochastic slope limiter (SL) for 2000 space cells, $K_\Omega = 2$ and $\WENOdegree=2$. Example \eqref{sweInitialConditions}.} 
	\label{table:LAerror}
\end{table}

\externaltikz{LA_surf}{\begin{figure}[htb]
\hspace*{-0.4cm}
\begin{subfigure}[c]{0.245\textwidth}
\begin{tikzpicture}
\begin{axis}[
width=.68\textwidth,
height=.68\textwidth,
scale only axis,
enlargelimits=false,
xlabel = $x$,
ylabel = $\xi$,
axis lines= box,
xmin=0.4,
xmax=2,
ymin=-0.33333,
ymax=0.3333,
ylabel style = {rotate=-90},
ylabel style={at={(axis description cs:-0.25,0.6)},anchor=north},
]
\addplot graphics [xmin=0.4,ymin=-0.33333,xmax=2,ymax=0.3333] {Images/LA3_analyt.png};
\end{axis}
\end{tikzpicture}
\subcaption{analytic}
\label{fig:surfLAa}
\end{subfigure}
\begin{subfigure}[c]{0.245\textwidth}
\begin{tikzpicture}
\begin{axis}[
width=.68\textwidth,
height=.68\textwidth,
scale only axis,
enlargelimits=false,
xlabel = $x$,
ylabel = $\xi$,
axis lines= box,
xmin=0.4,
xmax=2,
ymin=-0.33333,
ymax=0.3333,
ylabel style = {rotate=-90},
ylabel style={at={(axis description cs:-0.25,0.6)},anchor=north},
]
\addplot graphics [xmin=0.4,ymin=-0.33333,xmax=2,ymax=0.3333] {Images/LA3_SL0.png};
\end{axis}
\end{tikzpicture}%
\subcaption{sG}
\end{subfigure}
\begin{subfigure}[c]{0.245\textwidth}
\begin{tikzpicture}
\begin{axis}[
width=.68\textwidth,
height=.68\textwidth,
scale only axis,
enlargelimits=false,
xlabel = $x$,
ylabel = $\xi$,
axis lines= box,
xmin=0.4,
xmax=2,
ymin=-0.33333,
ymax=0.3333,
ylabel style = {rotate=-90},
ylabel style={at={(axis description cs:-0.25,0.6)},anchor=north},
]
\addplot graphics [xmin=0.4,ymin=-0.33333,xmax=2,ymax=0.3333] {Images/LA3_SLM.png};
\end{axis}
\end{tikzpicture}%
\subcaption{\schemeshort{}}
\label{fig:surfLAc}
\end{subfigure}
\begin{subfigure}[c]{0.245\textwidth}
\begin{tikzpicture}
\begin{axis}[
axis on top,
width=.68\textwidth,
height=.68\textwidth,
scale only axis,
enlargelimits=false,
xlabel = $x$,
ylabel = $\xi$,
axis lines= box,
xmin=0.4,
xmax=2,
ymin=-0.33333,
ymax=0.3333,
ylabel style = {rotate=-90},
ylabel style={at={(axis description cs:-0.25,0.6)},anchor=north},
colorbar,
point meta min=-0.4913,
point meta max=1.4897,
colorbar style={
	width= .06\textwidth,
	height = .68\textwidth,
	at={(1.15,1)},
	anchor=north west,
	ytick = {0,1},
	yticklabel style={font=\footnotesize,anchor=south, /pgf/number format/.cd, fixed,precision = 4,/tikz/.cd},
	yticklabel shift = 5pt,
}
]
\addplot graphics [xmin=0.4,ymin=-0.33333,xmax=2,ymax=0.3333] {Images/LA3_2DWENO.png};
\end{axis}
\end{tikzpicture}
\subcaption{2D WENO}
\label{fig:surfLAd}
\end{subfigure}

 \caption{Comparison of the analytical solution, sG, \schemeshort{} and 2D WENO (from left to right) with 3 \ME{}s in the $\x$ - $\uncertainty$ plane for the linear advection problem at $t=0.5$, using 2000 space cells, $\SGtruncorder=2$, $\WENOdegree=2$. Example \eqref{sweInitialConditions}.}
\label{fig:LAsurf}
\end{figure}}

\subsubsection{Burgers' Equation}
In this numerical example we consider the Burgers' equation
\begin{equation}
 \frac{\partial}{\partial t} \solOned +  \frac{\partial}{\partial x} \Big(\frac{\solOned^2}{2}\Big)=0.
\end{equation}
We compare the \schemeshort{} scheme to standard stochastic Galerkin and the \fullWENO{} method from \algoref{algo:fullWENO} as we did for the linear advection problem in the previous subsection. For the nonlinear Burgers' equation, we consider the following initial state
\begin{equation}\label{eq:BurgerOscil}\solution(0,\x,\uncertainty) = \sin(2\pi(\x + 0.1\uncertainty)),\end{equation}
where $\x\in\Domain=[0,1]$ and the uncertainty is uniformly distributed, i.e. $\uncertainty\sim\mathcal{U}(-1,1)$. Then we compute the solution on a $\ncells=2000$ space grid until $t=0.4$ and derive the $L_1$ error to a reference solution obtained by Monte Carlo sampling at the quadrature nodes in $\uncertainty$ and the total variations \eqref{eq:TVx} and \eqref{eq:TVxi}. 

The results can be found in \tabref{table:Burgererror}, where we now observe similar values for \schemeshort{} and \fullWENO{}. The $L_1$ errors are again slightly higher for these two modified methods compared to plain stochastic Galerkin. They all converge while increasing the number of \ME{}s. However, the total variation with respect to $\x$ is almost the same as in the reference solution, while standard stochastic Galerkin yield an increase of around 14\%. Thus, the overshoots could be completely eliminated which can be verified in \figref{fig:Burgererror}, showing the solution for 10 \ME{}s in the $\x-\uncertainty$ plane. Here, the sG approach in the left picture has huge oscillations that vanish in the \schemeshort{} approximation presented  in the right picture. The total variation wrt. $\uncertainty$ shows a similar behavior as for the linear advection example in \tabref{table:LAerror} and \figref{fig:LAsurf}, hence, it is reduced by the modified schemes due to the lack of information that is transported along the uncertainty. We have used $\nbRnodes=1000$ quadrature nodes in $\uncertainty$  and $\nbxnodes=4$ quadrature nodes in $\x$ for the calculations of the total variations. 
 
 \begin{table}[htbp]
 \centering
\begin{tabular*}{0.47\textwidth}{C{2.2cm}ccc}
	\toprule
	\begin{tabular}{@{}c@{}} \mbox{}\\\mbox{}\end{tabular} $L_1$ error & $\MEElements =3$ & $\MEElements=10$ & $\MEElements =30$\\[0.2cm]
	\midrule
	reference &--&--&--\\[0.1cm]
	\midrule 
	sG&  0.0227 & 0.0064 & 0.0018 \\[0.1cm]
	\midrule
	\schemeshort{}&  0.0283 &  0.0085 &  0.0028 \\[0.1cm]
	\midrule
	\fullWENO{} &  0.0284 & 0.0143  &  0.0053 \\[0.1cm]
	\bottomrule
\end{tabular*}
\hspace*{0.35cm}
\begin{tabular*}{0.49\textwidth}{C{2.5cm}ccc}
	\toprule
	\begin{tabular}{@{}c@{}} \mbox{}\\\mbox{}\end{tabular} $TV_\uncertainty$ & $\MEElements =3$ & $\MEElements=10$ & $\MEElements =30$\\[0.2cm]
	\midrule
	reference  &  0.2257 &  0.0677 & 0.0226\\[0.1cm]
	\midrule
	sG&  0.3187 & 0.0949 & 0.0306 \\[0.1cm]
	\midrule
	\schemeshort{} &  0.1132 &  0.0339 &  0.0112  \\[0.1cm]
	\midrule
	\fullWENO{} &  0.1134 &  0.0102 & 0.0012 \\[0.1cm]
	\bottomrule
\end{tabular*}
\begin{tabular*}{0.47\textwidth}{C{2.2cm}ccc}
 		\multicolumn{4}{c}{} \\[0.2cm]
 		\toprule
 		\begin{tabular}{@{}c@{}} \mbox{}\\\mbox{}\end{tabular} $TV_\x$ & $\MEElements =3$ & $\MEElements=10$ & $\MEElements =30$\\[0.2cm]
 		\midrule
 		reference &  3.3846 &  3.3846 & 3.3846 \\[0.1cm]
 		\midrule
 		sG&  3.9526 & 3.9476 & 3.9491\\[0.1cm]
 		\midrule
 		\schemeshort{}& 3.3817 &  3.3836 &  3.3867\\[0.1cm]
 		\midrule
 		\fullWENO{} &  3.4391 &  3.3897 & 3.3899\\[0.1cm]
 		\bottomrule
\end{tabular*}
\hspace*{0.35cm}
\begin{tabular*}{0.49\textwidth}{cccc}
 		\multicolumn{4}{c}{} \\[0.2cm]
 		\toprule
 		\begin{tabular}{@{}c@{}}percentage above \\ $TV_\x$ -- reference\end{tabular} & $\MEElements =3$ & $\MEElements=10$ & $\MEElements=30$\\[0.2cm]
 		\midrule
 		reference &--  &--   &--\\[0.1cm]
 		\midrule
 		sG& 14.3\%  & 14.2\% &  14.2\%\\[0.1cm]
 		\midrule
 		\schemeshort{}&  0\% & 0\% & 0\%\\[0.1cm]
 		\midrule
 		\fullWENO{} &  1.6\% & 0.1\% & 0.1\% \\[0.1cm]
 		\bottomrule
\end{tabular*}
	\caption{$L_1$ error and total variation for Burgers' equation with and without stochastic slope limiter (SL) for 2000 space cells, $K_\Omega = 2$ and $\WENOdegree=2$. Example \eqref{eq:BurgerOscil}.} 
	\label{table:Burgererror}
\end{table}
\externaltikz{Burger_surf}{
\begin{figure}[htb]
\begin{tikzpicture}

\begin{axis}[
width=0.35\textwidth,
height=0.35\textwidth,
scale only axis,
xlabel = $x$,
ylabel = $\xi$,
zlabel = $\solOned $,
axis lines= box,
xmin = 0.000000,
xmax = 1.000000,
ymin = -0.040000,
ymax = 0.040000,
zmin = -1.600000,
zmax = 1.600000,
scaled ticks=false, 
zticklabel style={/pgf/number format/fixed},
xlabel style={at={(axis description cs:0.1,0)},anchor=north},
ylabel style = {rotate=-90},
ylabel style={at={(axis description cs:0.61,-0.05)},anchor=north},
zlabel style = {rotate=-90},
xtick = {0, 0.5, 1},
ytick = {-0.04, 0.04},
yticklabels = {-0.1, 0.1},
ztick = {-1, 0, 1}
]

\addplot3 graphics[points={%
(0.694829,-0.004900,-1.002008) => (347.283461,136.246869)
(0.317099,-0.009475,-0.032754) => (259.480800,205.762894)
(0.950222,0.021241,-0.174124) => (371.381630,226.254238)
(0.034446,0.023616,0.468202) => (137.905593,260.351407)
}] {Images/Burger_SL0.png};

\end{axis}
\end{tikzpicture}
\hspace*{0.5cm}
\begin{tikzpicture}
\begin{axis}[
width=0.35\textwidth,
height=0.35\textwidth,
scale only axis,
xlabel = $x$,
ylabel = $\xi$,
zlabel = $\solOned$,
axis lines= box,
xmin = 0.000000,
xmax = 1.000000,
ymin = -0.040000,
ymax = 0.040000,
zmin = -1.600000,
zmax = 1.600000,
scaled ticks=false, 
zticklabel style={/pgf/number format/fixed},
xlabel style={at={(axis description cs:0.1,0)},anchor=north},
ylabel style = {rotate=-90},
ylabel style={at={(axis description cs:0.61,-0.05)},anchor=north},
zlabel style = {rotate=-90},
xtick = {0, 0.5, 1},
ytick = {-0.04, 0.04},
yticklabels = {-0.1, 0.1}
]

\addplot3 graphics[points={%
(0.678735,0.012438,-0.713846) => (316.713591,169.988169)
(0.757740,-0.026305,-1.452252) => (395.825243,87.696349)
(0.743132,0.016484,-1.289178) => (326.685963,125.868007)
(0.392227,-0.037453,1.035065) => (321.146890,280.153278)
}]  {Images/Burger_SLM.png};
\end{axis}

\end{tikzpicture}%

 \caption{Comparison of the sG approach (left) for the Burgers' equation with 10 \ME{}s, 2000 space cells, $\SGtruncorder=2$, $\WENOdegree=2$ and the \schemeshort{} approximation (right) at $t=0.4$. Example \eqref{eq:BurgerOscil}.}
\label{fig:Burgererror}
\end{figure}
}

\subsubsection{Compressible Euler Equations}
The one-dimensional compressible Euler equations for the flow of an ideal gas are given by \eqref{eulereq}.
We consider the Euler Equations with an uncertain shock test case as in \cite{Poette2009}, which is given by the following initial conditions:
\begin{subequations} \label{eq:initialEulerSod}
\begin{align}
\density(0,\x,\uncertainty) &= \begin{cases}
1, \qquad &\x < 0.5 + 0.05\uncertainty, \hspace*{2cm}\\
0.125, \qquad &\x\geq 0.5 + 0.05\uncertainty,
\end{cases}\\
\momentum(0,\x,\uncertainty) &= 0,\\
\energy(0,\x,\uncertainty) &= \begin{cases}
0.25, \qquad &\x < 0.5 + 0.05\uncertainty,\\
2.5, \qquad &\x\geq 0.5 + 0.05\uncertainty,
\end{cases}
\end{align}
\end{subequations}
where $\uncertainty\sim\mathcal{U}(-1,1)$ and $\x\in\Domain=[0,1]$. Moreover we choose $\eulergamma=1.4$ and compute the solution on a $\ncells=2000$ grid until $t=0.1$.  For this test case, we again compare the plain stochastic Galerkin approach to \schemeshort{} and \fullWENO{} described in \algoref{algo:hSGWENO} and \algoref{algo:fullWENO}, respectively. Note that we have to apply the hyperbolicity-preserving limiter from \cite{Schlachter2017a} in order to ensure the hyperbolicity of the underlying system. \refone{Otherwise, the computation would already crash in the first time step. This has already been observed for this specific test case in \cite{Poette2009}. Our numerical implementations show that the application of the stochastic slope limiter is not sufficient to preserve hyperbolicity of the solution. Therefore, the additional usage of the hyperbolicity-preserving limiter is inevitable for systems of conservation laws. } The $L_1$ error and total variations \eqref{eq:TVx} and \eqref{eq:TVxi} for the density $\density$ are shown in \tabref{table:Eulererror}. 

As before, we observe the smallest $L_1$ error for standard stochastic Galerkin and that each method is converging if we increase the number of \ME{}s. Only this time, the \fullWENO{} method has a better error than \schemeshort{}. The results for the total variation in $\uncertainty$ are similar to the previous scalar examples. In this test case, the modified methods only marginally improved the total variation with respect to $\x$, \tabref{table:Eulererror} shows almost the same percentage for each of the schemes and even an increase for \fullWENO{} and 10 \ME{}s. This is additionally demonstrated in \figref{fig:Eulererror}, where the plain stochastic Galerkin approach and the \fullWENO{} method are illustrated for 10 \ME{}s in the $\x-\uncertainty$ plane.  They only show minor oscillations compared to the previous examples, however, the overshoots at the boundaries of the \ME{} in sG could be eliminated using the full WENO reconstruction. We have used $\nbRnodes=1000$ quadrature nodes in $\uncertainty$  and $\nbxnodes=4$ quadrature nodes in $\x$ for the calculations of the total variation. 

 \begin{table}[htbp]
\begin{tabular*}{0.47\textwidth}{C{2.2cm}ccc}
	\toprule
	\begin{tabular}{@{}c@{}} \mbox{}\\\mbox{}\end{tabular} $L_1$ error & $\MEElements =3$ & $\MEElements=10$ & $\MEElements =30$\\[0.2cm]
	\midrule
	reference &--&--&--\\[0.1cm]
	\midrule 
	sG&  0.0040 & 0.79e-03 & 1.23e-04\\[0.1cm]
	\midrule
	\schemeshort{}&  0.0092 &  2.75e-03 &  8.06e-04 \\[0.1cm]
	\midrule
	\fullWENO{} &  0.0062 &  2.51e-03& 8.17e-04 \\[0.1cm]
	\bottomrule
\end{tabular*}
\hspace*{0.35cm}
\begin{tabular*}{0.49\textwidth}{C{2.5cm}ccc}
		\toprule
		\begin{tabular}{@{}c@{}} \mbox{}\\\mbox{}\end{tabular} $TV_\uncertainty$ & $\MEElements =3$ & $\MEElements=10$ & $\MEElements =30$\\[0.2cm]
		\midrule
		reference  &  0.0291 &  0.0068 & 0.0028\\[0.1cm]
		\midrule 
		sG&  0.0322 & 0.0092 & 0.0030 \\[0.1cm]
		\midrule
		\schemeshort{} &  0.0033 &  0.0008 &  0.0002  \\[0.1cm]
		\midrule
		\fullWENO{} &  0.0179 &  0.0018 & 0.0002 \\[0.1cm]
		\bottomrule
\end{tabular*}
\begin{tabular*}{0.47\textwidth}{C{2.2cm}ccc}
		\multicolumn{4}{c}{} \\[0.2cm]
		\toprule
		\begin{tabular}{@{}c@{}} \mbox{}\\\mbox{}\end{tabular} $TV_\x$ & $\MEElements =3$ & $\MEElements=10$ & $\MEElements =30$\\[0.2cm]
		\midrule
		reference &  0.8781 &  0.8781 & 0.8781 \\[0.1cm]
		\midrule
		sG&  1.0622 &  0.9407 & 0.9124 \\[0.1cm]
		\midrule
		\schemeshort{}& 1.0501 &  0.9225 &  0.8912\\[0.1cm]
		\midrule
		\fullWENO{} &  1.0134 &  1.0218 & 0.9068\\[0.1cm]
		\bottomrule
\end{tabular*}
\hspace*{0.35cm}
\begin{tabular*}{0.49\textwidth}{cccc}
		\multicolumn{4}{c}{} \\[0.2cm]
		\toprule
		\begin{tabular}{@{}c@{}}percentage above \\ $TV_\x$ -- reference\end{tabular} & $\MEElements =3$ & $\MEElements=10$ & $\MEElements=30$\\[0.2cm]
		\midrule
		reference &--  &--   &--\\[0.1cm]
		\midrule
		sG& 17.3\%  & 6.6\% &  3.5\%\\[0.1cm]
		\midrule
		\schemeshort{}&  16.4\% & 4.5\% & 1.3\%\\[0.1cm]
		\midrule
		\fullWENO{} &  13.3\% & 13.9\% & 2.4\% \\[0.1cm]
		\bottomrule
\end{tabular*}
	\caption{$L_1$ error and total variation of density $\density$ for Euler equations with and without stochastic slope limiter (SL) for 2000 space cells, $K_\Omega = 2$ and $\WENOdegree=2$. Example \eqref{eq:initialEulerSod}.} 
	\label{table:Eulererror}
\end{table}

\externaltikz{Euler_surf}{
\begin{figure}[htb]
\begin{tikzpicture}

\begin{axis}[
width=0.35\textwidth,
height=0.35\textwidth,
scale only axis,
xlabel = $x$,
ylabel = $\xi$,
zlabel = $\density $,
axis lines= box,
xmin = 0.000000,
xmax = 1.000000,
ymin = -0.100000,
ymax = 0.100000,
zmin = 0.000000,
zmax = 1.100000,
scaled ticks=false, 
zticklabel style={/pgf/number format/fixed},
xlabel style={at={(axis description cs:0.35,-0.05)},anchor=north},
ylabel style = {rotate=-90},
ylabel style={at={(axis description cs:0.9,0)},anchor=north},
zlabel style = {rotate=-90},
xtick = {0, 0.5, 1},
ytick = {-0.1, 0.1},
ztick = {0, 1}
]

\addplot3 graphics[points={%
	(0.814724,0.026472,1.053258) => (386.078769,339.168318)
	(0.905792,-0.080492,1.061377) => (343.467113,314.153703)
	(0.126987,-0.044300,0.173374) => (170.153081,125.028783)
	(0.913376,0.009376,1.067652) => (400.375130,336.434250)
}] {Images/euler_sl0.png};

\end{axis}
\end{tikzpicture}
\hspace*{0.5cm}
\begin{tikzpicture}
\begin{axis}[
width=0.35\textwidth,
height=0.35\textwidth,
scale only axis,
xlabel = $x$,
ylabel = $\xi$,
zlabel = $\density$,
axis lines= box,
xmin = 0.000000,
xmax = 1.000000,
ymin = -0.100000,
ymax = 0.100000,
zmin = 0.000000,
zmax = 1.100000,
scaled ticks=false, 
zticklabel style={/pgf/number format/fixed},
xlabel style={at={(axis description cs:0.35,-0.05)},anchor=north},
ylabel style = {rotate=-90},
ylabel style={at={(axis description cs:0.9,0)},anchor=north},
zlabel style = {rotate=-90},
xtick = {0, 0.5, 1},
ytick = {-0.1, 0.1},
ztick = {0, 1}
]

\addplot3 graphics[points={%
	(0.814724,0.026472,1.053258) => (386.078769,339.168318)
	(0.905792,-0.080492,1.061377) => (343.467113,314.153703)
	(0.126987,-0.044300,0.173374) => (170.153081,125.028783)
	(0.913376,0.009376,1.067652) => (400.375130,336.434250)
}] {Images/euler_weno2d.png};
\end{axis}

\end{tikzpicture}%

 \caption{Comparison of the sG approach (left) for the density $\density$ in the Euler equations with 10 \ME{}s, 2000 space cells, $\SGtruncorder=2$, $\WENOdegree=2$ and the full \fullWENO{} approximation (right) at $t=0.1$. Example \eqref{eq:initialEulerSod}.}
\label{fig:Eulererror}
\end{figure}
}

\section{Conclusions and Outlook}
In this article, we demonstrated the propagation of Gibbs phenomenon into the stochastic domain of the stochastic Galerkin system based on an uncertain linear advection example. This lead to the formulation of our modified stochastic Galerkin scheme, including the stochastic slope limiter, which is supposed to reduce overshoots in the  solution manifold. The scheme is combined with a Multielement ansatz, a WENO finite volume method and Runge Kutta time stepping, giving the so called WENO stochastic Galerkin scheme and altogether a stable high order approximation of the solution of the conservation law. Combined with the hyperbolicity-preserving limiter from \cite{Schlachter2017a}, the method is able to be used on any hyperbolic system of equations. We additionally considered a similar numerical scheme using full \fullWENO{} reconstruction in both the physical and stochastic domain, which is motivated by \cite{Tokareva2014}.

We applied the two methods to the scalar linear advection and nonlinear Burgers' equation as well as to the system of Euler Equations and compared the results to the standard stochastic Galerkin scheme using a WENO reconstruction in the physical space.
An analysis of the total variations within the physical and stochastic domain verified the reduction of Gibbs oscillations for example up to 23\% for the linear advection problem compared to plain stochastic Galerkin, coming with the price of an slightly higher $L_1$ error. This underlines the necessity of our stochastic slope limiter for discontinuous solutions in the $x-\xi$ plane. The \schemeshort{} and \fullWENO{} methods behave differently on our numerical test cases which indicates variable choices on the investigated problem setting.

\refone{Another approach to dampen oscillations is to filter the gPC coefficients in the stochastic Galerkin expansion, which is presented in \cite{Kusch2018}. The method is so far applied  in combination with first-order numerical schemes and it would be interesting to be combined with WENO reconstructions as proposed in this article in order to obtain high-order approximations and to further reduce the impact of Gibbs phenomenon. 
Moreover, our stochastic slope limiter can be constructed as a maximum principle satisfying limiter, which is a strategy known from the theory of deterministic conservation laws and explained in \cite{Schneider2015,Zhang2011}, where it is, similarly to our hyperbolicity-limiter, combined with a positive-preserving limiter in order to ensure the maximum principle as well as admissible solutions which would yield comparable properties to IPM \cite{Poette2009,Kusch2017}.}

\reftwo{Due to the curse of dimensionality in the stochastic Galerkin system, this scheme is mainly applied to low dimensional random variables, where it is able to outperform non-intrusive methods such as Multi-Level Monte Carlo \cite{Mishra2014} and stochastic collocation \cite{Xiu2005}.} In this context, it is important to find a range of applicability for the general stochastic Galerkin method, see for example \cite{Meyer2019a}, \reftwo{where we can see that stochastic Galerkin is outperformed by stochastic collocation already in two or three dimensions.}
\reftwo{Thus, the method presented in this article is - similar to general intrusive uncertainty quantification methods - dedicated to systems of conservation laws in low dimensions. For a modification of the stochastic slope limiter to multi-dimensional uncertainties we refer to the approaches in deterministic conservation laws \cite{Cockburn1990,cockburn1997runge}, whereas other techniques such as the sub-cell limiter with JST troubled cell indicator in \cite{Sonntag2017} might be relevant to adopt into the framework of capturing shocks in the stochastic variable.
In addition to that, for multi-dimensional stochastic or spatial variables, we have to extend the WENO reconstruction to more than two dimensions. However, ENO or WENO techniques are mainly used for low dimensions and we require the implementation of sparse grids as for example in \cite{Kolb2018}. }

\section*{Acknowledgements}
Funding by the Deutsche Forschungsgemeinschaft (DFG) within the RTG GrK 1932 ``Stochastic Models for Innovations in the Engineering Science'' is gratefully acknowledged.

\section*{References}
\bibliographystyle{siam}
\bibliography{library,bibliography}

\end{document}